\documentclass[dvips]{conm-p-l}  %% FOR LATEX
\usepackage%[pdftex]
     {graphicx}
\usepackage{amssymb,amstext}

\newcommand{\texorpdfstring}[2]{#1}
\newcommand{\secLab}[1]{\label{section:#1}}
\newcommand{\thmLab}[1]{\label{thm:#1}}
\newcommand{\propLab}[1]{\label{prop:#1}}
\newcommand{\lemLab}[1]{\label{lemma:#1}}
\newcommand{\corLab}[1]{\label{cor:#1}}
\newcommand{\conjLab}[1]{\label{conj:#1}}
\newcommand{\defLab}[1]{\label{def:#1}}
\newcommand{\eqLab}[1]{\label{eq:#1}}
\newcommand{\figLab}[1]{\label{fig:#1}}
\newcommand{\secRef}[1]{\ref{section:#1}}
\newcommand{\thmRef}[1]{\ref{thm:#1}}
\newcommand{\propRef}[1]{\ref{prop:#1}}
\newcommand{\lemRef}[1]{\ref{lemma:#1}}
\newcommand{\corRef}[1]{\ref{cor:#1}}
\newcommand{\defRef}[1]{\ref{def:#1}}
\newcommand{\eqRef}[1]{\thetag{\ref{eq:#1}}}
\newcommand{\figRef}[1]{\ref{fig:#1}}
\def\SingleFig#1#2{%
    \begin{figure}[htb]
     \centering
     \includegraphics{Images/#1}
    \caption{#2\figLab{#1}}
    \end{figure}
    }

\newtheorem{theorem}{Theorem}[section]
\newtheorem{lemma}[theorem]{Lemma}
\newtheorem{prop}[theorem]{Proposition}
\newtheorem{cor}[theorem]{Corollary}

\theoremstyle{definition}
\newtheorem{definition}[theorem]{Definition}

\theoremstyle{remark}
\newtheorem{remark}[theorem]{Remark}
\newtheorem{conjecture}[theorem]{Conjecture}
\newtheorem{algorithm}[theorem]{Algorithm}

\usepackage{pict2e}

\newcommand{\findheight}[1]{\setbox0=\hbox{$#1X$}\unitlength=0,05\ht0 }
\newcommand{\mathPalett}[1]{\mathchoice
  {\linethickness{.45pt}\findheight\displaystyle      #1{4}}
  {\linethickness{.45pt}\findheight\textstyle         #1{4}}
  {\linethickness{.4 pt}\findheight\scriptstyle       #1{5.5}}
  {\linethickness{.35pt}\findheight\scriptscriptstyle #1{6}}}

\newcommand{\nnonpointed}{{\mathPalett\makeX}}
\newcommand{\npointed}{{\mathPalett\makeY}}
\newcommand{\makeX}[1]{% argument 1 is size of the dot.
    \begin{picture}(15,20)
      \put(0,0){\line(3,4){15}}
      \put(7.5,10){\line(3,-5){6}}
      \put(7.5,10){\line(-3,6){5}}
      \put(7.5,10){\circle*{#1}}
    \end{picture}}
\newcommand{\makeY}[1]{%
    \begin{picture}(16,20)(-2.5,0)
      \linethickness{.4pt}
     \put(7.5,10){\line(3,5){6}}
      \put(7.5,10){\line(-3,-3){10}}
      \put(7.5,10){\line(-3,5){6}}
      \put(7.5,10){\circle*{#1}}
    \end{picture}}
\def\nrVert{n}
\def\nrEdges{m}

\def\nrCH{n_B}
\def\nrInt{n_I}
\def\nrReflex{r}
\def\nrCorners{k}

\def\nrNonPointed{\nrVert_\nnonpointed}

\def\pts{P}
\def\points{\{ p_1, \ldots, p_{\nrVert}\}}
\def\GP{G(\pts)}

\newcommand{\PT}{{\textit{PT}}}
\newcommand\reals{{\mathbb{R}}}
\newcommand\eps{\varepsilon}

\newcommand\psT{pseudo-trian\-gu\-la\-tion}
\newcommand\PsT{Pseudo-Trian\-gulation}
\newcommand\ppsT{pointed pseudo-trian\-gu\-la\-tion}
\newcommand\PpsT{Pointed Pseudo-Trian\-gulation}
\newcommand\rnp{relative nonpointed}
\newcommand\rp{relative pointed}

\newcommand{\complaint}[1]{}

\def\withcomplaints{  %  \input{Latex/complaint}
\newcounter{mycomplaints}
\def\complaint##1{\stepcounter{mycomplaints}%
\ifhmode\unskip
{\dimen1=\baselineskip \divide\dimen1 by 2 %
\raise\dimen1\llap{\tiny -\themycomplaints-}}\fi
\marginpar{\sloppy\hbadness 10000 \tiny [\themycomplaints]: ##1}}%
}

\begin{document}

\title{Pseudo-Triangulations --- a Survey}

\author{G{\"{u}}nter ~Rote}
\address{Institut f\"ur Informatik, Freie Universit\"at
Berlin, Takustra{s}e%
~9, D-14195~Berlin, Germany.
}
\email{rote@inf.fu-berlin.de}
\thanks{First author partly supported by the Deutsche
Forschungsgemeinschaft (DFG) under grant RO~2338/2-1.}

\author{Francisco ~Santos}
\address{Departamento de Matem\'aticas, Estad\'{\i}stica y
Computaci\'on, Universidad de Canta\-bria, E-39005 Santander, Spain}
\email{francisco.santos@unican.es}
\thanks{Second author supported by grant MTM2005-08618-C02-02
of the Spanish Ministry of  Education and Science.}

\author{Ileana ~Streinu}
\address{Department of Computer Science, Smith College,
Northampton, MA 01063, USA.} \email{streinu@cs.smith.edu}
\thanks{Third author supported by NSF grants CCR-0430990 and NSF-DARPA CARGO-0310661.}

%    General info
\date{\today}
%\date{September 24, 2006 and, in revised form, October 15, 2006.}

%\complaint{To agree on AMS classification, if OK: Done.
%Here is a selection: 05C10 Topological
%graph theory, imbedding; 05C62 Graph representations (geometric and
%intersection representations, etc.); 52B11 $n$-dimensional
%polytopes; 52C25 Rigidity and flexibility of structures; 68U05
%Computer graphics; computational geometry; }


%-------------------------------------------------------------------------------         Abstract           --------

\begin{abstract}
A pseudo-triangle is a simple polygon with exactly three convex vertices,
and a pseudo-triangu\-lation is a face-to-face tiling of a planar region into
pseudo-triangles. Pseudo-triangulations appear as data structures in
computational geometry, as planar bar-and-joint frameworks in
rigidity theory and as projections of locally convex surfaces. This
survey of current literature includes combinatorial properties and
counting of special classes, rigidity theoretical results,
representations as polytopes, straight-line drawings from abstract
versions called combinatorial pseudo-triangu\-lations, algorithms
and applications of pseudo-triangulations.
\end{abstract}

\maketitle

\setcounter{tocdepth}{3} % for draft version
\setcounter{tocdepth}{2} % for draft version
\setcounter{tocdepth}{1} % for arxiv version !!
\tableofcontents

%-------------------------------------------------------------------------------         Main        --------

%!TEX root = article.tex

%-------------------------------------------------------------------------------         Introduction        --------

\section{Introduction}
\secLab{introduction}

A \emph{pseudo-triangle} is a simple polygon in the plane  with
exactly three convex vertices, called \emph{corners}.  A
\emph{\psT} is a tiling of a planar region into
pseudo-triangles. In particular, a triangle is a pseudo-triangle and
\psT s are generalizations of triangulations.
\SingleFig{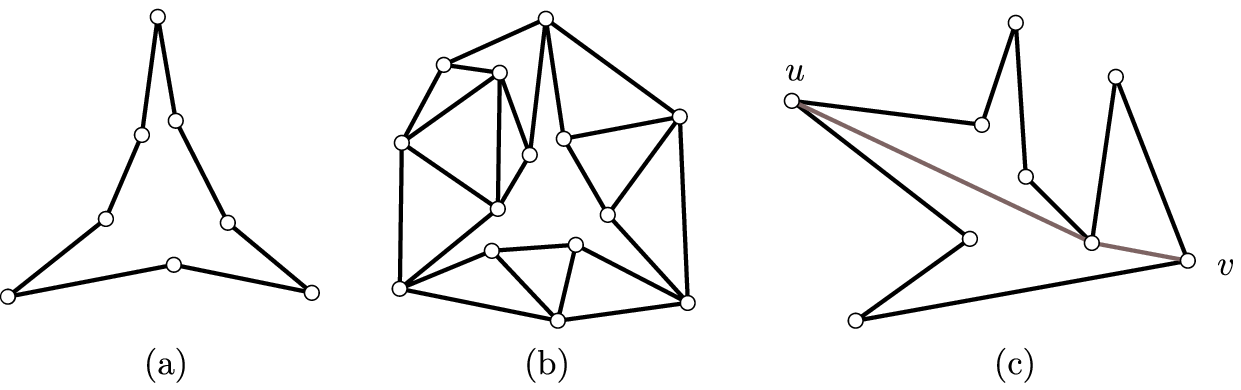}{(a) A pseudo-triangle, (b) a
pseudo-triangu\-lation of a point set and (c) a
pseudo-triangu\-lation of a simple polygon, including
 a geodesic path from $u$ to $v$.}
Special cases include the \emph{\psT\  of a finite
point set} and that of a \emph{simple polygon}, which partition the
convex hull of the point set, resp.~the interior of the polygon,
into pseudo-triangles and use no additional vertices. See
Figure~\figRef{pseudoTriang}.

Pseudo-triangulations have arisen in the last decade as interesting
geometric-combinatorial objects with connections and applications in
visibility, rigidity theory and motion planning.

\subsection*{Historical Perspective}
The names \emph{pseudo-triang\-le} and
\emph{pseudo-triangu\-la\-tion} were coined by Pocchiola and Vegter
in 1993~\cite{pocchiola:vegter:visibilityComplex:1996,pocchiola:vegter:pt:96}, inspired by a connection with pseudoline arrangements
\cite{pocchiola:vegter:orderType:94}. They were studying the
\emph{visibility complex} of a set of disjoint convex obstacles in the
plane~\cite{pocchiola:vegter:topoSweep:1996,pocchiola:vegter:visibilityComplex:1996},
and defined pseudo-triangulations by taking a maximum number of
non-crossing and free bitangents to pairs of objects. 
We review their work in Sections~\secRef{obstacles},~\secRef{visibility-complex}, and \secRef{pseudo-lines}.

For polygons, pseudo-triangulations
had already appeared in the computational geometry literature in the
early 1990's, under the name of \emph{geodesic
triangulations}~\cite{chazelle:etAl:rayShootingGeodesicTriangulations:1994,
goodrich:tamassia:dynamicRayShooting:1997}, and were obtained by
tiling a polygon via non-crossing geodesic paths joining two polygon
vertices, as in Figure~\figRef{pseudoTriang}(c). Compactness and
ease of maintenance led to their use as efficient \emph{kinetic
data structures} for collision detection of
polygonal obstacles~%
\cite{agarwal:basch:guibas:hershberger:zhang:kineticCollisionDetection:2003,
beghz-kcdts-04,
%basch:erickson:guibas:hershberger:zhang:kineticCollisionDetection:1999,
speckmann:kirkpatrick:kineticMaintenance:2002,speckmann:kirkpatrick:snoeyink:kineticCollisionDetection:2002}.
See Section~\secRef{geodesic}.

In 2000, the work of Streinu
\cite{streinu:pseudoTriangRigidityMotionPlanning-confAndJour:2005}
on the Carpenter's Rule Problem 
(see Section~\secRef{carpenter}) brought in an entirely different
perspective from rigidity theory. She showed that \emph{pointed}
pseudo-triangu\-lations, when viewed as \emph{bar-and-joint
frameworks} (or \emph{linkages with fixed edge-lengths}) are
minimally rigid, and become expansive mechanisms with the removal of
a convex hull edge.
(A \psT\ is pointed if every vertex is incident
to an angle larger than~$\pi$, see Section~\secRef{basicProperties}
for more definitions.)
%\complaint{GR. inserted early def. otherwise the reader can make no sense}
 Expansive motions were a crucial ingredient in
the solution to the Carpenter's Rule Problem by Connelly, Demaine
and Rote earlier that year
\cite{connelly:demaine:rote:carpenterRule:2003}. This newly
discovered combinatorial expression was  further exploited in
\cite{streinu:pseudoTriangRigidityMotionPlanning-confAndJour:2005}
for a second, pseudo-triangu\-lation-based, algorithmic solution of
the same problem.

These results not only hinted for the first time to the deep
connections between pseudo-triangulations and rigidity theory, but
also highlighted their nice combinatorial properties and emphasized
the importance of the concept of \emph{pointedness}. They also led
to the use of pseudo-triangulations in the investigation of the cone
of all expansive infinitesimal motions of a point
set~\cite{rote:santos:streinu:polytopePseudoTriangulations:2003},
which resulted in the definition of the polytope of pointed
pseudo-triangu\-lations. This appears as a natural generalization of
the well-studied associahedron~\cite{lee:associahedron:1989}, which
corresponds to triangulations of a convex point set in the plane and
thus, indirectly, to a long list of other combinatorial objects with
ubiquitous applications in computer science and combinatorics 
(\emph{Catalan structures} such as binary trees, lattice paths, stacks,
etc.).

This work triggered several lines of research on
pseudo-triangulations in the last five years. Here are most of those
we are aware of:
\begin{itemize}
\item Two more {\em polytopes of pseudo-triangulations} have
been found: one is a direct generalization of the polytope
from~\cite{rote:santos:streinu:polytopePseudoTriangulations:2003}
but covers all (not necessarily pointed)
pseudo-triangulations~\cite{orden:santos:polytope:2005}; the other
is more an analogue of the {\em secondary polytope} of
triangulations (see~\cite{billera:filliman:sturmfels:constructionsSecondaryPolytopes:1990,
gelfand:kapranov:zelevinsky:discriminantsResultants:1994}), and
stems from the work of Aichholzer, Aurenhammer, Krasser, and Bra\ss\
relating pseudo-triangulations
to locally convex
functions~\cite{aichholzer:aurenhammer:brass:krasser:pseudoTriangulationsNovelFlip:2003}.
\item There has been an increased interest in the study of combinatorial properties of
pseudo-triangulations: their number, vertex degrees, and how these compare for
different point sets or with respect to the analogous concepts in triangulations~%
\cite{aichholzer:hoffmann:speckmann:toth:degreeBoundsConstrainedPseudoTriangulations:2003,
aichholzer:orden:santos:speckmann:numberPseudoTriangulationsCertainPointSets:2006+,
aichholzer:aurenhammer:krasser:speckmann:convexityMinimizes:2004,
speckmann:kettner:etAl:tightDegreeBounds:2003,
santos:randall:rote:snoeyink:countingTriangulationsWheels:2001,
rote:wang:wang:xu:constrainedPseudoTriangulations:2003}.

\item Related to this, but with algorithmic applications in mind,
the \emph{diameter} of the graphs of
flips~\cite{aichholzer:aurenhammer:krasser:adaptingPseudoTriangulations:2003,
aahk-tstpt-06, bereg:transforming-pseudo-triangulations:2004}, and
methods for the \emph{efficient enumeration} of
pseudo-triangulations
\cite{streinu:aichholzer:rote:speckmann:zigZagPath-wads:2003,bereg:enumerating-pseudo-triangulations:2005,
bronnimann:kettner:pocchiola:snoeyink:CountingAndEnumeratingPointedPseudoTriang:2006}
have been studied.

\item A stronger connection between planar graphs and
pseudo-triangu\-lations came with the proof that not only are
pseudo-triangu\-lations rigid (and pointed pseudo-triangulations
minimally rigid), but the converse is also true: every planar
(minimally) rigid graph admits a drawing as a (pointed)
pseudo-triangu\-lation
\cite{streinu:haas:etAl:planarMinRigidPseudoTriang-confAndJour:2005,
orden:santos:servatiusB:servatiusH:combinatorialPT:2006}. To prove this result, the concept of
\emph{combinatorial pseudo-triangu\-lations} is introduced. They are defined as
plane maps in which each internal face has three specified \emph{corners}. 

\item One of the key tools used in
\cite{connelly:demaine:rote:carpenterRule:2003,
streinu:pseudoTriangRigidityMotionPlanning-confAndJour:2005} for the
Carpenter's Rule Problem was Maxwell's Theorem from 1864, relating
projections of polyhedral surfaces to plane self-stressed frameworks
and to the existence of reciprocal diagrams. In the same spirit is
the work of Aichholzer et al.~\cite{aichholzer:aurenhammer:brass:krasser:pseudoTriangulationsNovelFlip:2003},
where a special type of locally convex piecewise-linear surface is
related, via projections, to pseudo-triangu\-lations of polygonal
domains.
Maxwell's reciprocal diagrams of (necessarily non-pointed)
pseudo-triangu\-lations are also considered in
\cite{orden:roteEtAl:nonCrossingFrameworks:2004}.

\item As a further connection with rigidity theory, Streinu's study of
  pointed   pseudo-triangulations~\cite{streinu:pseudoTriangRigidityMotionPlanning-confAndJour:2005}
  was extended to \emph{spherical \psT s}, with
  applications to the spherical Carpenter's Rule Problem and
  single-vertex origami%
  ~\cite{streinu:whiteley:origami:2005}.  This paper also contains
  partial work on combinatorial descriptions of expansive motions in
  three dimensions.
\item In the theory of rigidity with fixed edge-directions
(rather than fixed edge-lengths), pointed pseudo-triangu\-lation
mechanisms have been shown
to have a kinetic
behavior, linearly morphing tilings while remaining
non-crossing and pointed~\cite{streinu:parallelPseudoTriangulations:2005}.

\item Finally, pseudo-triangulations have found applications
as a tool for proofs:
in the area of art galleries
(illumination by
floodlights)~\cite{speckmann:toth:vertexGuardsPseudoTriang:2003}; and
in an area which,
at first sight, may seem unrelated to discrete geometry:
the construction of counter-examples to a conjecture of A. D. Alexandrov
characterizing the sphere among all smooth
surfaces~\cite{p-bcah-05}.

\end{itemize}

\subsection*{Overview}  This survey presents several points of view on
pseudo-triangu\-lations.
%, consistent with the ways in which they are
%used and analyzed in the literature.
First, as a tiling of a planar
region, they are related to each other by local changes called
flips. This is in several ways analogous to the ubiquitous
triangulations which appear almost everywhere in Combinatorial
Geometry, and has led to the investigation of similar questions:
counting, enumeration, flip types, connectivity and diameter. We
cover these topics in  Sections~\secRef{basicProperties} and~\secRef{setAll}. Next, in~Section~\secRef{liftings}, we study their
relationship with projections of locally convex surfaces in space.
This serves as a bridge between the combinatorial and the rigidity
properties of pseudo-triangulations, when viewed as bar-and-joint
frameworks, which are presented in Sections~\secRef{self-stresses-0},~\secRef{rigidity} and \secRef{planarLaman}.
Section~\secRef{polytope} describes polytopes of
pseudo-triangulations, whose construction relies on properties
studied in the preceding Sections~\secRef{liftings} and~\secRef{rigidity}. Finally, in Section~\secRef{applications} we
briefly sketch several applications of pseudo-triangu\-la\-tions
that have appeared in the literature, a preview of which appears
above in the historical introduction (ray shooting, visibility
complexes, kinetic data structures, and the Carpenter's Rule problem).

The emphasis of this survey is on concepts and on the logical flow
of ideas, and not so much on proofs or on the historical
developments. But we sometimes have found, and included, shorter proofs than
those in the literature.
In particular, in Section~\secRef{liftings} we provide
for the first time  a uniform treatment for lifted surfaces in
connection with \psT s, which appeared independently in the context
of the locally convex functions in
\cite{aichholzer:aurenhammer:brass:krasser:pseudoTriangulationsNovelFlip:2003},
and in the rigidity investigations of
\cite{streinu:pseudoTriangRigidityMotionPlanning-confAndJour:2005}.
The results in Sections~\secRef{number-k-gon} and~\secRef{graph-drawing} are published here for the first time.

%!TEX root = article.tex

%-------------------------------------------------------------------------------         Basic Properties        --------

\section{Basic Properties of Pseudo-Triangulations}
\secLab{basicProperties}

Pseudo-triangu\-lations generalize and inherit certain properties
from triangu\-lations. This section and the next one address their
similarities in a comparative manner. 

In this section, after fixing the basic terminology and notation to be used throughout, we
 exhibit the
simple relationships that exist among several parameters of a
pseudo-triangu\-lation: numbers of vertices, edges, faces, pointed
vertices, convex hull or outer boundary vertices. They lie at the
heart of the more advanced combinatorial properties presented later.

\subsection{Definitions}
A \emph{graph} $G=(V,E)$ has $\nrVert$ vertices, $V= %[\nrVert]:=
\{1,\ldots, \nrVert\}$ and $|E|=\nrEdges$ edges. A \emph{geometric graph}
is a drawing of $G$ in the plane with straight-line edges. The
mapping $V\to \reals^2$ of the vertices $V$ to a set of points $\pts
= \points$ is referred to as the (straight-line) \emph{drawing}, \emph{embedding} or \emph{realization} of $G$, and denoted $\GP$. With few exceptions, we will
consider realizations on point sets with distinct elements (which
induce edges of non-zero length), and in general position
(which permits the analysis of \emph{pointed} graph embeddings,
defined below).

\subsubsection*{Plane graphs} A geometric graph $G$
is \emph{non-crossing} or \emph{a plane graph} if two disjoint edges $ij,
kl\in E$, ($i,j \notin \{k,l\}$) are realized as disjoint (closed)
line segments. The complement of the points and edges is a
collection of planar regions called faces, one of which is
unbounded. When $G$ is connected, the bounded faces are
topologically disks, and the unbounded face is a disk with a hole. A
graph is \emph{planar} if it admits a plane embedding. With few
exceptions, the graphs considered in this paper are planar and
connected.

\subsubsection*{Polygons and corners}
A simple polygon is a non-crossing embedding of a cycle. It
partitions the plane into two connected regions:
an interior region and an exterior,
unbounded one. More generally, we may encounter degenerate disk-like
open polygonal \emph{regions}, which have non-simple (self-touching
but non-crossing) polygonal boundaries, called \emph{contours}, as
in Fig.~\figRef{degenerate}~(right) on p.~\pageref{fig:degenerate}.

A vertex of a polygonal region (simple or degenerate) is called
\emph{convex}, \emph{straight} or \emph{reflex} depending on whether
the angle spanned by its two incident edges, facing the polygonal
region, is strictly smaller, equal to or strictly larger than $\pi$,
respectively. General position for the vertices, which we usually
assume, implies the absence of straight angles. Convex vertices
incident to a face are also called \emph{corners} of that face.

\subsubsection*{Pseudo-$\nrCorners$-gons, pseudo-triangles%, interior bitangent.
}
A simple polygon with exactly $\nrCorners$ corners is called a
pseudo-$\nrCorners$-gon. The special cases $\nrCorners=3$ and
$\nrCorners=4$ are called \emph{pseudo-triangles} and %, resp.\
\emph{pseudo-quadrilaterals}. A bounded face in a geometric
graph must have at least $3$
corners, but the unbounded face may be a pseudo-$k$-gon with $k\leq
2$. %%
\subsubsection*{Point sets, polygons and pointgons} We will work with
point sets (denoted $\pts$), polygons (denoted $R$) and what we call
\emph{pointgons}. A \emph{pointgon} $(R,\pts)$ is a polygon $R$
and a specified finite set of points $\pts$, including
the vertices of $R$ and (perhaps) additional points in
its interior. A polygon $P$ is a special case of a pointgon, with no
interior points. Similarly, a point set $\pts$ can naturally be
considered a pointgon, in which $R$ is the convex hull of $P$.

\subsubsection*{Pointed graph embeddings} A vertex of an embedded
graph is called \emph{pointed} if some pair of consecutive edges (in
the cyclic order around the vertex) span an angle larger than~$\pi$, and \emph{non-pointed} otherwise. Here we are making use of the general position assumption.
 %
 %\complaint{Added parenthetical remark for refere 1. Paco, September 30, 2007}
 %
  The two edges
incident to the reflex angle are called the \emph{extreme} edges of
the pointed vertex. A pointed (planar) graph embedding is one with
all its vertices pointed.

\subsubsection*{Pseudo-triangulations} A pseudo-triangu\-lation is a
planar embedded connected graph whose bounded faces are
pseudo-triangles. The following three variants have been considered
in the literature, depending on whether the boundary is allowed to
be non-convex or whether interior points are allowed as vertices:
\begin{itemize}
\item A \emph{pseudo-triangu\-lation of a simple polygon $R$} is a
subdivision of the interior of $R$ into pseudo-triangles, using only
the polygon vertices.
\item A \emph{pseudo-triangu\-lation of a pointgon} $(R,\pts)$ partitions the
interior of the polygon $R$ into pseudo-triangles using as vertices
all of the points $\pts$.
\item A \emph{pseudo-triangu\-lation of a finite point set} $\pts$ is a
\psT\ of the pointgon $(R,\pts)$, where $R$ is the convex hull of
$\pts$. In particular, a triangulation of $\pts$ (using all
vertices) is a pseudo-triangu\-lation.
\end{itemize}

%%%%%%
\SingleFig{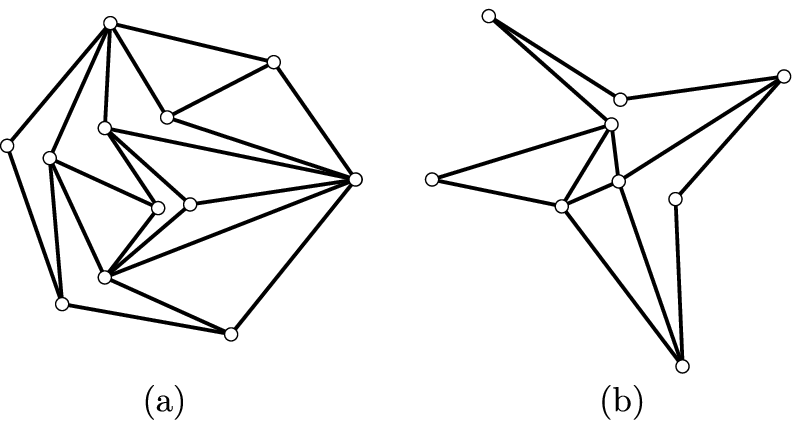}{(a) a non-pointed pseudo-triangu\-lation of a point
set, (b) a non-pointed pseudo-triangu\-lation of a pointgon.}
%%%%%%

In any of the variants, a \emph{pointed pseudo-triangu\-lation} is
one in which every vertex is pointed. See
Figure~\figRef{pseudoTriang}(b,c) for examples of pointed
pseudo-triangu\-lations, and Figure~\figRef{ppt} for
non-pointed ones.

\subsection{Pseudo-Triangulations of Polygons}
\secLab{polygons}
Throughout most of this paper we deal with pseudo-triangulations of pointgons, because they are the most general case, or with pseudo-triangulations of point sets, because they have nicer properties, specially related to rigidity. But let us start with the study of pseudo-triangulations of simple polygons as an introduction to the subject.

The \emph{geodesic path} between two points (typically, but not
necessarily, two vertices) of a polygon $R$ is the shortest path
from one to the other \emph{in} $R$ (with $R$ understood as  a
region).
Pseudo-triangu\-lations of a simple polygon $R$ are also
called \emph{geodesic triangulations}, because they arise by
inserting non-crossing  geodesic paths in $R$.
For example, the pseudo-triangulation of
Figure~\figRef{pseudoTriang}c is obtained by inserting a single geodesic, 
between the corners $u$ and $v$.
A special case of such a geodesic triangulation
is the \emph{shortest path tree}, consisting of the geodesic paths 
from any chosen vertex to all other vertices.
%, or from a corner to 
%all other corners if one wants a pointed pseudo-triangulation.
%vertices,
%or only to all corners if one wants a pointed pseudo-triangulation.
Geodesic triangulations are historically the first pseudo-triangulations considered in the literature,
introduced in~\cite{chazelle:etAl:rayShootingGeodesicTriangulations:1994,goodrich:tamassia:dynamicRayShooting:1997} for ray shooting and shortest path queries in changing structures.
See~\secRef{geodesic} for more details on these applications.

We will focus on geodesics between \emph{corners}.
If the corners are consecutive (in the cyclic order of all corners around $R$),
then the geodesic path between them is a sequence of polygon edges called a \emph{pseudo-edge}.
If they are not consecutive, the geodesic path consists  of (perhaps zero) polygon edges and 
(at least one) interior diagonals. Here, a \emph{diagonal} is
any segment joining two vertices of $R$ through the interior of
$R$. For a diagonal $e$ to be part of some geodesic between corners 
it needs to have the special property that the graph $R\cup e$ is
pointed. Such diagonals are called (interior)~\emph{bitangents}\secLab{tangent} of $R$,
because they are \emph{tangent} at both endpoints. 
We say that a line segment $l$ ending in a vertex $p$ of
a pseudo-triangle $\Delta$ is \emph{tangent} to
$\Delta$ at $p$ if either $p$ is a corner of $\Delta$ and $l$ lies in 
the convex angle at $p$, or $p$ is a reflex vertex
and the two incident edges
% the reflex angle 
at~$p$ lie on the same side of the supporting line of~$l$.
%both % sides
%directions of 
%the line through $l$ lie in the reflex angle at~$p$.
\complaint{Defined tangent and bitangent. Paco, october 14.
Modified after Guenter. Remodified GR since ``both sides'' of a line
is misleading.
Modified again Ileana, "directions" is also misleadings.}%
The geodesic path of Figure~\figRef{pseudoTriang}c  consists 
of two bitangents.

A bitangent may be part of several geodesics, but there is always a
canonical one:

\begin{lemma}
\lemLab{tangents} For every bitangent of $R$ there is a unique pair of
corners such that the geodesic between them consists only of this
bitangent plus \textup(possibly\textup) some
boundary edges of $R$.
\end{lemma}

\begin{proof}
Extend the bitangent by following at each end, in the direction of
tangency, a (possibly empty) sequence of edges into the next corner.
\end{proof}

There is an important analogy between pointed pseudo-triangulations of a pseudo-$k$-gon and triangulations of a convex $k$-gon.
Every pseudo-$k$-gon with $k>3$ has at least one bitangent (since the geodesic path between any two non-consecutive corners includes bitangents) and the insertion of a bitangent in a pseudo-$k$-gon divides it into a pseudo-$i$-gon and a pseudo-$j$-gon with $i+j=k+2$. With the same arguments as for triangulations of a convex $k$-gon it is easy to conclude that
by recursively inserting bitangents in a $k$-gon $R$ one eventually gets a pointed pseudo-triangulation of it. When we say ``recursively'' we imply that the edges included are bitangents  not (only) of the original polygon $R$, but of the polygon they are inserted in. This ensures that the graphs obtained are always pointed.

By the splitting formula $i+j=k+2$ we need  $k-3$ bitangents and get $k-2$ pseudo-triangles 
to arrive to a pointed pseudo-triangulation. Moreover, all pseudo-triangulations of $R$ arise in this way:

\begin{theorem}
\thmLab{geodesic}
Every pointed pseudo-triangulation of a pseudo-$\nrCorners$-gon consists of $\nrCorners-2$ pseudo-triangles and
uses $\nrCorners-3$ interior bitangents of $R$.
\qed
\end{theorem}

The maximum and minimum number of pointed pseudo-triangulations that a pseudo-$k$-gon can have are obtained in Section~\secRef{number-k-gon}.
A case of special interest is a pseudo-quadrilateral:

\begin{lemma}
\lemLab{quadrilateral}
A pseudo-quadrilateral has exactly two interior bitangents. Hence, it also has two pointed pseudo-triangulations, obtained by inserting one or the other bitangent.
\end{lemma}

\begin{proof}
That there are at least two bitangents follows from the fact that the two geodesics between opposite corners must use different sets of bitangents. That there are only two follows from Lemma~\lemRef{tangents}: given any bitangent, the geodesic of that lemma must be one of the two geodesics between opposite corners.
\end{proof}

As in the case of triangulations, this property of pseudo-quadrilaterals is at the heart of the concept of \emph{flip} that will be introduced in Section~\secRef{flips}.

\subsection{Vertex and Face Counts}
\secLab{count}
One of the basic properties of
triangulations in the plane
%, which follows easily from Euler's Formula, 
is that all triangulations of the same region
and with the same set of vertices have the same number of edges (and
of faces). The following theorem generalizes this to pseudo-triangulations.

\begin{theorem}
\thmLab{count}
Let $(R,\pts)$ be a  pointgon on $|\pts|=\nrVert$ points and
$\nrReflex$ reflex vertices in the polygon $R$. Let $T$ be a
pseudo-triangu\-lation of $(R,\pts)$ with $\nrNonPointed$ non-pointed
vertices. Then $T$ has $2\nrVert- 3 + (\nrNonPointed - \nrReflex)$
edges and $\nrVert-2 + (\nrNonPointed - \nrReflex)$
pseudo-triangles.
\end{theorem}

\begin{proof}
Let $\nrEdges$ denote the number of edges. Then $2\nrEdges$ equals
the total number of angles in $T$, since a  vertex of degree $d$ is
incident with $d$ angles. Now we count separately the number of
convex and reflex angles. The reflex angles are $\nrVert -
\nrNonPointed$ (one at each pointed vertex). The convex angles are
the three in each pseudo-triangle plus the exterior angle of each
reflex vertex of $R$. Hence,
\[
2 \nrEdges = 3 t + \nrReflex+ \nrVert - \nrNonPointed,
\]
where $t$ is the number of pseudo-triangles.
By Euler's Formula, $\nrEdges + 1 = \nrVert + t$.
Eliminating $t$ (respectively $\nrEdges$) from these two formulas
gives the statement.
\end{proof}

%\medskip
%\noindent
 Here are some interesting special cases:

%\noindent
{\bf Triangulations.}
In this case the only pointed vertices are the convex vertices of $R$,
so that  $\nrNonPointed = \nrReflex + \nrInt
$, where $\nrInt
$
is the number of interior vertices.
This leads to the well-known
relation $|E|=2 \nrVert - 3 + \nrInt
$.

%\noindent
{\bf Pseudo-triangu\-lations of a point set.} The polygon $R$ has no
reflex vertices, so that $\nrReflex = 0$. The number of
non-pointed vertices can go from zero in the case of a pointed pseudo-triangulation
to the number  $\nrInt$
of points in the interior of $R$.
%In any case:

\begin{theorem}
\thmLab{count2}
Let $\pts$ be a point set with $\nrVert$ elements. Then, every
pseudo-triangulation $T$ of $\pts$ has $2\nrVert -3 +\nrNonPointed$
edges, where $\nrNonPointed$ is the number of non-pointed vertices in
it.
\qed
\end{theorem}

In particular, pointed pseudo-triangulations have the minimum
possible number of edges, namely $2\nrVert -3$,  among all
pseudo-triangulations of $\pts$. This motivated the term
\emph{minimum} pseudo-triangu\-lations
in~\cite{streinu:pseudoTriangRigidityMotionPlanning-confAndJour:2005},
for what are now called pointed pseudo-triangulations.

%\noindent
{\bf Geodesic Triangu\-lations.}
Geodesic triangulations have 
$\nrReflex=\nrVert-\nrCorners$, where $\nrCorners$ is the number of
corners. Theorem~\thmRef{count} gives $\nrVert + \nrCorners- 3 + \nrNonPointed$
edges and $\nrCorners-2 + \nrNonPointed$
pseudo-triangles. This reduces to Theorem~\thmRef{geodesic}
in the pointed case.

\subsection{Pointedness}
\secLab{pointedPT}
Pseudo-triangulations may be regarded as maximal non-crossing graphs
with a prescribed set of pointed vertices:
\begin{theorem}
\thmLab{maximal} A non-crossing geometric graph $T$  is a
pseudo-triangu\-lation of its underlying point set $\pts$ if and
only if its edge set is maximal among the non-crossing geometric
graphs with vertex set $\pts$ and with the same set of pointed
vertices as $T$.
\end{theorem}

The same result is true for pseudo-triangulations of a pointgon $(R,\pts)$, under the additional hypothesis that $T$ contains all the boundary edges of $R$.
\complaint{added this sentence. Paco, october 14; fixed a bit, Ileana oct 15}%

\begin{proof}
\emph{Only if:}  since $T$ is a pseudo-triangu\-lation of $\pts$,
any additional edge will go through the interior of a
pseudo-triangle. But pseudo-triangles have no interior bitangents, so this edge
 creates a non-pointed vertex.

\emph{If:} suppose that no edge can be inserted without making some
pointed vertex non-pointed. In particular, all convex hull edges
of $\pts$ are in $T$, since convex hull vertices cannot be made non-pointed
by the addition of any edge. We prove that every interior face $R$
is a pseudo-triangle. A priori, the face may not even be simply
connected if $T$ is not connected, but it will always have a
well-defined outer contour. Number the corners of $R$ along this
contour from $v_1$ to $v_{\nrCorners}$, $\nrCorners\geq 3$. Consider
two paths $\gamma_+$ and $\gamma_-$ from vertex $v_1$ to $v_3$
through the interior of $R$, close to the contour of $R$ and in
opposite directions. Now \emph{shorten} them continuously as much as possible,
i.\,e., consider geodesic paths $\gamma'_1$ and $\gamma'_2$ homotopic
to them.
Adding these paths will maintain pointedness at all pointed vertices.
The only possibility for these geodesic paths not to add any
edges to $T$ is that they coincide (hence $R$ is simply connected)
and go along the boundary of $R$ (hence $v_1$ and $v_3$ are
consecutive corners, and $R$ is a pseudo-triangle).

%A similar argument works for the outer face, if is not a convex polygon.
\complaint{Ileana oct 15, deleted last sentence, was repeated above}
\end{proof}

In particular, pointed pseudo-triangulations can be reinterpreted as the
maximal non-crossing \emph{and pointed} graphs,
in the same way as triangulations
are the maximal non-crossing graphs. This leads to the following list of
equivalent characterizations of pointed pseudo-triangulations. Another one, which
shows how to incrementally build pointed pseudo-triangulations adding one vertex at a time,
will appear in Theorem~\thmRef{henneberg}.

\begin{theorem}[Characterization of pointed pseudo-triangu\-lations
\cite{streinu:pseudoTriangRigidityMotionPlanning-confAndJour:2005}]
\thmLab{charactPPT}
Let $T$ be a graph embedded on a set $\pts$ of $\nrVert$ points. The
following properties are equivalent.
\begin{enumerate}
\item \label{minimum} %\emph{\bf (Minimum pseudo-triangu\-lation)}
$T$ is a pseudo-triangu\-lation of $\pts$ with the {\bf minimum}
possible number of edges.
\item \label{acypseudo} %\emph{\bf (Pointed pseudo-triangu\-lation)}
$T$ is a {\bf pointed} pseudo-triangu\-lation  of $\pts$.
\item \label{lamanP} %\emph{\bf ($2\nrVert-3$ pseudo-triangu\-lation)}
$T$ is a pseudo-triangu\-lation of $\pts$ with {\bf $\nrEdges =
2\nrVert-3$} edges \textup(equivalently, with $f=\nrVert-2$ faces\textup).
\item \label{lamanA} %\emph{\bf ($2\nrVert-3$ Planar and Pointed)}
$T$ is {\bf non-crossing, pointed} and has {\bf $2\nrVert-3$
edges}.
\item \label{maximal} %\emph{\bf (Maximal Planar and Pointed)}
$T$ is {\bf pointed, non-crossing}, and {\bf maximal} \textup(among the
pointed non-crossing graphs embedded on $\pts$\textup).
\end{enumerate}
\end{theorem}

\begin{proof}
The first three equivalences, and the implication
(2)$\Rightarrow$(4),  follow from Theorem~\thmRef{count}. The
equivalence (2)$\Leftrightarrow$(5) is Theorem~\thmRef{maximal}. The
implication (4)$\Rightarrow$(5) is a combination of both theorems:
every non-crossing pointed graph can be completed to a maximal one,
which is, by  Theorem~\thmRef{maximal}, a pointed
pseudo-triangulation and has, by Theorem~\thmRef{count}, $2n-3$
edges. If that was already the number of edges we started with, then
the original graph was already maximal.
\end{proof}

The above theorem is similar to the well-known characterization of
trees as graphs that are connected and minimal, cycle-free and maximal, 
or that have any two of the properties of being connected, cycle-free and
having $n-1$ edges.

{}From Conditions (4) and~(5), we get:
\begin{cor}
\corLab{max-pointed-noncrossing}
A non-crossing pointed graph on $n$ points contains at most $2n-3$
edges.
\qed
\end{cor}

\subsubsection*{Pseudo-triangulations and Laman graphs}
A graph $G$ is called a \emph{Laman graph} if it has % exactly 
$2\nrVert-3$
edges and % satisfies the following hereditary property: 
every subset
of $\nrVert' \ge 2$ vertices spans at most $2\nrVert'-3$ edges. By
Theorem~\thmRef{count2} the graph of every pointed
pseudo-triangulation satisfies the first property. By
Theorem~\thmRef{maximal} it also satisfies the second, since every
subgraph will itself be pointed. Hence:
\begin{cor}
[Streinu
\cite{streinu:pseudoTriangRigidityMotionPlanning-confAndJour:2005}]
\corLab{Laman}
The underlying graphs of pointed pseudo-trian\-gu\-lations of a point set are Laman graphs.
\end{cor}

Laman graphs are crucial in 2-dimensional rigidity theory: in almost all embeddings, they are
minimally rigid. Corollary~\corRef{Laman} points to the deep connections between
pseudo-triangu\-lations and rigidity, which will be developed in
Section~\secRef{rigidity}. Similar arguments applied to arbitrary pseudo-triangulations lead to:

\begin{cor}
\corLab{generalized-Laman}
Let $T$ be a pseudo-triangulation of $\nrVert$ points, $\nrNonPointed$ of them non-pointed. Then, $T$ has  $2\nrVert -3 + \nrNonPointed$ edges and  for every subset of $\nrVert'\ge 2$ vertices, $\nrNonPointed'$ of them non-pointed, the induced subgraph has at most $2\nrVert'-3 + \nrNonPointed'$ edges.
\end{cor}

Note that this \emph{generalized Laman property} is not a property of
the abstract graph, but of a geometric one, since we need to know
 which vertices are pointed.

\medskip

For a non-crossing geometric graph, let us define the \emph{excess
of corners} $\bar{\nrCorners}$ as the number of convex angles minus three times the
number of bounded faces. Since every
bounded face has at least three corners, the excess of corners is always 
at least zero, with equality if and only if the graph is a
pseudo-triangulation of a point set. The following statement is the
most general form of the formula for the number of edges of a
non-crossing graph in terms of pointedness.

\begin{theorem}
\thmLab{generalized-Laman} A connected geometric graph with
$\nrEdges $ edges, $\nrVert$ vertices, $\nrNonPointed$ of which are
non-pointed and with excess of corners $\bar{\nrCorners}$ satisfies\textup:
\[
\nrEdges = 2\nrVert -3 + (\nrNonPointed - \bar{\nrCorners}).
\]
\end{theorem}

The assumption of connectivity can be removed if $\bar{\nrCorners}$ is
defined additively on connected components, as the number of convex
angles minus three times the number of \emph{bounded face cycles}.

\begin{proof}
Use the same counts as in Theorem~\thmRef{count}. There are  $2\nrEdges$ angles,
 $\nrVert - \nrNonPointed$ of them reflex and  $3f + \bar{\nrCorners}$
 convex, where $f$ is the number of bounded faces.
Hence,
\[
2 \nrEdges = 3 f + \bar{\nrCorners} + \nrReflex+ \nrVert - \nrNonPointed.
\]
Euler's Formula $\nrEdges + 1 = \nrVert + f$ finishes the proof.
\end{proof}

In particular, the difference $\nrNonPointed - \bar{\nrCorners}$ does
not depend on the particular non-crossing embedding of a given
planar graph $G$. Later (Theorem~\thmRef{graph-drawing}
in Section~\secRef{graph-drawing}) we show that the rigidity properties of $G$ determine
how small the parameters $\nrNonPointed$ and $\bar{\nrCorners}$ can be.

%------------------------------------------------------------------------      Flips between pseudo-triangu\-lations    ------

\subsection{Flips in Pseudo-Triangulations}
\secLab{flips}

Flipping is the transformation of a \psT\ into another one
by inserting and/or removing an edge.
More precisely, we have the following types
of \emph{flips} in a pseudo-triangulation $T$ of a pointgon. See the examples in
Figure~\figRef{ptflips}.

\begin{itemize}
\item (\emph{Deletion, or edge-removal, flip})  The removal of an interior edge $e\in T$,
if the result is a pseudo-triangu\-lation.
%\complaint{AAKB call this ``edge-removing flip'',
%OS call it ``deletion flip''.
%Section 4 is still full of ``edge-removal'' flips! Should we switch
%to edge removal? GR}
%
\item (\emph{Insertion flip})  The insertion of a new edge $e\not\in
T$, if the result is a pseudo-triangu\-lation.

\item (\emph{Diagonal flip}) If $e$ is an interior edge whose removal does not
produce a pseudo-triangulation, then
there exists a unique edge $e'$
different from $e$ that can be added to obtain a new
pseudo-triangulation $T-e+e'$.
%: the other diagonal of the pseudo-quadrilateral
%$\Gamma$
%created in $T \setminus e$.
%
%
\end{itemize}
%%%%%%
\SingleFig{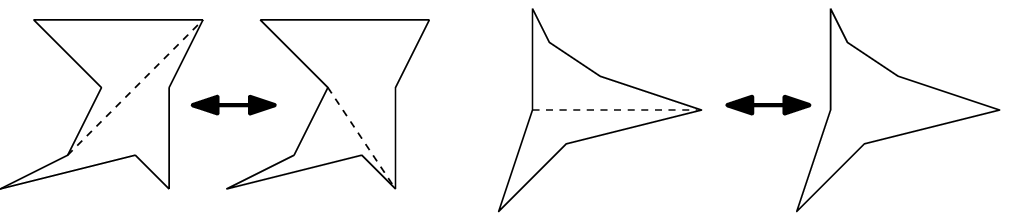}{Left, a diagonal flip. Right, an
insertion-deletion flip.
These pseudo-triangles might form part of a larger \psT.}
%%%%%%
%The following result will be generalized in Theorem~\thmRef{ptgraph}:

To see that these are the only possibilities,
let us analyze what happens when we remove an edge of a \psT.
\complaint{Ileana to Paco: please put the next sentences in the statement of the lemma 2.12 or 2.13.}%
Let $T$ be a pseudo-triangulation of a pointgon $(R,\pts)$ and let
$e$ be an interior edge, common to two pseudo-triangles. 
We have to prove that, if $T-e$ is not a \psT,
then there is a unique edge $e'\ne e$
such that $T-e+e'$ is again
a \psT.
%, as stated in the definition of diagonal flip.

When $e$ is removed from $T$, its two incident pseudo-triangles
become a single region $\Gamma$ that we can regard as a (perhaps
degenerate, see Figure~\figRef{degenerate}) polygon.

\begin{prop}
\thmLab{flips} This
\complaint{Ileana to Paco: style complaint; I'd prefer that you define the objects in the statement of the Lemma, not before it.}%
 polygon $\Gamma$ is:
\begin{itemize}
\item a pseudo-quadrilateral if both endpoints of $e$ preserve their pointedness with the removal;
this happens when $e$ is a bitangent of $\Gamma$.
\item a pseudo-triangle, otherwise. In this case, exactly one of the endpoints changes from non-pointed to
pointed with the removal.
\end{itemize}
\end{prop}

%%%%%%
\SingleFig{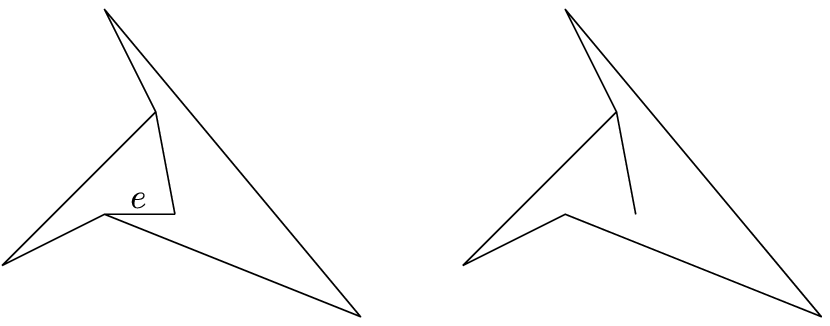}{The removal of an edge may produce a
degenerate pseudo-quadrilateral.}
%%%%%%
%

\begin{proof}
The statement can be proved geometrically, by looking at the old and new
angles at the endpoints of $e$.
Here we offer a
counting argument based on Theorem~\thmRef{generalized-Laman}.
Applied to $T$ and to $T\setminus e$, the theorem gives
\[
\nrEdges=2\nrVert - 3 + (\nrNonPointed- \bar{\nrCorners})
\]
and
\[
\nrEdges-1=2\nrVert - 3 + (\nrNonPointed' -\bar\nrCorners'),
\]
where $\nrEdges$, $\nrVert$, $\nrNonPointed$ and $\bar{\nrCorners}$ are the number of edges, vertices,
non-pointed vertices and
excess of corners in $T$, and $\nrNonPointed'$ and $\bar\nrCorners'$ are the same in $T\setminus e$. Hence,
\[
\nrNonPointed- \bar{\nrCorners} = \nrNonPointed'- \bar{\nrCorners}' -1.
\]

If both endpoints keep their (non-)pointedness, then the excess of
corners increases by one, which implies that $\Gamma$ is a
pseudo-quadrangle. If one endpoint passes from non-pointed to
pointed then the excess of corners is preserved, and $\Gamma$ is a
pseudo-triangle. It is impossible for both endpoints to pass from
non-pointed to pointed, since it would imply $\Gamma$ being a
pseudo-$2$-gon.
\end{proof}

Since every pseudo-quadrilateral has exactly two pointed
pseudo-triangulations (Lemma~\lemRef{quadrilateral}), we have:

\begin{prop}
\propLab{flip-count}
In every pseudo-triangulation $T$ of $(R,\pts)$ there is one flip
for each interior edge and one for each pointed vertex that is not a corner of $R$.
These are all possible flips.
\end{prop}

\begin{proof}
By Proposition~\thmRef{flips},
we have exactly one deletion or diagonal flip on every interior edge,
that deletes the edge and (if needed) inserts the other diagonal of the pseudo-quadrilateral formed.

An insertion flip is the inverse of a deletion flip and, by
Proposition~\thmRef{flips}, it turns a pointed vertex $p$ to non-pointed.
Moreover, as long as the reflex angle at $p$ is in a pseudo-triangle $\Delta$
(that is, if $p$ is not a corner of $R$),
there is one insertion flip possible at $p$, namely the insertion of the diagonal
that is part of the geodesic from $p$ to the
opposite corner of  $\Delta$.
\end{proof}

This statement suggests that we introduce a name for a
pointed vertex that is not a corner of $R$.
Equivalently, for a point that is a reflex vertex of some 
(and then a unique) region of $G$ contained in $R$.
We call them  \emph{pointed in $R$}, or \emph{\rp}.
The rest of vertices are called \emph{\rnp}.  This classification
was first introduced in~\cite{aichholzer:aurenhammer:brass:krasser:pseudoTriangulationsNovelFlip:2003}, were these two types of vertices were called
 ``incomplete'' (or ``pending'') and ``complete'', respectively.
%\complaint{Added this paragraph for referee 1. Edited the corresponding one in Section 4. Paco, September 30, 2007.}

%------------------------------------------------------------------------      Henneberg constructions     ------

\subsection{Henneberg Constructions of Pointed Pseudo-Triangulations} 
In pointed pseudo-triangulations, only
diagonal and insertion flips are possible. In particular, every
interior edge can be flipped, producing another pointed
pseudo-triangulation. This observation is the basis for an inductive procedure (called a \emph{Henneberg
construction}) for generating pseudo-triangulations of a point set.  It mimics the one devised by Henneberg \cite{henneberg:graphischeStatik:1911-68} for minimally rigid graphs (not necessarily planar).

%%%%%%
\SingleFig{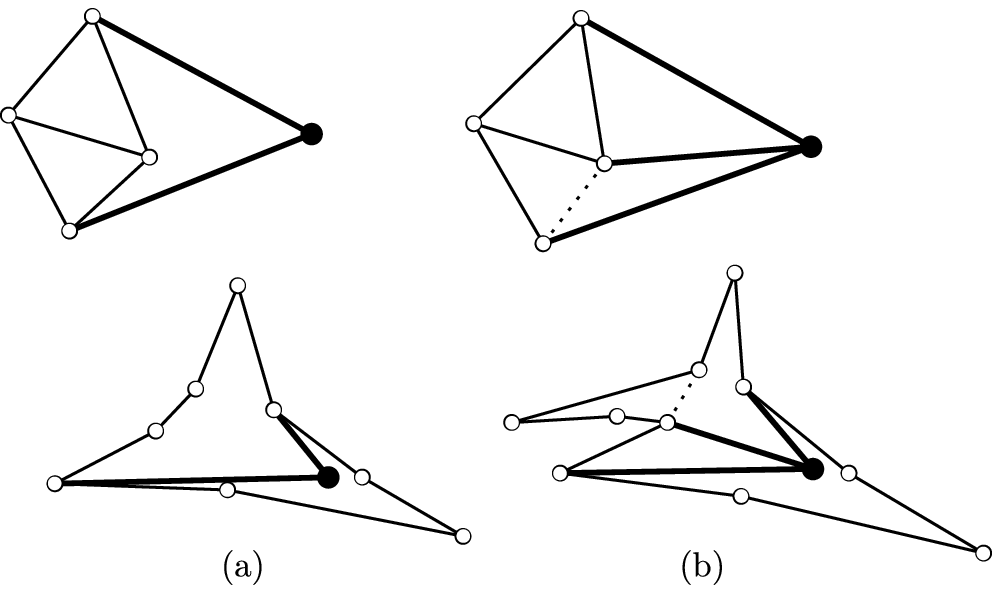}{Henneberg steps. (a) type 1 and (b) type 2.
Top row: the new vertex is added on the outside face. Bottom
row: it is added inside a pseudo-triangular face.
The added edges are thick. The dotted edge is the one that is removed
in the type-2 step.}
%%%%%

\begin{theorem}[Streinu~\cite{streinu:pseudoTriangRigidityMotionPlanning-confAndJour:2005}]
\thmLab{henneberg} Let $T$ be a pointed pseudo-triangulation of a
point set $\pts$. Then there exists an ordering
$p_1,p_2,\dots,p_\nrVert$ of the points $\pts$ and a sequence of
pointed pseudo-triangulations $T_i$, on the point set $\{p_1,\dots,p_i\}$, for
$i=3,\dots, \nrVert$, such that each $T_{i+1}$ is obtained from
$T_i$ in one of the following two ways
\textup(see Fig.~\figRef{henneberg}\textup)\textup:
\begin{enumerate}
\item\emph{Type 1 (add a vertex of degree $2$):} Join the vertex $p_{i+1}$ by two
segments. If $p_{i+1}$ is in the outer face of $T_i$ the segments are tangent to the boundary of
 $T_i$. Otherwise, the two segments are parts of geodesics to two of the three corners of the pseudo-triangle of $T_i$ containing $p_{i+1}$.
\item\emph{Type 2 (add a vertex of degree $3$):} Add the vertex $p_{i+1}$ with degree 2 as before, then flip an edge in the pseudo-edge opposite to $p_{i+1}$ in the unique triangle that has $p_{i+1}$ as a corner.
\end{enumerate}
\end{theorem}

\begin{proof}
Since $T$ has $2\nrVert -3$ edges, the average degree of a vertex is $4 - 6/\nrVert$.
In particular, there must be a vertex of degree two or three.

If there is a vertex of degree two, consider it the last vertex in the ordering, $p_\nrVert$.
Removing the two edges incident to it leaves a pointed non-crossing graph on $\nrVert-1$ vertices and with
$2(\nrVert-1)-3$ edges, hence a pseudo-triangulation that we call $T_{n-1}$.
Then, $T$ is obtained from $T_{n-1}$ by a type 1 step as described in the statement.

If there is no vertex of degree two, then there is a vertex of
degree three, that we take as $p_n$. Since $p_n$ is pointed, one of
its edges lies within the convex angle formed by the other two (and,
in particular, it is an interior edge). Let $T'$ be the
pseudo-triangulation obtained by flipping that edge, in which $p_n$
has degree two. Let $T_{n-1}$ be the pseudo-triangulation obtained
from $T'$ be removing the two edges incident to $p_n$, as before.
Then, $T$ is obtained from $T_{n-1}$ by a step of type 2.
\end{proof}

%!TEX root = article.tex

%------------------------------------------------------------------------      Set of all pseudo-triangu\-lations    ------

\section{The Set of all Pseudo-Triangu\-lations}
\secLab{setAll}

In this section we consider the set of {\em all} pseudo-triangulations of
a given point set or pointgon and look at it {as a
whole}. We address questions about their number, enumeration,
as well as pseudo-triangulations with extremal properties within the set.
Flips turn this \emph{set} into a {\em graph}.

%------------------------------------------------------------------------      Graph of pseudo-triangulations    ------

\subsection{The Graph of Pseudo-Triangulations}
\secLab{graph-of}

% Given $(R,\pts)$, 
The  \emph{graph of \psT s} of a pointgon $(R,\pts)$
has one node for each
pseudo-triangu\-lation of $(R,\pts)$ and an arc joining $T$ and
$T'$  if there is a flip producing one from the other. This is an undirected graph, 
since every flip has an inverse.

\begin{theorem}
\thmLab{ptgraph}
Let $(R,\pts)$ be a pointgon with $\nrVert$ vertices and $\nrInt$
interior points.
\begin{enumerate}
\item Its graph of pseudo-triangu\-lations is regular of degree $\nrVert+2\nrInt-3$.

\item The subgraph
induced by \ppsT s is also regular, of degree
$\nrVert-\nrReflex+\nrInt-3=\nrCorners+2\nrInt-3$, where $\nrReflex$
and $\nrCorners$ are the numbers of
reflex vertices and corners in $R$.
\end{enumerate}
In both cases the graph is connected.
\end{theorem}
\begin{proof}
By Proposition~\propRef{flip-count} the number of flips in a pseudo-triangulation
equals the sum of its interior edges plus its \rp\ vertices.
These two numbers are, respectively:
\[
2\nrVert - 3 + (\nrNonPointed-\nrReflex) - \nrCH
=
\nrVert  + \nrInt - 3 + (\nrNonPointed-\nrReflex)
\]
and
\[
\nrVert - \nrNonPointed - \nrCorners.
\]
Their sum is
\[
\nrVert  + \nrInt - 3 + (\nrNonPointed-\nrReflex)
+
\nrVert - \nrNonPointed - \nrCorners
=
2\nrVert + \nrInt - 3  - \nrReflex - \nrCorners
=
\nrVert + 2 \nrInt - 3 .
\]

This proves part (1). For part (2), deletion and insertion flips do not apply. Since we have
$\nrNonPointed=0$, we only need to count
interior edges. Hence the number of diagonal flips is:
\[
2\nrVert - 3 -\nrReflex - \nrCH =
\nrVert - 3 -\nrReflex + \nrInt.
\]

We prove connectivity only in the case of
pseudo-triangu\-lations of a point set $\pts$. For pointgons, a
proof can be found in
\cite{aichholzer:aurenhammer:brass:krasser:pseudoTriangulationsNovelFlip:2003}.
Let $p_i$ be  a point on the convex hull of $\pts$. The crucial
observation is that pseudo-triangu\-lations of $\pts
\setminus\{p_i\}$ (and flips between them) coincide with the
pseudo-triangu\-lations of $\pts$ that have  degree~$2$ at $p_i$, if
in the latter we forget the two tangents from $p_i$ to 
the convex hull of $\pts \setminus\{p_i\}$. 
%$\conv(\pts \setminus\{p_i\})$. 
By induction, we assume the
pseudo-triangu\-lations of $\pts \setminus\{p_i\}$ to be connected
in the graph. On the other hand, in pseudo-triangu\-lations with
degree greater than $2$ at $p_i$ all interior edges incident to $e$
can be flipped and produce pseudo-triangu\-lations with smaller
degree at $e$.
%\complaint{GR removed: Remark: this would not be
%necessarily true if $p_i$ was not a vertex of $\conv(\pts)$). }
%Decreasing one by one the number of edges incident to $p_i$ will
%eventually lead to a pseudo-triangu\-lation with degree $2$ at
%$p_i$.
\end{proof}

%%%%%%
\SingleFig{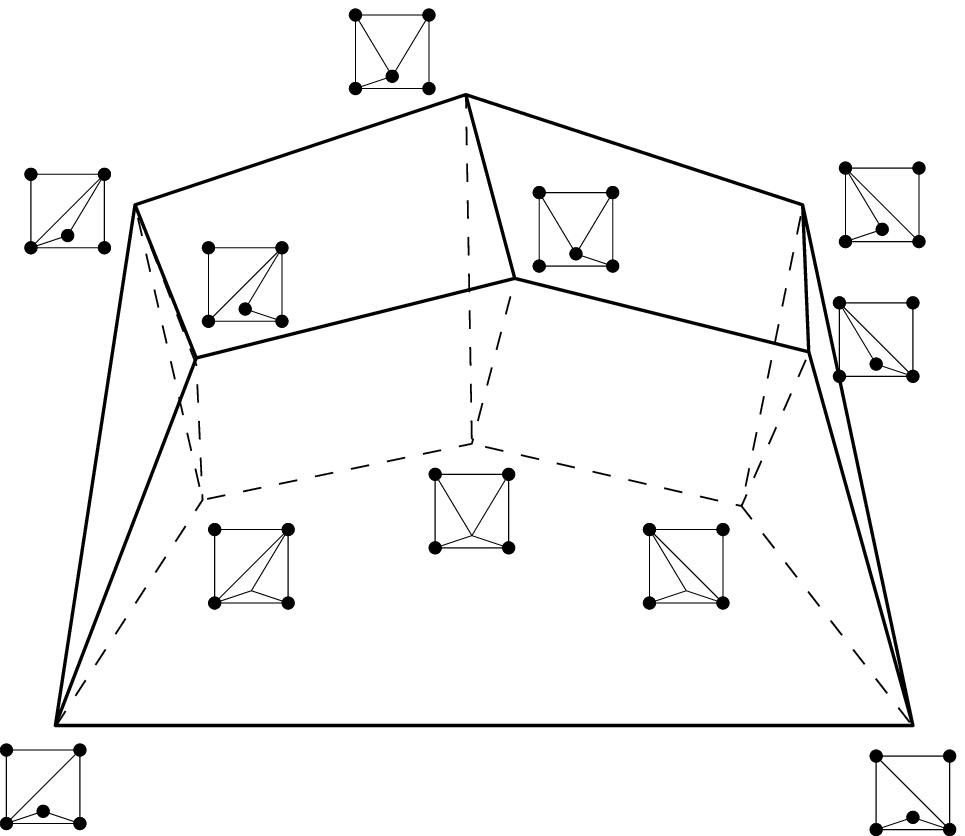}{The graph of all pseudo-triangu\-lations of
this point set, connected by flips, forms the 1-skeleton of a 4-polytope.
Pointed pseudo-triangulations form the 1-skeleton of a 3-polytope (solid lines).}%{0.8}
%%%%%%

As an example, Figure~\figRef{graphsq} shows the graph of pseudo-triangulations
of a set of five points, one of them interior. As predicted by
Theorem~\thmRef{ptgraph}, the whole graph is 4-regular and the graph of pointed pseudo-triangulations
(the solid edges) is 3-regular. In the picture, both the solid and the whole graph are 1-skeleta of simple
polytopes, of dimensions 3 and 4 respectively. We will prove in Section~\secRef{polytope}
(Theorems~\thmRef{expansion-polytope} and~\thmRef{pt-polytope}) that this is always the case 
for pseudo-triangulations of arbitrary pointgons. These polytopes generalize  the well-known \emph{associahedron}, whose 1-skeleton is the graph of
flips in triangulations of a convex polygon.

\subsubsection*{The diameter of the graph of \psT s}
The number of flips necessary to go from one pseudo-triangulation to another is called  the \emph{diameter} of the graph of pseudo-triangulations.

\begin{theorem}
\thmLab{diameter}
For every set $\pts$ of $\nrVert$ points:
\begin{enumerate}
\item 
The graph of all 
\psT s of  $\pts$
has diameter at most $O(n\log n)$
\emph{(Aichholzer et
al.~\cite{aichholzer:aurenhammer:krasser:adaptingPseudoTriangulations:2003,aichholzer:aurenhammer:brass:krasser:pseudoTriangulationsNovelFlip:2003})}.
\item
The subgraph induced by \ppsT s has diameter
at most $O(n\log n)$
 \emph{(Bereg~\cite{bereg:transforming-pseudo-triangulations:2004})}.
\end{enumerate}
\end{theorem}

Observe that part (2) does not follow from part (1) since the
distance between two pointed pseudo-triangulations can increase
when only pointed pseudo-triangulations are allowed as intermediate
steps. An example of this is given in~\cite{aahk-tstpt-06},
where it is also shown that the diameter bound in part (1) can be
refined to $O(n\log c)$ for a point set with $c$ convex layers.

\begin{proof}
  We present Bereg's divide-and-conquer proof, which was originally
  devised for part~(2) but works for both parts.

Sort the  points in clockwise order about
the left-most point $p_0$ and split them in half by a ray $l$ of
median slope (see Fig.~\figRef{bereg}). 
 Take the convex hulls of the two halves, with $p_0$ included in both
  halves.  This defines the \emph{median}
  pseudo-triangle from $p_0$,
 bounded by the common tangent to the two convex hulls
and the two  convex hull chains starting from $p_0$. 

%%%%%%
\SingleFig{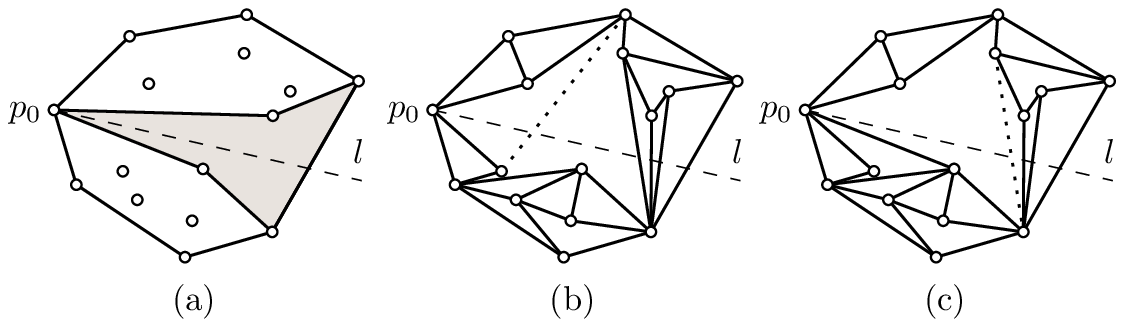}{(a) The median pseudo-triangle
from $p_0$, defined by the splitting line $l$ (dashed).
(b) The first edge intersected by $l$ (dotted) is to be flipped.
(c) After the flip, the number of intersections with $l$ is reduced by
one, and another edge is to be flipped.} 
%%%%%%

To obtain the desired flip diameter, it suffices to show that it takes $O(n)$ flips to go from any \ppsT\ $T$ to one with  the median pseudo-triangle (then apply recursion on the two
  halves). For this, we 
  flip one by one the $O(n)$ edges of $T$ that intersect the median
  line $l$, taking always the edge whose intersection with $l$ is
closest to $p_0$. Each flip may be a deletion flip or 
 a diagonal flip but:
  \begin{itemize}
  \item It is never an insertion flip. In particular, if the original triangulation $T$ is pointed, all intermediate ones are pointed, too.
  \item If it is a diagonal flip, the inserted edge is part of a geodesic from $p_0$ to the opposite corner of a pseudo-quadrilateral; hence it does not intersect the median line $l$. 
  \end{itemize}
  This last remark guarantees that at each step we have one
  intersection less between $l$ and our pseudo-triangulation. When no
  interior edge intersects $l$, the pseudo-triangulation must
  necessarily use the median pseudo-triangle.
\end{proof}

These upper bounds are not known to be tight; no better bound than the trivial lower bound of $\Omega(\nrVert)$ is known. However, they are much better than the (worst-case) diameter of the graph of
diagonal flips between triangulations of a point set, which is known to be quadratic. This implies that one can flip much faster between triangulations of a point set if pseudo-triangulations are allowed as intermediate steps. As a side remark, we mention that 
an even better, linear bound is known for triangulations using the so-called \emph{geometric bistellar flips} (see~\cite[Section 1.3]{santos-icm06}). In this case, 
only triangulations appear as intermediate steps, 
but in addition, deletions/insertions of interior \emph{vertices} of degree three are allowed.

\subsubsection*{Constrained pseudo-triangulations}
Constrained subdivisions require the usage of certain prescribed edges. We
extend this notion to include  a subset of vertices prescribed to be pointed,
while the rest are free to be either pointed or not.
Let $V$ be a subset of $\pts$ containing no corners of $R$ and
let $E$ be a set of interior edges in $(R,\pts)$
with the property that every $p\in V$ is pointed in $E$. We call \emph{pseudo-triangulations
of $(R,\pts)$ constrained by $E$ and $V$} all pseudo-triangulations
whose graph contains $E$ and whose pointed vertices contain $V$.

\begin{theorem}
\thmLab{constrained}
%In the above conditions, let
%Then,
 The graph of pseudo-triangulations  of $(R,\pts)$ constrained by $E$
%\complaint{turned G to E. GR} 
and $V$
is non-empty, connected, and regular of degree
\[
\nrVert+2\nrInt-3 - c,
\]
where $c =|E|+|V|$.
\end{theorem}

This statement generalizes both parts of Theorem~\thmRef{ptgraph} ($c=0$ and
$c=\nrInt + \nrReflex$, respectively).

\begin{proof}
Theorem~\thmRef{maximal} implies that the graph is not empty. Regularity follows
from part (1) of Theorem~\thmRef{ptgraph}, since each constraint forbids exactly one flip.
Connectedness can be proved with arguments similar to those in Theorem~\thmRef{ptgraph},
and is also a consequence of Theorem~\thmRef{pt-polytope}.
\end{proof}

\begin{cor}
%In the above conditions, 
If only vertex constraints are present \textup(that is, if $E=\emptyset$\textup)
then the diameter of the graph of flips between constrained pseudo-triangulations
is bounded by $O(n\log n)$.
\end{cor}
\begin{proof}
%Also, since in our proof of Theorem~\thmRef{diameter} the paths of
%flips obtained never turn a pointed vertex to non-pointed, we have:
   Following our proof of Theorem~\thmRef{diameter}, one can flip
from any two constrained \psT s to a common \psT\ without ever
 turning a pointed vertex to non-pointed.
\end{proof}
It is not known whether % tempting to conjecture that
 the same bound holds when edge constraints are allowed.

%\complaint{Do we conjecture that the diameter of this graph is $O(d\log d)$,
%where $d=\nrVert+2\nrInt-3 - c$ is the degree of regularity? Paco.
%Don't know. GR.}

%------------------------------------------------------------------------                Degree bounds        ------

\subsection{Vertex and Face Degree Bounds}
Pseudo-triangles can have arbitrarily many edges. However, with a simple
argument one can show that every point set has pointed pseudo-triangulations
with bounded face-degree:

\begin{theorem}
[Kettner et al.\ \cite{speckmann:kettner:etAl:tightDegreeBounds:2003}]
Every point set in general position has a pointed
pseudo-triangu\-lation consisting only of triangles and four-sided
pseudo-triangles.
\end{theorem}

\begin{proof}
Triangulate the convex hull
 of $\pts$ and then insert the interior points one by one via
two edges each. It is easy to see that if $p_i$ is a point in the
interior of a triangle or  four-sided pseudo-triangle $\Delta$, then it
is always possible to divide $\Delta$ into two triangles or four-sided
pseudo-triangles by two edges incident to $p_i$.
\end{proof}

More surprising is the result %, from the same paper, 
that the min-max
vertex degree can also be bounded by a constant. Observe that for
triangulations  the situation is quite different: In \emph{every}
triangulation of the point set  in
Figure~\figRef{circle-chain}c in Section~\secRef{number}
below, the top vertex has
 degree $\nrVert -1$.

\begin{theorem}[Kettner et al.~\cite{speckmann:kettner:etAl:tightDegreeBounds:2003}]
\thmLab{degreebound}
Every point set $\pts$ in general position has a pointed
pseudo-triangu\-lation whose maximum degree is at most five.

The bound five cannot be improved.
\end{theorem}
The method used in the following proof
gives rise to an algorithm which, with appropriate data
structures, runs in $O(\nrVert \log \nrVert)$ time.
\begin{proof}[Proof (Sketch)]
We construct the pointed pseudo-triangulation by successively
refining a \emph{partial pseudo-triangulation}, by which we mean a
partition of the convex hull of $\pts$ into some (empty)
pseudo-triangles and some convex pointgons.

We start with the edges of the convex hull of the given point set,
which defines a convex pointgon, as in Figure
\figRef{partitionPrune}(a). At each subsequent step, one of the
following two operations is used to subdivide one of the current
convex pointgons $(R',\pts')$:

\paragraph{\bf{Partition}} Choose
a vertex $p_i$ and an edge $p_jp_k$ of $R'$ not incident to $p_i$.
Choose also a line passing through $p_i$ and crossing $p_ip_j$. If
generic, this line splits $\pts'$ into two subsets with $p_i$ as
their only common point. Then, subdivide $R'$ into the convex hulls
of these two subsets (two convex pointgons) plus the pseudo-triangle
with corners $p_i$, $p_j$ and $p_k$ that gets formed in between.
Except for the degenerate case described below, which produces only
one pointgon, the degree of $p_i$ increases by $2$ and the degrees
of $p_j$ and $p_k$ by one. See Figure \figRef{partitionPrune}(b).
\paragraph{\bf Prune} A degenerate situation of partitioning arises
when one of the two subsets consists only of $p_i$ and one of $p_j$
and $p_k$ (say, $p_j$). (For this it is necessary, but not
sufficient, that  $p_ip_j$ is also a boundary edge of $R'$). The
resulting partitioning produces a pseudo-triangle and only one new
pointgon; the other one degenerates to a line segment and is
ignored. The degrees of $p_i$ and $p_k$ increase by one. See Figure
\figRef{partitionPrune}(c).

\SingleFig{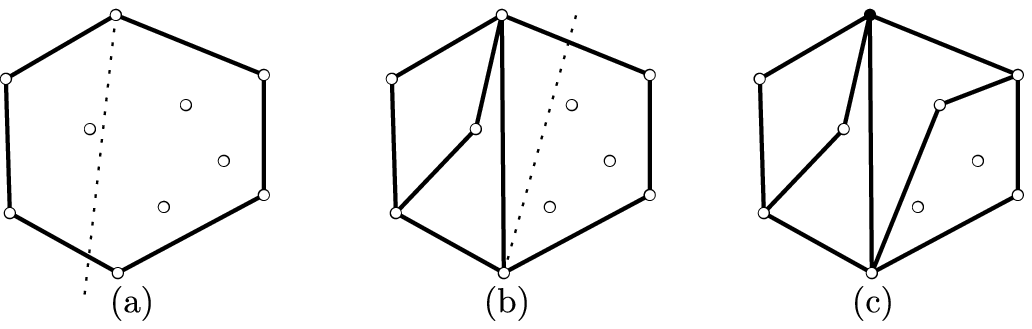}{(a) Initial convex pointgon, (b) a \emph{partition} step and (c) a \emph{prune} step on the right convex
pointgon from the previous step, pruning the black vertex on top.}

Pruning and partitioning  maintain both pointedness and planarity, and
%If performed again and again starting with the whole convex pointgon,
eventually they must lead to a pointed pseudo-triangulation.
The rest of the proof consists in selecting
 these operations
in the right order to satisfy some cleverly chosen
invariants on the degrees of the boundary points
of each convex subpolygon,
and in this way guarantee that the degree does not exceed~five.

To show that the bound five in the theorem cannot be reduced, Kettner et
al.~\cite{speckmann:kettner:etAl:tightDegreeBounds:2003} proved that, for the
vertex set of a regular $(2n+1)$-gon ($n\ge 5$) together with
its center, every pointed pseudo-triangu\-lation has some vertex of degree at
least five.
\end{proof}

% By Theorem~\thmRef{count2}, the average degree of vertices in a
% pointed pseudo-triangu\-lation of $\nrVert$ points is $4-
% 6/\nrVert$. So, for  $\nrVert>6$ the previous result is tight unless
% the set has \emph{very regular} pointed pseudo-triangulations with
% maximum degree equal to $4$. In the same paper, Kettner at
% al.~construct arbitrarily large point sets without them, proving
% that the  \emph{five}  in Theorem~\thmRef{degreebound} cannot be
% improved.

% \begin{prop}[Kettner et al.~\cite{speckmann:kettner:etAl:tightDegreeBounds:2003}]
% Let $\pts$ be the vertex set of  a regular $(2n+1)$-gon \textup($n\ge
% 5$\textup) together with its center. Every pointed
% pseudo-triangu\-lation of $\pts$  has some vertex of degree at least five.
% \qed
% \end{prop}

%------------------------------------------------------------------------      Enumeration and counting  ------

\subsection{Algorithms for Enumeration and Counting}
\secLab{enumeration}

In order to perform computer experiments that support or disprove
statements, it is useful to have algorithms that \emph{enumerate}
all \psT s of a given point set $\pts$ explicitly. There are two
algorithms for doing this in the literature.  Both traverse an
enumeration tree that is implicitly built on top of the graph of
\ppsT s.

The algorithm of
Bereg~\cite{bereg:enumerating-pseudo-triangulations:2005} is based on the
reverse search paradigm of Avis and
Fukuda~\cite{avis:fukuda:reverse-search:1996}. It takes $O(n)$ space, and its
running time is $O(\log n)$ times the number of \ppsT s.

Another enumeration algorithm has been given by Br{\"{o}}nnimann, Kettner,
Pocchiola, and Snoeyink~\cite{bronnimann:kettner:pocchiola:snoeyink:CountingAndEnumeratingPointedPseudoTriang:2006}.
They developed the \emph{greedy flip} algorithm,  based
on the analogous algorithm by Pocchiola and
Vegter~\cite{pocchiola:vegter:topoSweep:1996} for the case of \psT s of convex
objects (cf.~Figure~\figRef{pseudoTriangPV} and Section~\secRef{visibility-complex}),
and on its generalization by Angelier and Pocchiola~\cite{angelier:pocchiola:sumofsquares:2003}.
%
%\complaint{Added reference~\cite{angelier:pocchiola:sumofsquares:2003} for Referee 1. Paco, Septemeber 30, 2007}
%
  % binary
The enumeration tree that the algorithm uses is a binary tree, and may contain
``dead ends'', whose number can only be analyzed very crudely.  The algorithm
takes $O(n)$ space and the proved upper bound on the running time is
$O(n\log n)$ times the number of \ppsT s.
But the algorithm has been implemented and, in practice, it seems to
need only $O(\log n)$ time per \ppsT.  It can also be adapted
to constrained \ppsT s, where a subset of the edges is held fixed.

The most stringent bottleneck to the applicability of these enumeration algorithms is not
the time per \psT, but the exponential growth of the number of \ppsT s, see
Section~\secRef{number}.

Other approaches to enumeration are conceivable.  In particular the known
enumeration algorithms for vertices of polytopes can be applied to the
polytopes of \psT s that are mentioned in Section~\secRef{polytope}.  This
would also lead to algorithms for enumerating all (pointed and non-pointed)
\psT s of a pointgon, or of pseudo-triangulations constrained
in the sense of Theorem~\thmRef{constrained}.  These approaches have not been developed so far.

If one just wants to \emph{count} \psT s, it is not necessary to enumerate
them one by one.  A divide-and-conquer algorithm for counting (pointed or
arbitrary) \psT s is given by Aichholzer
et~al.~\cite{streinu:aichholzer:rote:speckmann:zigZagPath-wads:2003}.  A
constraint set $V$ of vertices which must be pointed can be specified.

%------------------------------------------------------------------------      Number of pseudo-triangu\-lations  of point sets  ------

\subsection{The Number of Pseudo-Triangu\-lations of a Point Set}
\secLab{number}
What is the minimum and maximum number of pseudo-triangu\-lations of
a point set $\pts$, for a fixed cardinality $n$ of $\pts$? Before
going on, let us summarize what is known about the analogous
question for triangulations. We use the notations $\Theta^*$,
$\Omega^*$ and $O^*$ to indicate that a polynomial factor has been
neglected.

\begin{itemize}
\item For points in convex position, the number of triangulations
(and of pseudo-triangu\-lations) is the Catalan number
$C_{n-2}=\frac{1}{n-1}\binom{2n-4} {n-2}$. Asymptotically, this
grows as $\Theta(4^n n^{-3/2})$, or $\Theta^*(4^n)$.

\item The number of triangulations of an arbitrary point set
in general position is at most $O^*(43^n)$
\cite{Sharir+Welzl:RandomTriangulations:2006} and at least
$\Omega^*(2.33^n)$
\cite{aichholzer:hurtado:noy:lowerBoundTriangulations:2004}. Refined
versions, for $i$ interior and $h$ convex hull points, are known: an
upper bound of $O^*(43^i 7^h)$ from
\cite{Sharir+Welzl:RandomTriangulations:2006} and a lower bound
of $\Omega(2.72^h 2.2^i)$  (or $\Omega(2.63^i)$, for fixed $h$)~\cite{mccabe:seidel:newLowerBoundsTriangulations:2002}.

\item The point sets with the minimum and maximum number of triangulations
known have asymptotically $\Theta^*(\sqrt{12}^n)$ and $\Theta^*(\sqrt{72}^n)$
triangulations~\cite{aichholzer-etal:number-planegraphs:2007}.
The first one is the so-called \emph{double circle}, consisting of a
convex $n/2$-gon and a point very close to the interior of every
edge of it. The second one is a variation of the so-called {\em
double chain}, consisting of two convex $n/2$-gons \emph{facing each
other} so that each vertex of one of them sees all but one edges of
the other. See these point sets in Figure~\figRef{circle-chain}.
\SingleFig{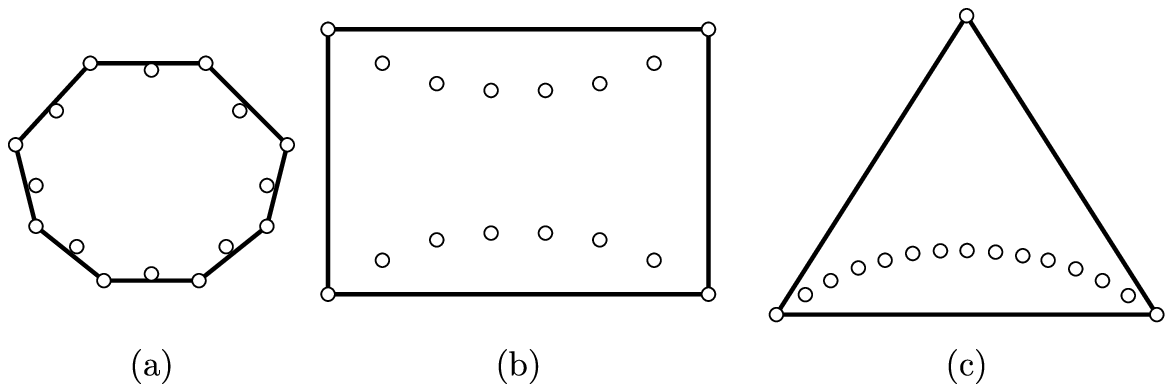}{(a) A double circle; (b) a double chain; (c) a single chain, all with $16$ points.}
\end{itemize}

In the case of pseudo-triangulations, 
the minimum possible number is attained by points in convex position, even if we only count \emph{pointed} pseudo-triangulations:

\begin{theorem}
[%Aichholzer et al.
\cite{aichholzer:aurenhammer:krasser:speckmann:convexityMinimizes:2004}]
\thmLab{convexityminimizes} 
Every point set in general position has at least as many pointed
pseudo-triangu\-lations as the convex polygon with the same number
of points.
\end{theorem}

This follows from the following lemma, taking into account that each Catalan number is less than four times the next one:

\begin{lemma}
Let $\pts$ be a set of at least five
points, at least one of them interior and let $p_0$ be an interior
point of $\pts$. Then, $\pts$ has at least four times as many
pointed pseudo-triangulations as $\pts\setminus \{p_0\}$. 
\end{lemma}

\begin{proof}[Proof (Sketch)]
From each pointed pseudo-triangulation $T$ of $\pts\setminus
\{p_0\}$ one can construct (at least) four pointed
pseudo-triangulations of $\pts$ as follows:
\begin{itemize}
\item[(a)] Three by a Henneberg step of type~1,
as introduced in Theorem~\thmRef{henneberg}. 
%That is, inserting two
%of the three geodesics that join $p_0$ to corners of the pseudo-triangle of $T$ that contains $p_0$
%in its interior.
\item[(b)] One  obtained by a Henneberg step of type~2.
That is, by performing a diagonal flip of an edge $e$ opposite to $p_0$
in one of the three ``Henneberg~1" pseudo-triangulations of the previous paragraph.
\end{itemize}

The tricky part of the proof, which we omit, is to show that for at least one of the three \psT\ in part (a)
there is at least one choice of $e$ that indeed increases the degree of $p_0$ from 2 to 3
to get the \psT\ of part (b).
(Observe that, contrary to what happens in triangulations, a diagonal flip may not increase
the degree of the opposite corners, since the diagonals of a pseudo-quadrilateral may not be incident
to the corners; see Figure~\figRef{ptflips}).

We also omit the proof that the list contains no repetition. 
%That is, that from the four pseudo-triangulations
%assigned to $T$ we can actually recover $T$. For the three with degree two at $p_0$ this is obvious. For the fourth one, call it $T''$ it is also true:
% $p_0$ is pointed in $T''$ and has degree three, hence one of its edges is ``between''
%the other two. $T$ is recovered by first flipping this edge and then removing the two edges of $p_0$ that are left.
\end{proof}

%\complaint{Changed notation to $V_B$ and $V_I$, to match $\nrCH$ and
%$\nrInt$. Paco.}
To get finer statements,
it is convenient to stratify the set of
\psT s of a point set $\pts$ according to the set of
pointed vertices. For this, let  $V_I $ be the set of interior points of
$\pts$. For each subset $V_\npointed \subseteq V_I$ let
$\PT(V_\npointed)$ denote the set of pseudo-triangu\-lations of $V$ in
which the points of $V_\npointed$
are pointed and the remaining vertices $%V_{np} =
V_I \setminus V_\npointed$ are
non-pointed. For example,  $\PT(\emptyset)$ and $\PT(V_I)$ are the
triangulations and the pointed pseudo-triangu\-lations of $\pts$,
respectively. The following is easy to prove:

\begin{prop}
[%Santos et al.
\cite{santos:randall:rote:snoeyink:countingTriangulationsWheels:2001}]
\propLab{monotone}
For every point set in general position, for any
subset $V_\npointed$ of interior points designated as pointed, and
for every point $p_0\in V_\npointed$:
\[
|\PT(V_\npointed)|\le 3 \, |\PT(V_\npointed\backslash \{p_0\})|
\]
\end{prop}
%
%\begin{proof}
%Let us consider the graph of insertion/deletion flips that relate
%$\PT(V_\npointed)$ and $\PT(V_\npointed\backslash \{p_0\})$. This is a bipartite graph
%in which a pseudo-triangulation of $\PT(V_\npointed)$ is joined to the unique one
%obtained by the insertion flip that turns $p_0$ from pointed to non-pointed.

%The statement follows from the claim that no pseudo-triangulation of
%$\PT(V_\npointed\backslash \{p_0\})$ has degree more than three in this graph. That is,
%that no more than three edges incident to any given vertex $p_0$ produce
%deletion flips that turn $p_0$ from non-pointed to pointed. This holds since
%for such an edge $e$, the two angles incident to $e$ at $p_0$ must add to more than 180 degrees,
%and this cannot happen for more than three edges (in fact, it can only happen for two edges
%unless $p_0$ has degree three).
%\end{proof}

In other words, the number of \psT s \emph{does not increase too much} if the prescription for a point
changes from non-pointed to pointed.
Experience and partial results show that
 the number  actually \emph{decreases}:
%the opposite usually happens. More precisely, there is the following conjecture:
%
%
\begin{conjecture}
\conjLab{monotone}
For every point set in general position, for any subset $V_\npointed$ of
interior points designated as pointed, and for every point $p_0\in
V_\npointed$:
\[
|\PT(V_\npointed)|\ge \, |\PT(V_\npointed\backslash \{p_0\})|
\]
\end{conjecture}

This conjecture is known to hold for sets with a single interior point~\cite{santos:randall:rote:snoeyink:countingTriangulationsWheels:2001} and for the following
three specific families of point sets~\cite{aichholzer:orden:santos:speckmann:numberPseudoTriangulationsCertainPointSets:2006+}: the double circle,
the double chain, and the third point set of Figure~\figRef{circle-chain}. This last set, called a ``single-chain''
consists of a convex $n-1$-gon together with a point that sees all of its edges except one.
% Besides proving the conjecture,
The asymptotic numbers of pseudo-triangulations of these point sets are also computed
in \cite{aichholzer:orden:santos:speckmann:numberPseudoTriangulationsCertainPointSets:2006+},
and summarized in
the following table.
\medskip
\begin{center}
\begin{tabular}{|c||c|c|c|}
\hline
   & double  & double  & single  \\
   & circle    & chain   & chain   \\
\hline triangulations
    &$\Theta^*(\sqrt{12}^n)$ &   $\Theta^*(8^n)$ & $\Theta^*(4^n)$    \\
\hline pointed pseudo-triangu\-lations
    &$\Theta^*(\sqrt{28}^n)$   &    $\Theta^*(12^n)$ &   $\Theta^*(8^n)$ \\
\hline all pseudo-triangu\-lations
    &$\Theta^*(\sqrt{40}^n)$  &   $\Theta^*(20^n)$ &   $\Theta^*(12^n)$ \\
%\hline Conjecture \conjRef{monotone} holds?
%      &YES       &      YES  &        YES       \\
\hline
\end{tabular}
\end{center}
\medskip
The number of triangulations of the single chain is just a Catalan
number.
It may come as a surprise that the double circle, which has as few
triangulations as known so far, still has much more pointed
pseudo-triangu\-lations than the single chain, or the convex
$n$-gon. But this is a consequence of Theorem~\thmRef{convexityminimizes}.

To finish this section, as a joint application of Theorem~\thmRef{convexityminimizes} and
Proposition~\propRef{monotone} we obtain the following lower bound on the
size of $\PT(V_\npointed)$:

\begin{cor}
For a point set with $h$ points on the convex hull and $i$ in the
interior, and for every set $V_\npointed$ of $k=|V_\npointed|$ interior points
designated to be pointed:
\[
\PT(V_\npointed) \ge \frac{\PT(V_I)}{3^{i -k}}\ge \frac{C_{h+i-2}}{3^{i-k}}
= \Theta^*(4^h (4/3)^i 3^{k}).
\]
In particular, the total number of pseudo-triangu\-lations is at
least $\Omega^*(4^h ({16}/{3})^i)$.
\end{cor}

\begin{proof}
The first inequality comes from applying
Proposition~\propRef{monotone} one by one to the non-pointed vertices
in $V_I\setminus V_\npointed$. The second inequality is
Theorem~\thmRef{convexityminimizes}. The total number of
pseudo-triangu\-lations equals
$$
\sum_{V_\npointed\subseteq V_I} \PT(V_\npointed) \ge \sum_{V_\npointed\subseteq V_I}
\frac{C_{h+i-2}}{3^{i-|V_\npointed|}} = \frac{C_{h+i-2}}{3^{i}}
\sum_{V_\npointed\subseteq V_I} {3^{|V_\npointed|}} = \frac{C_{h+i-2}}{3^{i}}
{4^{i}}.
\eqno \qedhere
$$
\end{proof}

%------------------------------------------------------------------------      Number of pseudo-triangu\-lations  of polygons  ------

\subsection{The Number of Pointed Pseudo-Triangulations of a Polygon}
\secLab{number-k-gon}

The largest and smallest possible number of pointed
pseudo-triangulations of a polygon with $\nrCorners$ corners) are
 easy to obtain:

\begin{theorem}
A pseudo-$\nrCorners$-gon has between $2^{\nrCorners-3}$ and $C_{\nrCorners-2}$
\textup(the Catalan number\textup) \ppsT s. Both bounds are achieved.
\end{theorem}

\begin{proof}
The upper bound is achieved by a convex $\nrCorners $-gon.
The lower bound  is achieved by the pseudo-$\nrCorners $-gon of Figure~\figRef{pseudo-k-gon}
whose diagonals  come in $\nrCorners-3$ crossing pairs.
Every choice of one diagonal from each pair gives a \ppsT.

\SingleFig{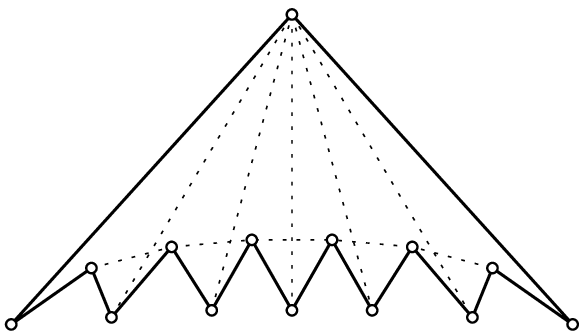}{A pseudo-$\nrCorners $-gon with $2^{\nrCorners-3}$ pseudo-triangulations.}%{0.5}

%\complaint{Changed (corrected) this paragraph. Paco, November 27,2006.o.k.GR}
To prove the lower bound, let $e$ be a diagonal in $R$.
The diagonal $e$ divides
$R$ into two polygons $R'$ and $R''$ with $\nrCorners'$ and $\nrCorners''$ corners respectively, 
with $\nrCorners' + \nrCorners'' =  \nrCorners + 2$. Then,
the number of pointed pseudo-triangulations of $R$ that contain this diagonal equals the product of the numbers of pointed pseudo-triangulations of $R'$ and $R''$. By inductive hypothesis this
gives at least
\[
2^{\nrCorners'-3}\cdot 2^{\nrCorners''-3} = 2^{\nrCorners-4}
\]
pointed pseudo-triangulations. But
the number of \ppsT s that do not use $e$ is at least the same number:
to each pseudo-triangulation $T$ that uses $e$ we associate the one obtained by the flip at $e$, and no two choices of $T$ produce the same $T'$, by
Lemma~\lemRef{switch}  below.

For the upper bound, consider the $\nrCorners $ corners of $R$
corresponding cyclically to the $\nrCorners$ vertices of a convex
$\nrCorners$-gon. To every triangulation $\hat T$ of the $\nrCorners
$-gon we associate the \ppsT\ $T$ that uses \emph{the
same} geodesics. (This correspondence will be important again in
Section~\secRef{geodesic}, see
Figure~\figRef{geodesicTriangulation}).
That every pseudo-triangulation $T$ of $R$ arises in this way can be proved
using Lemma~\lemRef{tangents}: To each bitangent of $T$ we associate its canonical
geodesic, and consider the corresponding set of diagonals in the $k$-gon. These diagonals
are mutually non-crossing, hence there is a triangulation $\hat T$ containing all of them.
\end{proof}

\begin{lemma}
\lemLab{switch}
Let $T$ be a pseudo-triangulation of a pointgon $(R,\pts)$ and $e$ a possible edge that is not used in
$T$. Then, there is at most one diagonal flip in $T$ that inserts $e$,
unless one of the end-points of $e$ is an interior vertex of degree
two in $T$, in which case there may be two.
\end{lemma}

\begin{proof}
If an edge of
$T$ crosses $e$, then only the flip on that edge can insert $e$.
If no edge of $T$ crosses $e$, then $e$ lies within a certain pseudo-triangle $\Delta$ of $T$.
We regard the diagonal flip as obtained by first inserting $e$
and then deleting another edge.
(We do not need the intermediate
graph to be a pseudo-triangulation, although it follows from our proof that indeed it is.)

Since a pseudo-triangle has no bitangent, $e$ is not tangent to $\Delta$ at (at least) one of its end-points. On the other
hand, if $e$ can be inserted by a diagonal flip,
$e$ must be tangent to $\Delta$  at one of the ends, because otherwise the insertion of $e$ turns both end-points from
pointed to non-pointed and it will be impossible to make them both pointed again by the removal of a single edge.
Hence, $e$ is tangent to $\Delta$ at exactly one of its end-points. The other one, let us call it $p_0$,
is a reflex vertex of $\Delta$. The edge removed by the flip must be one of its extremal edges,
 $e_1$ and $e_2$. We now have two cases:
\begin{enumerate}
\item If $p_0$ has degree two in $T$ then any of the two edges incident to it produces a flip that inserts $e$.
\item If there is another edge $e'$ incident to $p_0$ in $T$ besides $e_1$ and $e_2$, then only the flip at the $e_i$ that lies in the
reflex angle formed by $e$ and $e'$ can possibly insert $e$ by a diagonal edge, since $e$, $e'$ and that $e_i$ make $p_0$
non-pointed.
\qedhere
\end{enumerate}
\end{proof}

This lemma also shows that every pseudo-$k$-gon has at least $2(k-3)$ diagonals: the $k-3$ forming
a pointed pseudo-triangulation plus the $k-3$ different (by the lemma) ones inserted by flips in it. So,  the pseudo-$k$-gon of Figure~\figRef{pseudo-k-gon} is also minimal in this sense.

%\medskip

Concerning possibly non-pointed pseudo-triangulations of a polygon, the analogue of
Proposition~\propRef{monotone} is true, with a similar proof but a better constant:

\begin{prop}
For every polygon, for any subset $V_\npointed$ of reflex vertices of it
 designated as pointed, and for every point $p_0\in
V_\npointed$:
\[
|\PT(V_\npointed)|\le 2 \, |\PT(V_\npointed\backslash \{p_0\})|
\eqno \qed
\]
\end{prop}
%
%
%\begin{proof}
%Similar to that of Proposition~\propRef{monotone}.  Observe that in that proof we only needed a $3$
%in the bound in order to take into account the possibility of an interior non-pointed vertex of degree three, in which the removal of any of the three edges produces a pseudo-triangulation. Now we do not have that case, and only two edges incident to $p_0$ can at most produce edge-removing flips.
%\end{proof}

%It is also natural to think that Conjecture~\conjRef{monotone} extends to pseudo-triangulations of polygons (and of pointgons).

%!TEX root = article.tex

%------------------------------------------------------------------------      Liftings and Surfaces    ------

\section{3D Liftings and Locally Convex Functions}
\secLab{liftings}

We switch now to a three-dimensional geometric problem which leads naturally to
pseudo-triangulations of pointgons: locally convex polyhedral
surfaces.
This section is heavily based on Sections 3--5 of Aichholzer et
al.~\cite{aichholzer:aurenhammer:brass:krasser:pseudoTriangulationsNovelFlip:2003}, but
some of the proofs are new. Specially, that of Lemma~\lemRef{downflip} is more direct
than the original one.

\subsection{The Lower (Locally) Convex Hull}

Suppose a set of data points $p_i=(x_i,y_i)$ in the plane with
associated height values $h_i$ is given, and we look for the \emph{highest}
convex function $f\colon \reals^2\to \reals$ that remains below the given
height values:
\begin{equation}
  \eqLab{below}
f(x_i,y_i) \le h_i\text{, for all $i$}
\end{equation}
Then it is well-known that the function $f$ will be piecewise linear. It is
defined on the convex hull of the point set $P$, and its
graph is the lower convex hull of the points $(x_i,y_i,h_i)\in\reals^3$,
i.\,e., the part of the convex hull that is seen from below.
If the height are sufficiently generic,
the pieces where $f$ is linear are triangles, forming a triangulation.
(The resulting triangulations are usually called \emph{regular
  triangulations} of~$P$
\cite{gelfand:kapranov:zelevinsky:discriminantsResultants:1994,
billera:filliman:sturmfels:constructionsSecondaryPolytopes:1990,
edelsbrunner:shah:incrementalTopologicalFlipping:1996};
but observe they may use a proper subset of $P$ as set of vertices.)
%\complaint{Added last sentence, for referee 2; Paco, September 30, 2007}

%\subsection{Locally Convex Functions}
We can ask the same question for a function $f$ that is defined over a
non-convex  polygonal region $R$. But, now, it is natural to replace
the condition of convexity  by local convexity.
A function $R\to \reals$ is called \emph{locally convex}
if it is convex on every straight segment contained in $R$:
For a pointgon $(R,\pts)$ with given heights $h_i$ at the points $p_i\in P$,
we look for the \emph{lower locally convex hull},
the (graph of the) highest function $f\colon R\to \reals$ that
fulfills \eqRef{below}.
As shown by Aichholzer et
al.~\cite{aichholzer:aurenhammer:brass:krasser:pseudoTriangulationsNovelFlip:2003},
 the regions on which $f$ is linear form
a \psT\ of a
subset of~$P$ (Theorem~\thmRef{envelope}).

Thus, \psT s arise very naturally in this context:
we start with a pointgon and some height values and construct
the lower locally convex hull. The edges of this hull, when projected to the
plane, yield a \psT.

\subsection{Liftings of Plane Graphs}
\secLab{liftings-subsec}
To study the problem of the {lower locally convex hull}, we will proceed in the
opposite direction: we take a fixed \psT\ $T$ in the plane and ask
for the piecewise linear surfaces that project onto it (liftings of~$T$).
Along the way we get  a simple but fundamental result that describes
explicitly the space of liftings of a given \psT\
(Theorem~\thmRef{liftingbasis}).
These results are relevant also for the rigidity properties
of pseudotriangulations in Section~\secRef{rigidity}.

%
%
% Liftings of 2D triangulations to surfaces in 3D are a classic topic.
% It is clear that once a lifting of a point set $P$ has been fixed,
% every triangulation of $P$ can be lifted to a triangulated surface
% in a unique way. In particular, every choice of heights for the
% point set induces a partial ordering among the set of
% triangulations, which has proved useful in several contexts. The
% unique minimal element in this ordering, that is, (the projection
% of) the lower envelope of the lifted point set, is called a {\em
% regular} triangulation
% \cite{gelfand:kapranov:zelevinsky:discriminantsResultants:1994,
% billera:filliman:sturmfels:constructionsSecondaryPolytopes:1990,
% edelsbrunner:shah:incrementalTopologicalFlipping:1996} of $P$. The
% prototypical example is the Delaunay triangulation of $P$, obtained
% when the point set is lifted to lie in a paraboloid of vertical
% revolution. The space of regular triangulations has specially nice
% properties, such as the fact that between two of them there is
% always a ``monotone'' path of flips.

% Let $R$ be a polygonal domain, which we do not assume to be simply
% connected,

% For a pseudo-triangu\-lation of a point set, pending vertices are
% the pointed and interior vertices, while complete are the boundary
% or non-pointed vertices.

\begin{definition}
\defLab{lifting}
Let $G$ be a plane straight-line graph and let $R$ be a union of (closed) faces of $G$.
A \textup(3D\textup) \emph{lifting} of $(G,R)$ is the graph of a continuous
function $f\colon R\to \reals$ that restricts to an affine-linear function on
each face of $G$, see Figure~\figRef{locallyconvex}a.
 In other words, every face of
$G$ is lifted to a planar face in space.
\end{definition}

\SingleFig{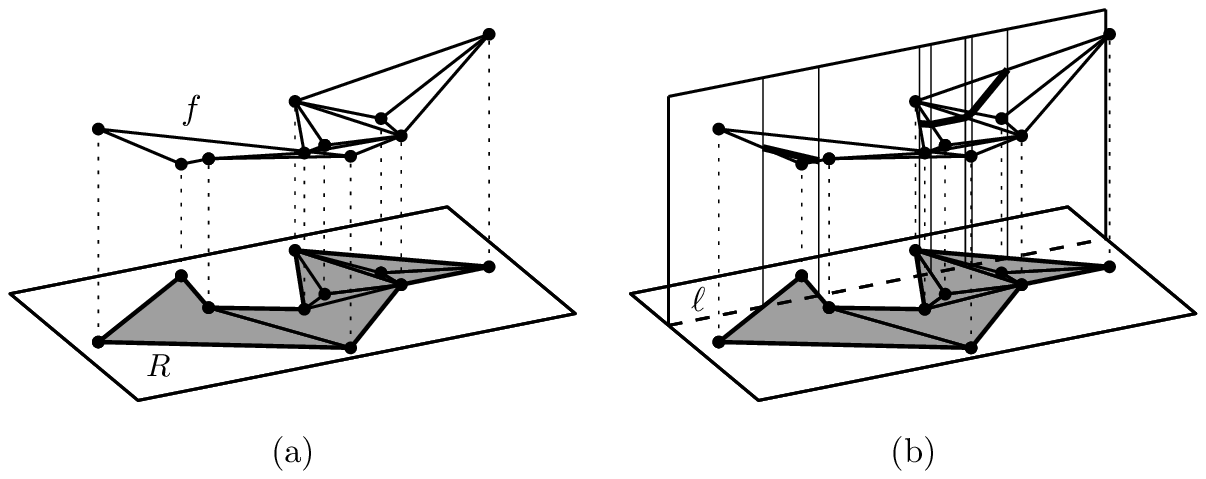}{(a) A 3D lifting $f$ of a geometric graph
  (in this case, a \psT) over a plane region $R$.
(b) This lifting is locally convex. Looking at the restriction of $f$ to a line~$\ell$,
one sees that $f$ cannot be extended to a convex function over~$\reals^2$.}

As an important special case we have liftings of the whole plane $R=\reals^2$.
In this case, we usually insist that $f$ is identically zero
on the outer face, which can always be done by subtracting from a given $f$ the
affine-linear function it coincides with in the outer face.
By the Maxwell-Cremona theorem  (see Theorem~\thmRef{MC}),
these 3D liftings have correspondences to other objects: % namely
reciprocal diagrams, which are treated in Section~\secRef{reciprocal},
and self-stresses, which are treated in Section~\secRef{self-stresses}.
%the results that we state in this section are basic.

On the other hand, when $G$ is a \psT\ of a pointgon $(R,P)$,
 we don't care about
the outer face, and we consider $f$ defined only on the domain $R\subset \reals^2$.
  In this case,
the boundary vertices need not be coplanar in the lifting. We will
not distinguish between $f$ as a function and the lifting as a
three-dimensional surface (the graph of the function).

The following easy observation lies at the heart of
many proofs that use liftings to prove properties of \psT s,
and it highlights the role of pointedness.
As we did in Section~\secRef{basicProperties},
we call a vertex $p$ in a plane graph $G$ over a region $R$
\emph{\rp} if it is a reflex
vertex of some (and then a unique) region of $G$ in $R$.
Otherwise it is called \emph{\rnp}. 
%\complaint{Edited this paragraph for referee 1. Paco, September 30, 2007.}

% in the context of liftings
\begin{lemma}
  \lemLab{lifting-basic}
  Let $f$ be a lifting of a plane graph $G$ over a region $R$.
  Let $p$ be a vertex of $G$ that is incident to a reflex vertex of some face
  $F\subset R$.
    If $f$ has a global maximum at $p$, then $f$ is constant on $F$,
    \textup(and every point of the interior of $F$ is a maximum point of~$f$\textup).
\end{lemma}
\begin{proof}
On any segment through $p$ contained in  $F$, $f$ is linear, and hence
must be constant if $p$ is a maximum,
see Figure~\figRef{noglobalmax}.
By considering two non-parallel segments through $p$, we conclude that
$f$ is constant on~$F$.
\end{proof}
\SingleFig{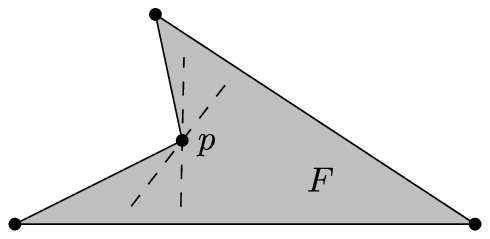}{The lifting must be horizontal on
every segment in $F$ passing through~$p$.}

If one walks across a lifted edge between two faces, the slope may increase,
decrease, or remain the same. Accordingly, we call the lifted edge a
\emph{valley edge}, a \emph{mountain edge} (or \emph{ridge}), or a \emph{flat
  edge}.   
At valley and mountain  edges, the function $f$ is (locally) convex
and concave, respectively. 
(The names valley and mountain should not be taken too literally.
The slope does not have to change from negative to positive when
we cross a valley. It must only increase.)  
%\complaint{Added sentence in
%  parenthesis, for referee 2. Paco, September 30, 2007. Modif. GR}

\begin{lemma} \lemLab{max-height}
The maximum and minimum height in every lifting is attained at some
\rnp{} vertex.
\end{lemma}
\begin{proof}
We take a convex hull vertex $p$ of the set of vertices
  where the maximum height is attained.  It follows from
  Lemma~\lemRef{lifting-basic} %(a) % part~\ref{global-max}
 that $p$ cannot be a \rp{}
  vertex.
\end{proof}

\subsection{Liftings of \PsT s}
\secLab{pst-liftings}

Let $T$ be a fixed \psT\ of a pointgon $(R,\pts)$.

Clearly, the liftings of $T$ form a vector space that can be
considered a subspace of the space $\reals^{\pts}$ of all maps $P\to
\reals$. The constraints for the heights of the vertices of $T$ to
define a lifting are linear equalities: for each \rp{} vertex $p_l$,
there is a linear equation specifying that the lifting of $p_l$ lies
in the plane containing the three lifted corners $p_i,p_j,p_k$ of
the unique pseudo-triangle in which $p_l$ is reflex.  More
precisely, since $p_l$ lies in the convex hull of $p_i,p_j,p_k$,
there are unique coefficients $\lambda_i,\lambda_j,\lambda_k$ with
$p_l=\lambda_ip_i+\lambda_jp_j+\lambda_kp_k$, with
$\lambda_i+\lambda_j+\lambda_k=1$ and
$0<\lambda_i,\lambda_j,\lambda_k<1$. Then the equation for the
heights $z$ is
\begin{equation}
 \eqLab{z}
z_l=\lambda_iz_i+\lambda_jz_j+\lambda_kz_k
\end{equation}
This equation is the algebraic reason
behind the geometric proof of
 Lemma~\lemRef{lifting-basic}: %(a): % part~\ref{global-max}:
$z_l$~is a convex combination of the heights $z_i,z_j,z_k$ of the three corners
of the pseudo-triangle in which $p_l$ is reflex.

The following result provides an explicit basis of the vector space of
liftings, hence it shows what its dimension is.

\begin{theorem}[The Surface Theorem, Aichholzer et al.~\cite{aichholzer:aurenhammer:brass:krasser:pseudoTriangulationsNovelFlip:2003}]
\thmLab{liftingbasis}
Let $T$ be a pseudo-triangu\-lation of a pointgon $(R,P)$.
\begin{enumerate}
\item[(i)]%\label{basis}
For every choice of heights for the  \rnp{} vertices of $T$,
 there is a unique lifting of $T$.
\item[(ii)]%\label{monotone}
The height of every \rp{} vertex is a linear function of
the height of the \rnp{} vertices, with nonnegative coefficients.
\end{enumerate}
\end{theorem}

\begin{proof}
We offer a geometric proof, different from the more
algebraic original proof of~\cite{aichholzer:aurenhammer:brass:krasser:pseudoTriangulationsNovelFlip:2003}.
%  To prove parts (\ref{basis}) and (\ref{monotone}),
 We use induction on the
  number of \rp{} vertices. If this number is zero, then every pseudo-triangle
  of $T$ is a triangle and the statement is obvious. Otherwise
% If the number of \rp{} vertices is at least 1,
  choose a \rp{} vertex $p$ of $T$. Let $\Delta$
be the pseudo-triangle of which $p$ is a reflex vertex and let $T'$
be the  \psT\ obtained by the edge-inserting flip at $p$.  The
\rnp{} vertices of $T'$ are those of $T$ plus $p$ itself.  Hence, by
the inductive hypothesis, for every choice of height at this new
vertex there is a unique lifting of $T'$. Moreover, the heights of
\rp{} vertices depend  linearly on this choice, and the maximum and
minimum heights are always  attained at some \rnp{} vertices.

We keep the heights of the original \rnp\ vertices of $T$ fixed and
vary the height of~$p$. When the height chosen for $p$ is very high,
$p$ is the global maximum of the lifting. In particular, it is above
the plane that contains the three lifted corners of $\Delta$.
Similarly, when the height of $p$ is very low, $p$ is below that
plane. Linearity implies that there is a unique height for $p$ that
makes  $p$ and the corners of $\Delta$ coplanar. This proves
part~(i). % (\ref{basis}). 
The linear dependence of the heights of \rp\
vertices on the given heights follows from the fact that the space
of liftings is a linear space.

To prove monotonicity in part~(ii), %\ref{monotone}),
a similar argument
works. If the dependence were not monotone, there would be a set of
initial heights for the \rnp{} vertices of $T$ such that an increase
in one of them ($p$) makes some \rp{} vertex $q$ go down. By
linearity, this process can be extrapolated, and a large increase in
the height of  $p$ would make the \rp{} vertex $q$ go below every
\rnp{} vertex, a contradiction to Lemma~\lemRef{max-height}.
\end{proof}

\begin{cor}
\corLab{lifts-dimension} Let $T$ be a pseudo-triangu\-lation of a
pointgon $(R,P)$. Then, the linear space of its lifts has dimension
equal to the number of \rnp{} vertices.
\end{cor}

\begin{remark}[``Non-projective'' pseudo-triangulations]
\label{rem:projective}
%\complaint{Added this remark. Paco, Nov 27, 2006. very nice, thanks GR}
As a particular case of Theorem~\thmRef{liftingbasis} we recover the familiar fact that if $T$ is a triangulation (i.e., all vertices are \rnp{}), then every choice of heights for the vertices induces a lift of $T$.
But the, also familiar, fact that \emph{sufficiently generic} choices of heights produce lifts in which no two 
adjacent faces are coplanar does not hold for pseudo-triangulations. 

Consider, for example, a pseudo-triangulation $T$ of $(R,P)$ which contains a
subset of interior pointed vertices such that: (1) the graph $T'$ obtained by
removing these interior vertices and their incident edges is still a
pseudo-triangulation (of a different pointgon $(R,P')$, where $P'\subset P$),
and (2) every \rnp{} vertex of $T$ is still \rnp{} in $T'$. An example is
shown in Figure~\figRef{projective}a, in which $T'$ is obtained by removing the
three interior vertices and the dashed edges. In the terminology
of~\cite{aichholzer:aurenhammer:brass:krasser:pseudoTriangulationsNovelFlip:2003},
such a $T$ is not ``face-honest''.
Theorem~\thmRef{liftingbasis} implies that, when this happens, $T$ and $T'$ have exactly the same lifts. 
%in bijection with the 
%possible choices of heights for their \rnp{} vertices.
 In particular,  in every lift of $T$ the faces of $T$ that form a single face in $T'$ are coplanar.
 
 Theorem~5.4
 in~\cite{aichholzer:aurenhammer:brass:krasser:pseudoTriangulationsNovelFlip:2003}
 is an attempt to characterize~\emph{projective} \psT s; that is, those that
 admit lifts with no coplanar adjacent faces. Besides showing that projective
 \psT s must be face-honest in the above sense, the authors give an example in
 which a face-honest \psT\ \emph{in special position} is not projective
 (Figure~\figRef{projective}b).  But they also make the statement that if the
 vertex set $P$ is sufficiently generic then every face-honest \psT\ of
 $(R,P)$ is projective. This statement is, unfortunately, wrong, as
 Figure~\figRef{projective}c shows.  A more precise approach to
  characterize \psT s that are projective for generic positions of the vertex set
  is made in~\cite{aichholzer:aurenhammer:hackl:preTriangulations:2006}.
  See, in particular, Theorem~\thmRef{envelope} below.
  %
  %\complaint{Edited this for referee 1. Paco, September 30, 2007}
  %
\end{remark}

\SingleFig{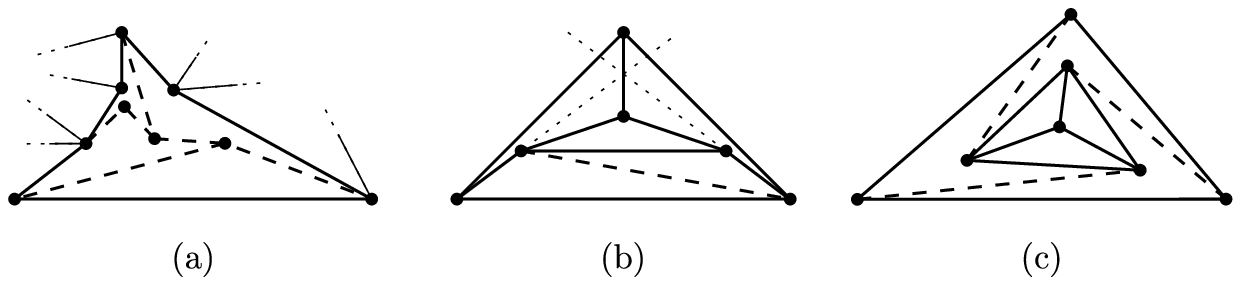}{Some ``non-projective'' pseudo-triangulations, that can only be lifted
with  flat edges (shown dashed).
(a)~A pointed \psT\ inside a pseudotriangular face (which might form
part of a larger graph) will always
be lifted flat.
(b)~In this special position, the dashed edge is flat in any lifting.
Perturbing the vertices will make the edge folded.
(c)~Even in generic vertex positions, the dashed edges are always lifted flat.}

A common trick that has also been used in the last proof is that we vary
the height of a single \rnp\ vertex, keeping the other heights fixed.
The following lemma describes the situation when this vertex is very high.
\begin{lemma} \lemLab{unique-max}
If, in some lifting, the \rnp\ vertex $p_i$ is higher than
all other \rnp{} vertices, then it is the unique global maximum
in the lifting.
\end{lemma}
\begin{proof}
As in the proof of Lemma~\lemRef{max-height},
we look at the convex hull of the set
  where the maximum height is attained. Every vertex of this convex hull
must be a \rnp\ vertex $p_j$, at its original height $h_j$.
This vertex can only be~$p_i$.
\end{proof}

The previous lemma can be rephrased as follows: if the global maximum
of a lifting is unique, then it is a \rnp\ vertex. The following
crucial local argument lies at the heart of its proof and of other proofs.
\begin{lemma} \lemLab{mountain-edge}
Let $p_i$ be a strict local maximum in some lifting
$f\colon R\to\reals$. Suppose that $p_i$ is not a corner of $R$.
Then, for every open half-plane $H$ 
with $p_i$ on its boundary, there is
either a boundary edge of $R$ or a mountain edge of 
the lift (or both) that is incident to $p_i$ and contained in $H$.
See Figure~\figRef{maximum}a\textup.
%\complaint{I re-rephrased this. Paco, Dec 21, 2006} 

In particular, if $p_i$ is an interior vertex, then the mountain edges
must surround $p_i$ in a nonpointed manner.
\end{lemma}

For vertices in general position,  this statement can be rephrased as follows:
``Then, $p_i$ is \rnp\ in the graph consisting of boundary edges of $R$ and
mountain edges of the lift''. We need the more careful statement for
cases were collinear edges arise naturally (see Lemma~\lemRef{downflip} and Figure~\figRef{downflip}b).

\begin{proof}
%Let $p_i$ be a strict local maximum. The statement is equivalent to the following:
%every open half-plane $H$ at $p_i$ that lies inside $R$ in the neighborhood of $p_i$ must
%contain at least one mountain edge incident to $p_i$ that leads into
%the interior of~$H$
%\textup(see Figure~\figRef{maximum}a\textup).
%
%\complaint{Paco, Oct 28: removed what follows, since that figure is not there anymore:
%GR:  reinserted a (new) figure. WE INCLUDE BOTH PROOFS}
\SingleFig{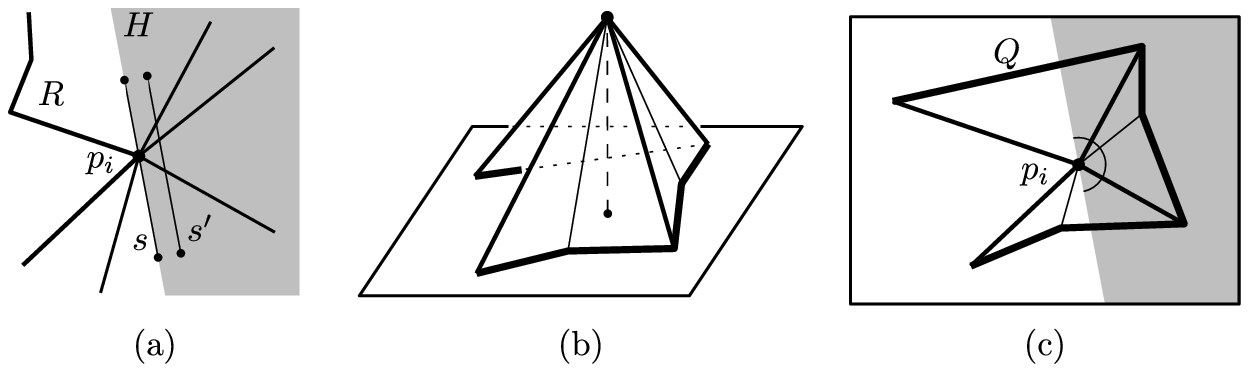}{(a)
At least one of the three edges pointing from the maximum point $p_i$
into $H$ must be a mountain edge.
(b) Cutting the surface below the maximum. Mountain edges are drawn
thicker than valley edges.
(c)~The projected intersection~$Q$.
Mountain edges become convex vertices
and valley edges become reflex vertices of~$Q$.}
%%%%%%
%

Let $s$ be a line segment through $p_i$ on the boundary of $H$.
Consider a segment $s'$ parallel to $s$ that is slightly pushed into $H$.
If $H$ contains no mountain edges incident to $p_i$,
the lifting $f$ must be a convex function on $s'$.
Pushing $s'$ towards $s$, we conclude that $f$ is convex on $s$, and hence
$p_i$ cannot be a strict local maximum.

Another, more visual proof is illustrated in
Figure~\figRef{maximum}b--c.  Let us cut the lifted surface with a
horizontal plane slightly below the maximum point~$p_i$. The
intersection projects to a polygonal chain $Q$ that has a vertex on
every edge incident to~$p_i$.  If follows that $Q$ must contain a
convex vertex in every $180^\circ$ angular range that lies within~$R$.
Since convex vertices of $Q$ result from mountain edges, the lemma is
proved.
\end{proof}

\subsection{Flipping to Local Convexity}
\secLab{locally-convex-flipping}

Let us now come back to the problem discussed at the start of this section: we
have a pointgon $(R,P)$ with given height values $h_i$ at the points $p_i\in
P$, and we look for the highest locally convex function $f$ above $R$ that
does not exceed these heights.  We have seen that such a function is uniquely
defined once we fix a \psT\ $T$ of $(R,P)$. (The given height values at the \rp\
vertices are simply ignored.)  $T$ may not use all interior points of $R$,
i.\,e., it can be a \psT\ of $(R,P')$ with $P'\subseteq P$.
The resulting function will, in general, not be locally convex,
and it may not respect the given heights at the
 \rp\ vertices and at the points of $P-P'$.
But, if $f$ happens to have these properties, it is the solution of our
problem.
\begin{lemma} \lemLab{highest}
  Let $T$ be a \psT\ $(R,P')$ with $P'\subseteq P$  and let $f\colon R\to
  \reals$ be the function that is uniquely defined by the heights of the \rnp\
  vertices of $T$ according to Theorem~\thmRef{liftingbasis}.  If
\begin{equation}
  \eqLab{respect-heights}
  f(p_i) \le h_i\text{, \ for all $p_i \in P$}
\end{equation}
and no interior edge of $T$ is lifted to a mountain edge, then
$f$ is the highest locally convex function that
satisfies \eqRef{respect-heights}.
%\qed
\end{lemma}

\begin{proof}[Proof (Sketch)]
%The proof is based on the fact that
In any locally convex lifting, the heights
$z_i = f(p_i)$ must satisfy \eqRef{z} as an inequality
\begin{equation}
 \eqLab{z-inequality}
z_l\le \lambda_iz_i+\lambda_jz_j+\lambda_kz_k
\end{equation}
if $p_l$ is a vertex of a pseudo-triangle
of $T$ with corners $p_i,p_j,p_k$.
One can show,
by a monotonicity argument similar to the proof of
Theorem~\thmRef{liftingbasis},
 that the highest values $z_i$ that fulfill
\eqRef{z-inequality} for all pseudo-triangles of $T$ and
\eqRef{respect-heights} for all \rnp\ vertices of $T$ must fulfill these
inequalities as equations, and hence they coincide with the function $f$
defined by Theorem~\thmRef{liftingbasis}.
\end{proof}

The following algorithm finds the appropriate \psT\ by a sequence of flipping
operations. We start with an arbitrary triangulation of $(R,P)$. Then all vertices are
\rnp, and \eqRef{respect-heights} is satisfied with equality.
We will maintain \eqRef{respect-heights} throughout.
As long as some interior edge $e$ of $T$ is lifted to a mountain edge,
 we
try to improve the situation by flipping the edge $e$, leading to
another pseudo-triangu\-lation~$T'$.
% It turns out that there is
% always a canonical lifting $F'$ of $T'$ that agrees with $F$ in all
% vertices that are \rnp{} in both $T$ and $T'$. Indeed:
\begin{enumerate}
\item If the flip is a diagonal-edge flip, $T$ and $T'$ have the same set of \rnp{} vertices.
%In particular, $F'$ is unique by Theorem \thmRef{liftingbasis}.

\item If the flip is an edge-removal, $T'$ has one \rnp{} vertex less than $T$.
Theorem \thmRef{liftingbasis} still applies, and the height of the
new \rp{} vertex is derived from the (given) heights of the
\rnp{} vertices.

% \item If the edge inserts an edge, then $T'$ has one more \rnp{} vertex than $T$ and
% $F'$ is not unique in principle, but we choose it as the lifting in
% which this new \rnp{} vertex has the same height it had in $F$.
% Actually, it is clear that in this case $F=F'$ (as a graph). The two
% pseudo-triangles formed by the insertion of the new edge are lifted
% coplanar.

\end{enumerate}

%We are specially interested in the flips of the first two types. In
%both of them there is an edge $e$ of $T$ that disappears.
We remark that in contrast to flipping in triangulations, an edge flip has a
non-local effect on the lifting $f$. The flipping operation affects not only
the faces incident to~$e$, but it modifies the system of equations
\eqRef{z} and may change the heights of other \rp\ vertices.
However, we can predict the direction of this change:
We say that the flip is a \emph{downward flip} if the edge $e$ was a mountain edge. That is,
if $f$ is locally concave in the neighborhood of~$e$.

\begin{lemma}
\lemLab{downflip}
After a downward flip from $T$ to $T'$, the new lifting $f'$
is everywhere weakly below the original lifting $f$. That is, for every  point $x\in R$,
$f'(x)\le f(x)$.
\end{lemma}

\SingleFig{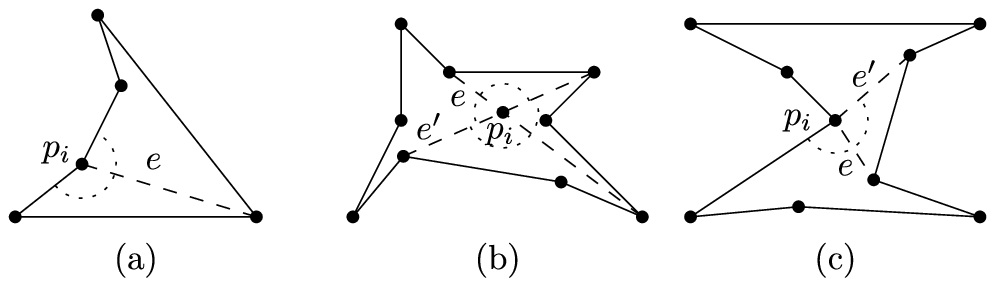}{(a) an edge removal flip;
(b) (c) two types of diagonal flip.}
\begin{proof}
Let $T''$ be the pseudo-triangu\-lation obtained by superimposing
the edges of $T$ and $T'$, see Figure~\figRef{downflip}.
There are several possible cases, but in
all of them $T''$ has exactly one more \rp{} vertex $p_i$ than $T'$:
in an edge removal flip, $T''=T$ and $p_i$ is the end-point  of
$e$ that changes from non-pointed to pointed. In a diagonal flip,
$p_i$ is the intersection of the deleted and inserted edge
(Figure~\figRef{downflip}a). The
intersection is either an interior point of both (then a new vertex
in $T''$, Figure~\figRef{downflip}b) or an end-point of both
(Figure~\figRef{downflip}c), then a pointed vertex of $T$ and
$T'$ that becomes non-pointed in $T''$).

By the definition of $f'$, $p_i$ is the only \rnp{} vertex of $T''$ that
has different height in $f$ and $f'$.

Consider the family of liftings obtained for
varying heights  $h_i$
of $p_i$, keeping the height of every other \rnp{} vertex of $T''$ fixed.
If we set
 $h_i = f(p_i) $
we get the original lifting: the edge $e'$ (if it exists) is lifted
to a flat edge.
For $h_i = f'(p_i) $,
we get the new lifting in which $e$ is flat.
 By the monotonicity stated in
Theorem \thmRef{liftingbasis}, we only  need to prove that
$f(p_i)>f'(p_i)$.

Start with
$h_i$ very high and gradually move
$h_i$ downward.
We know from Lemma~\lemRef{unique-max} that $p_i$ is initially
the unique maximum point.
If we look at the three cases of
Figure~\figRef{downflip} we see that $e$ lies always in some
$180^\circ$ region around $p_i$ with no other edges incident to $p_i$.
Hence, by Lemma~\lemRef{mountain-edge}, $e$ (or the two segments of $e$)
must be a mountain edge.
 The height of every point in the lifting depends linearly on $h_i$,
hence $e$ is a mountain edge for all
 $h_i > f'(p_i) $, it is a flat edge at
 $h_i = f'(p_i) $, and a valley edge below
 $f'(p_i) $.
Since $e$ was assumed to be a mountain edge at
 $h_i = f(p_i) $, we have $f(p_i)>f'(p_i)$.
%
%
% Assume without loss of generality that in the
%initial lifting $f$ the
%
% mountain edge %%!!
%
% that is going to be flipped is
%horizontal and a local maximum with respect to its two incident
%pseudo-triangles. In particular, the maximum of $f$ in these two
%pseudo-triangles is attained at $p_i$. As $p_i$ moves up, no other point can
%move up with higher speed than $p_i$ because the motion is linear and $p_i$ is
%the unique global maximum for big values of $h_i$.
%% (more precisely, for all
%% values bigger than the height of all other \rnp{} vertices).%
%
%In particular,
%for all values $h_i>f(p_i)$,
% $p_i$ is a strict local maximum, and hence the edge $e$ cannot be flat in the
% lifting.
%\complaint{This should be made somewhat more precise in the case of
%Figure~\figRef{downflip}a and c, perhaps another little lemma?}
%This implies that in the lifting $f'$, where $e$ is no longer an edge,
%we must have
%$h_i=f'(p_i) < f(p_i)$, as desired.
\end{proof}

This lemma implies that by performing downward flips
we will eventually arrive at a
locally convex function.

\begin{theorem}[\cite{aichholzer:aurenhammer:brass:krasser:pseudoTriangulationsNovelFlip:2003}]
\thmLab{convexfunction} For any given set of heights $h_i$ of
a pointgon $(R,P)$, and for any initial \psT\ of it, the process of flipping
downwards leads, in a finite number of steps, to a \psT\
$T_0$ and a lifting $f$ which is the highest locally convex function
on $R$ below the values $h_i$.
\end{theorem}

\begin{proof}
% In every  downward flip, the set of \rnp{} vertices either
% stays the same or decreases and the heights of the \rnp{} vertices
% that survive in the flip are kept constant. In particular,
In the
process of flipping downwards no pseudo-triangu\-lation can be
visited twice, by the monotonicity proved in Lemma~\lemRef{downflip}.
% Since the number of pseudo-triangu\-lations
%% of $R$ with the same vertex set as $T$
%is finite,
%% in a finite number of steps
 Eventually,
we must arrive at a pseudo-triangu\-lation $T_0$ whose
associated lifting function $f$ is locally convex.
We started with a triangulation of $(R,P)$, where every vertex of $P$
was lifted to the given height $h_i$.
By construction, $f$ can only decrease in each step, and thus it
is never above the given heights $h_i$.
At the \rnp\ vertices, we have $f(p_i)=h_i$. Hence,
by Lemma~\lemRef{highest}, $f$ is the desired lifting.
%
% To prove that $F_0$ equals the lower envelope, let $F'$ be any other
% locally convex function that lies point-wise below $F$. There is no
% loss of generality in assuming that $F$ and $F'$ agree on the
% \rnp{} vertices of $T$, since decreasing the heights of the
% \rnp{} vertices of $T$ to be those of $F'$ will only decrease the
% lifting $F$ of $T$. In this situation, $F'$ lies clearly below
% $F_0$, because $F_0$ can be characterized as the unique function
% that agrees with $F'$ on certain points (the \rnp{} vertices of
% $T_0$) and that is linear in certain subdomains of $R$ (the
% pseudo-triangles of $T_0$).
\end{proof}

%\complaint{MUST introduce vertex removing flip somewhere;
%Maybe this can come AFTER the above theorem?
%\\
%Done it. Paco.
%}

\SingleFig{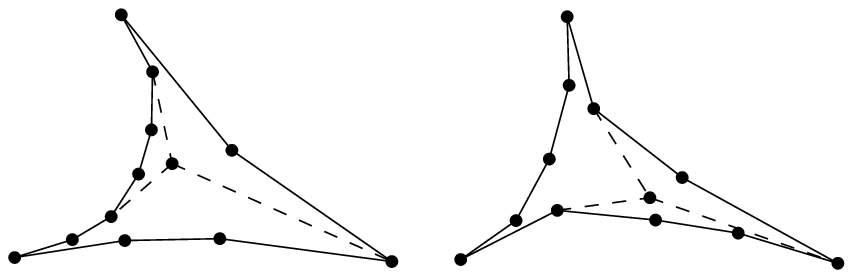}{
Two situations where an edge-removal flip for any
of the three dashed edges will cause the two other edges to become flat.
This is equivalent to removing the common vertex and merging the incident
faces into one.}

The lifting $f$ in the above statement is unique, but the pseudo-triangulation
$T_0$ may not be unique for the following reason: interior vertices may become flat
and will remain so for the rest of the process. Then, the flat union of pseudo-triangles
obtained in the final pseudo-triangulation may be different depending on the initial
pseudo-triangulation, or even on the chosen path of flips.
 To avoid this ambiguity,
we need
 a fourth type of flip that is performed as soon as a vertex gets degree two:
the \emph{vertex removal} flip
%\complaint{AAKB call this ``vertex-removing flip'', not ``deletion
%  flip''. GR changed accordingly. (hopefully everywhere)}
removes this vertex and its two incident edges, merging
two pseudo-triangles into one~%
\cite{aichholzer:aurenhammer:brass:krasser:pseudoTriangulationsNovelFlip:2003}.
%This is called a \emph{vertex deletion flip} (or what was
%the term used in [2]?, coordinate terminology with regular \psT in Chap7).
In particular, if we have a pointed interior vertex of degree~3, an
edge-removal
flip for any of the incident edges will cause this situation to occur,
and we might as well remove the degree-3 vertex right away.
%I will make a FIGURE!
Figure~\figRef{vertexremoval} shows this situation.
The converse of a vertex-removal flip inserts a degree-3 vertex into the
interior of a pseudo-triangular face, connecting it by geodesic paths to the three
corners of that face.
Observe that the result of a vertex removal is a pseudo-triangulation \emph{of a different pointgon},
with the same polygon $R$ but one point less in its interior.

% The last theorem can be read in two ways. On the one hand, it a
% priori gives an algorithm to compute the lower locally convex hull of any
% point set over a region: after triangulating the point set, one just needs
% to  perform downward flips until none is possible. However, it is
% not clear that
It is not known whether the  process of Theorem~\thmRef{convexfunction}
 terminates after a polynomial number of iterations.
% in the number of initial vertices.
In the proof, we have shown that no
pseudo-triangu\-lation can show up twice in the process, but
particular edges can in principle appear and disappear several times.
(An example is given in
\cite{aichholzer:aurenhammer:brass:krasser:pseudoTriangulationsNovelFlip:2003}.)
Only for a convex domain $R$, where we have just triangulations and no
pseudo-triangles, it can be guaranteed that no edge disappears and
reappears, implying a quadratic upper bound on the number of downward
flips needed to get the lower envelope.

%%%%%%%%%%%%%%%%%%%%%%%%%%%%%%%%%%%%%%%%%%%%%%%%%%%%%%%%%%

More significant for us is the fact that the the lower locally
convex hull over a domain $R$ gives rise to a \psT:

\begin{theorem}
[Aichholzer et
al.~\cite{aichholzer:aurenhammer:hackl:preTriangulations:2006}%
]
\thmLab{envelope}
If the points $p_i$ of a pointgon $(R,P)$ and their given heights
$h_i$ are generic, the regions of linearity of the lower locally
convex hull are pseudotriangles, and they form a \psT\ $T$ of
$(R,P')$ for a subset $P'\subseteq P$ of vertices.

Moreover, every interior vertex of $T$ is non-pointed.
\qed
\end{theorem}

\SingleFig{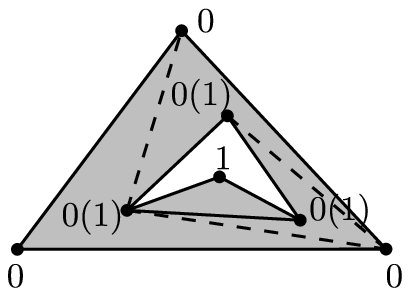}{
In a region~$R$ (shaded) that is not simply connected,
the folding edges of the highest locally convex function may not create a \psT.
The numbers indicate the heights $z_i$ in the lower locally convex hull.
The numbers in parentheses indicate the given heights $h_i$,
whenever they are different from the final heights~$z_i$.
There is a face at $z=0$ which is not a pseudo-triangle. (It is not
even a simple polygon.)
Perturbing the heights or the vertices will not change this situation.}

When the region~$R$ is not a simple polygon, we leave the realm of
\psT s, see Figure~\figRef{non-simply}.  The
characterization of the graphs that can arise (generically) as the
edge sets of locally convex functions (or of general piecewise
linear functions) over such general polygonal domains is still an
open problem.  Some results in this direction are given in
\cite[Section~8]{aichholzer:aurenhammer:hackl:preTriangulations:2006}.

%!TEX root = article.tex

\section{Self-Stresses, Reciprocal Diagrams, and the Maxwell-Cremona Correspondence}
\secLab{self-stresses-0}

\subsection{Maxwell Liftings of Pseudo-Triangulations} % of Point Sets}
\secLab{self-stresses}

A classical result of Maxwell (stated below) relates
three objects for a given geometric graph $G$: its 3D liftings, its
planar reciprocal diagrams, and its equilibrium stresses. Let $G$ be
a connected geometric non-crossing graph.

\emph{Reciprocal diagrams:} A geometric graph $G'$ is called a
\emph{dual} of $G$ if there is an incidence-preserving bijection
from faces and edges  of $G$ to vertices and edges of $G'$,
respectively: $G'$ has a vertex for each face of~$G$. For every edge
of~$G$ that is shared between two faces of~$G$, $G'$ has an edge
between the corresponding vertices.
 $G'$ is called a {\em reciprocal diagram}
of $G$ if each edge of $G$ is parallel to the corresponding edge of
$G'$. A reciprocal diagram $G'$ is not necessarily non-crossing. As
a boundary case we also allow vertices of $G'$ to coalesce. (A
zero-length edge of $G'$ is by definition considered to be always
parallel to the corresponding edge of $G$.)

\emph{Maxwell liftings:} A Maxwell lifting of $G$
is a 3D lifting in the sense of Definition~\defRef{lifting},
where the outer face is lifted to the horizontal plane $z=0$.

\emph{Equilibrium stresses:} Let $P$ and $E$ be the sets of vertices and edges
of $G$. An equilibrium stress (or self-stress) of $G$ is an assignment of a
scalar $\omega_e=\omega_{ij}=\omega_{ji}$ to every edge $e=ij$ of $G$
such that every vertex ``is in equilibrium'':
We think of the edge $ij$ as exerting a force
$\omega_{ij} ({p}_j-{p}_i)$ on the vertex $i$ (and an opposite force on~$j$).
The forces at every vertex $i$ must add up to zero:
\begin{equation}
  \eqLab{self-stress}
\sum_{j|\{i,j\}\in E}  \omega_{ij} ({p}_j-{p}_i)={0}.
\end{equation}
In Section~\secRef{rigidity}, equilibrium stresses will be related to rigidity.

\begin{theorem}[Maxwell~\cite{maxwell:Reciprocal:1864,maxwell:Reciprocal:1870}]
\thmLab{MC}
For every connected geometric non-crossing graph $G$ there is a one-to-one
correspondence \textup(bijection\textup)
between
\begin{enumerate}
\item reciprocal diagrams of $G$
in which the dual vertex of the outer face is at the origin\textup;
% (modulo translation).
\item equilibrium stresses on $G$\textup; and %that are non-zero on every edge.
\item Maxwell liftings of $G$.
%that are not affine-linear on any union of two adjacent faces of $G$ and lift
%the outer face horizontally to height zero.
\end{enumerate}
This bijection is a linear
isomorphism between the corresponding vector spaces.
\qed
\end{theorem}
The equivalence between reciprocal diagrams and equilibrium stresses
is very easy to formulate: From a given reciprocal diagram $G'$ we
associate to each edge $e$ the (signed) quotient between the length
of the edge $e'$ reciprocal to $e$ and the edge $e$ itself. The sign
must be chosen following the following rule (or the opposite one):
consider the edge $e$ oriented from $i$ to $j$, and give the same
orientation to the parallel edge $e'$. If the cell dual to $i$ is to
the right of $e'$ we choose a positive sign for the scalar.
Otherwise, we choose it negative. Conversely, if an equilibrium
stress is given in $G$, the graph $G'$ can be drawn as follows. If
the edges around a vertex $i$ are $e_1,\dots,e_k$, in cyclic order,
the boundary of the cell dual to $i$ in $G'$ consists of the cycle
of edges $\omega_1e_1,\dots, \omega_ke_k$. The equilibrium condition
on $i$ guarantees that this cycle of edges ends at the starting
point. It is not difficult to show that the cells obtained in this
way for the different vertices glue well, and give a non-crossing
graph $G'$. See an example in
Figure~\figRef{wheel4recippseudoeplode}.

\SingleFig{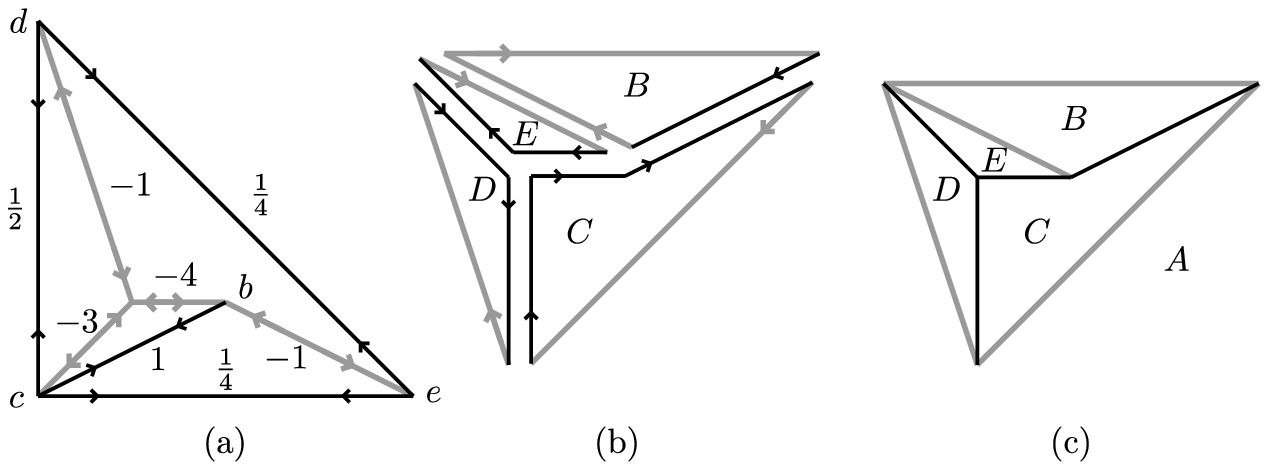}{Assembling the reciprocal.}

The relation between equilibrium stresses and Maxwell liftings is a
bit more difficult to show. The stress associated to a given edge
$e$ is related to the difference between the normal vectors to the
liftings of the two cells incident to $e$. In particular, the sign
of an edge in the stress indicates whether the edge is a mountain or
valley edge in the lifting.
 Edges with stress 0 are lifted to flat edges.

Let us look back at the Surface Theorem (Theorem~\thmRef{liftingbasis}) in the context of
Maxwell liftings, for the case when the boundary of the graph $G$
(i.e., the contour of the outer face) is a convex polygon.
Then the \rnp{} vertices are just
 non-pointed vertices (which lie necessarily in the
interior). Moreover, if $G$ is a pseudo-triangu\-lation then it is
``a pseudo-triangu\-lation of a point set''.

\begin{cor}
Let $T$ be a \psT\ of a point set $P$. Then, for
every choice of height on the non-pointed vertices of $T$ there is a
unique Maxwell lifting of $T$. Moreover, the height of every point
depends linearly and monotonically on the heights of the non-pointed
vertices, and the maximum of  the lifting is achieved either on the
boundary or at a non-pointed vertex, for every choice.
\end{cor}

Specially interesting
is the case of pseudo-triangulations with a unique non-pointed vertex. We call them ``almost-pointed''.

\begin{cor}\corLab{almost-reciprocal}
  \begin{itemize}
  \item [(a)]
 An almost-pointed pseudo-triangu\-lation of a point set has a
unique reciprocal diagram, and a unique Maxwell lifting, modulo
scaling and translation.

  \item [(b)]
 A pointed pseudo-triangulation of a point set has only the trivial Maxwell
lifting,
and it has only the ``degenerate'' reciprocal
diagram, where all vertices coalesce.
\qed
  \end{itemize}
\end{cor}

\subsection{Non-Crossing Graphs with Non-Crossing Reciprocals}
\secLab{reciprocal}

%If one draws an example of an almost-pointed pseudo-triangu\-lation
%and computes its unique reciprocal, the outcome will be another
For any almost-pointed pseudo-triangu\-lation,
the unique reciprocal will be another
almost-pointed pseudo-triangu\-lation. To understand this
phenomenon, Orden et
al.~\cite{orden:roteEtAl:nonCrossingFrameworks:2004} studied the
precise conditions that are sufficient and necessary for a
non-crossing graph with a given stress to produce a non-crossing
reciprocal. Their main result is the following characterization of
when this happens via the type of vertices (pointed or not) and the
sign changes in the equilibrium stress. In the statement, a {\em
sign change} at a face or a vertex is a pair of consecutive edges
(in the cyclic order along the boundary of the face or around the
vertex) whose stress has opposite sign.

\begin{theorem}
  [\protect\cite{orden:roteEtAl:nonCrossingFrameworks:2004}]
{\rm\bf Vertex conditions for a planar reciprocal.}
 \thmLab{planarreciprocal}
    Let $G$ be a non-crossing geometric graph with given self-stress $\omega$.
    Then, in order for
the reciprocal diagram $G'$ to be also non-crossing, the following
\emph{vertex conditions} on its vertex cycles are necessary and
sufficient\textup:
\begin{enumerate}
\item there is a non-pointed vertex with no sign changes.
\item all other vertices are in one of the following three cases:
\begin{enumerate}
\item pointed vertices with two sign changes, none of them at the big angle.
\item pointed vertices with four sign changes, one of them at the big angle.
\item non-pointed vertices
with four sign changes.
\end{enumerate}
\item the face cycles reciprocal to the vertices of type 2.c are themselves non-crossing.
\end{enumerate}
{\rm\bf Face conditions for a planar reciprocal.}
The four types of vertices produce, respectively, the following types
of faces in $G'$\textup:
\begin{enumerate}
\item the \textup(complement of\textup) the exterior face, which is strictly
convex with no sign changes.
\item the internal faces of $G$, which are either
\begin{enumerate}
\item pseudo-triangles with two sign changes, both occurring at corners.
\item pseudo-triangles with four sign changes, three occurring at corners.
\item pseudo-quadrangles with four sign changes, all occurring at corners.
\end{enumerate}
\end{enumerate}
In particular,  for a non-crossing framework to have a non-crossing
reciprocal it is necessary that its faces fall into these four
categories.
\end{theorem}

\begin{proof}[Proof (Sketch)]
The proof of the two sets of conditions is intertwined, and consists of the
following steps. First, necessity of the face conditions is shown by
a local argument: for a given face, the reciprocal vertex can
potentially produce a non-crossing reciprocal only if
 the sum of angles reciprocal to those of the face equals
$2\pi$, and this can be seen to be equivalent to satisfying one of
the face conditions. From this, necessity of the vertex conditions
is also derived, since they are reciprocal to the face
%\complaint{GR replaced vertex by face}
conditions.
Also, a counting argument shows that the vertex conditions (1) and (2)
imply the corresponding face conditions (for the original graph, not
only for the reciprocal). In vertices of types 1, 2.a and 2.b, the
reciprocal face is automatically non-crossing, because it is either
a convex polygon or a pseudo-triangle. However, the vertex condition
2.c can in principle produce a {\em self-intersecting
pseudo-quadrilateral}, see Figure~\figRef{signvsframe}b,
\SingleFig{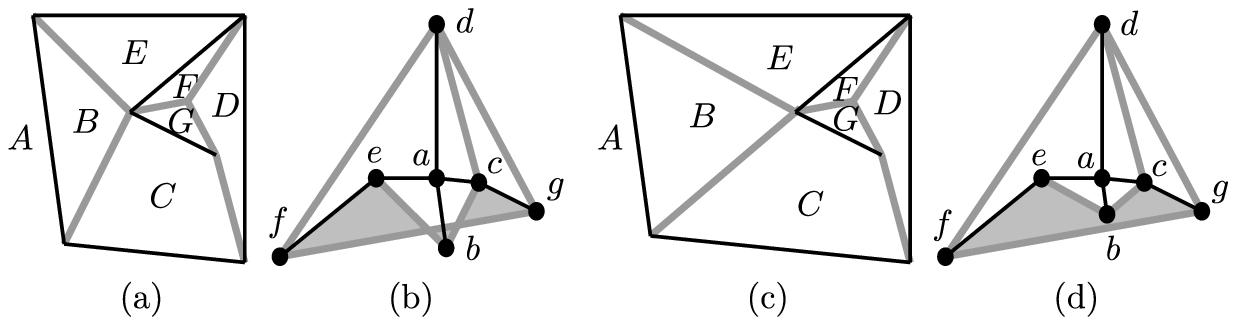} {Sign conditions are not enough to guarantee
a planar reciprocal. The two graphs in (a) and (c) have both the
self-stress with signs represented in the figure by grey and black
edges. The reader can visualize the  Maxwell lifting: the outer face
is at height~0, and if a height is chosen for the common vertex of
faces $B$, $C$ and $E$, then there is a unique Maxwell lifting
compatible with that choice. The signs in the stress indicated
whether the edges are mountain or valley edges in the lifting. The
reciprocal diagrams are shown in (b) and (d). One is non-crossing
but the other is not.}
 but this is
ruled out by condition (3). % in Theorem \thmRef{planarreciprocal}.
Finally, a topological argument shows that if all face cycles
reciprocal to the vertices of $G$ are non-crossing, then the
reciprocal diagram is globally non-crossing, proving sufficiency
of the vertex conditions. % in Theorem \thmRef{planarreciprocal}.
\end{proof}

Let us come back to pseudo-triangu\-lations.

\begin{cor}[\protect\cite{orden:roteEtAl:nonCrossingFrameworks:2004}]
Let $G$ be a pseudo-triangu\-lation with a unique non-pointed
vertex, and let $\omega$ be its unique self-stress. If
$\omega$ is non-zero on every edge, the reciprocal diagram
$G'$ is again an almost-pointed pseudo-triangu\-lation. The
reciprocals of the outer face and non-pointed vertex of $G$ are the
non-pointed vertex and outer face of $G'$.
\end{cor}

\begin{proof}[Proof (Sketch)]
The equilibrium
condition at a pointed vertex implies that around the vertex there
must be at least two sign changes, and at least four if one of them
is at the reflex angle.
A careful accounting of these sign changes proves that these bounds must be exact,
and the possible patterns of sign changes (and sign non-changes) around each
vertex and around each face are exactly equivalent to
conditions 1 and 2 for vertices and for faces in Theorem \thmRef{planarreciprocal}.
\end{proof}

If the self-stress $\omega$ has zeroes,
%? In this case Maxwell
% Theorem still applies, with some care. 
the corresponding zero edge becomes flat
in the corresponding Maxwell lifting (its two incident faces are
coplanar) and it degenerates to length zero in the reciprocal
diagram. % (Its two end-points become identified.)
It turns out that
the previous corollary generalizes to this case:
%, although we do not include a proof of it.

\begin{prop}[\protect\cite{orden:roteEtAl:nonCrossingFrameworks:2004}]
Let $G$ be a pseudo-triangu\-lation with a unique non-pointed
vertex, and let $\omega$ be its unique self-stress. Then, the
reciprocal diagram $G'$ of $G$ is a pseudo-triangu\-lation.

More precisely, let $G^*$ be the subgraph consisting of the edges of
$G$ with non-zero $\omega$. Then, $G^*$ is a pseudo-quadrangulation
\textup(subdivision of a convex polygon into pseudo-triangles and
pseudo-quadrangles\textup). % and $G'$ is a pseudo-triangu\-lation.
% with
%$n-1$ pseudo-triangles and $k+1$ non-pointed vertices, where $n$ and
%$k$ are the numbers of vertices and pseudo-quadrangles in $G^*$,
%respectively.
If $G^*$ has $n$ vertices and $k$ pseudo-quadrangles, $G'$ will 
be a \psT\ with $n-1$
pseudo-triangles and $k+1$ non-pointed vertices.  \qed
\end{prop}
There is a sort of converse of this
statement:
\begin{prop}[\protect\cite{orden:roteEtAl:nonCrossingFrameworks:2004}]
If a pseudo-triangu\-lation $G$ with a non-zero self-stress~$\omega$
produces a non-crossing reciprocal $G'$, then $G'$ can be extended
to an almost-pointed pseudo-triangu\-lation whose unique self-stress
is non-zero exactly at the edges of $G'$.
\qed\end{prop}

We finally mention several interesting properties of the unique,
up to scaling in the vertical direction,
Maxwell lifting of a non-crossing framework with a
non-crossing reciprocal.
We can assume, by changing the sign if necessary, that the lifting
contains points above the zero plane.
These properties are proved in~\cite{orden:roteEtAl:nonCrossingFrameworks:2004}:

\begin{itemize}
\item
 The local and global minima are precisely the points in outer face, including
 the boundary.
\item
 The unique local maximum is the distinguished non-pointed vertex whose reciprocal
is the outer face of $G'$. Moreover, all edges around this vertex
are mountain edges.
\item The maximum is the only point having a supporting plane
which leaves the surface $f$ (locally) on one
side of it.
\item For every height $h$ between the minimum and the maximum, the level curve of $f$ at height
$h$ is a simple cycle. In particular, $f$ has no saddle
points.
\item
 However, in every vertex  $v$ other than the maximum,
the surface is negatively curved in the following sense:
 there is a plane passing through
$v$ that cuts $f$ into 4 pieces in a neighborhood of $v$.
 (In other words, $v$ is a {\em tilted saddle point}.)
\end{itemize}

%!TEX root = article.tex

%------------------------------------------------------------------------------------         Rigidity       ------

%\section{Rigidity of Pseudo-Triangu\-lations}
\section{Pseudo-Triangulations and Rigidity}
\secLab{rigidity}

In this section, we look at geometric graphs as bar-and-joint \emph{frameworks}. The edges are treated as \emph{rigid bars} (they maintain
their lengths) and are allowed to rotate about their incident
vertices (called \emph{joints}). Intuitively, such frameworks are
\emph{flexible} or \emph{rigid}, depending on whether they admit motions
that change their shape or not (precise definitions are given below.)
The results of this section unravel a remarkable behaviour: not only are all
pseudo-triangulations of a point set rigid, but  when a convex hull edge is 
removed from a \emph{pointed} one, it becomes a flexible
\emph{mechanism} with one degree of freedom and which moves expansively
(never decreasing the distance between any pair of joints).
%\complaint{Rephrased a bit. Paco, oct 11}

The proofs rely on concepts and results from rigidity theory
\cite{graver:servatius:servatius:CombinatorialRigidity:1993,whiteley:surveyHandbook:2004,whiteley:Matroids:1996}, which we briefly survey below in Section \secRef{rigidity-prelim}.
%\complaint{Removed ``and involve several subtle technicalities''.
%It unnecessarily scares the reader. Paco, oct 11, 2007}
We start by observing in Section \secRef{generically-rigid} that the
combinatorics is right. As already pointed out in Section
\secRef{basicProperties} (Corollary~\corRef{Laman}), the graphs of pointed pseudo-triangulations
have a specific number and distribution of edges that, for \emph{some}
embedding, guarantees their rigidity, as well as their flexibility
when any edge is removed. Such graphs are known in the Rigidity Theory literature 
under the various technical names: \emph{generic infinitesimally minimally rigid}, {\emph isostatic}, or shortly, \emph{Laman} graphs. The Laman property doesn't guarantee
rigidity in all embeddings, as illustrated in Fig.~\figRef{prism}. In
Section \ref{pt-infRigid}  we prove that the \emph{particular}
pseudo-triangular embeddings are always \emph{infinitesimally rigid}
(defined below), and this will imply their rigidity. When an edge is
removed from a pointed pseudo-triangulation, minimality (with respect
to the number of edges) leads to an infinitesimally
flexible mechanism. A much stronger property, proven in Section
\secRef{pointed-mech}, is that when the removed edge is on the convex
hull, the mechanism moves \emph{expansively}.
% Finally, to guarantee that no undesired branching points are encountered (this may occur for other types of mechanisms), one can prove that the entire expansive phase is free of singularities.

\SingleFig{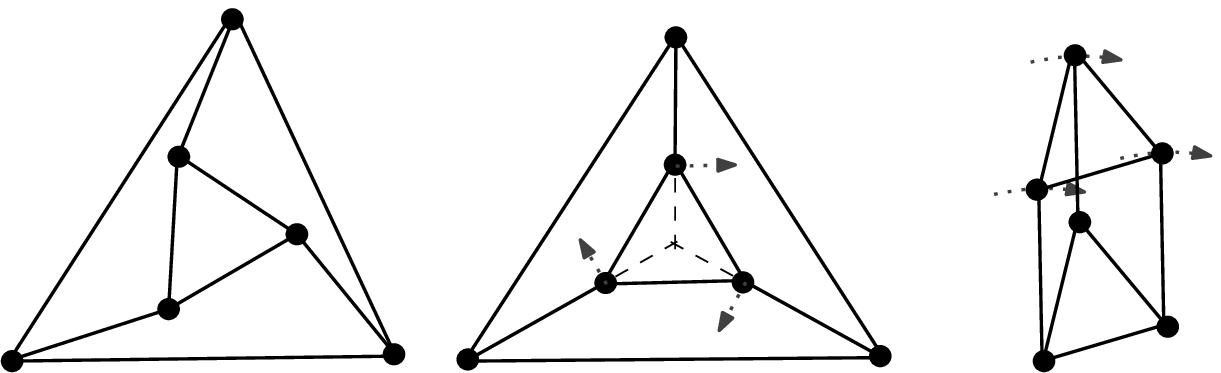}{A minimally rigid graph and three realizations with distinct rigidity properties: (left) a generic framework (infinitesimally rigid, and hence rigid); (center) rigid but infinitesimally flexible; (right) flexible.
}

\subsection{Rigidity of Frameworks}
\secLab{rigidity-prelim} There are three fundamental concepts to be defined:  \emph{rigid}, \emph{infinitesimally rigid} and \emph{generically infinitesimally rigid}.

\subsubsection*{Rigid and flexible frameworks} 
A geometric
graph $G$ embedded on the point set $P=\{p_1,\ldots,p_n\}$, denoted by
$\GP$, is called a \emph{framework} throughout this section,
in accordance with standard usage in rigidity theory.
%\complaint{Added this to justify our use of a new name for an object
 % that already had one (geometric graph). Paco, oct 11. modified GR}
The collection of all realizations $G(P')$ producing the same edge
lengths $\{\ell_{ij}, ij\in E\}$ as $\GP$ is called the \emph{configuration space} of
the framework. It is an algebraic subset of $\reals^{2\nrVert}$, 
described  as the set of all {real}
solutions of the quadratic system
\[
(x_i-x_j)^2 + (y_i-y_j)^2 = \ell_{ij}^2, \ ij \in E,
\]
where the pair of unknowns $(x_i,y_i)$ denotes the position of
the $i$th vertex of $P'$. Since an isometry applied to any solution leads to another solution, 
it is customary to fix (``pin down'') an arbitrary edge of $G$ and to consider the configuration space of the pinned framework. All throughout, we will assume without loss of generality that this is the edge between vertices $1$ and $2$. Algebraically, this is done by adding
the equations $(x_1, y_1)=p_1$ and $(x_2, y_2)=p_2$ to the quadratic system and it results in a reduction of the dimension of the configuration space by three (not four, since it makes the equation between vertices 1 and 2 redundant).

A framework $\GP$ is \emph{rigid} if its pinned version is an isolated point in its configuration 
space. Otherwise, it is \emph{flexible}.

\subsubsection*{Continuous flexes}
A \emph{flex} or \emph{reconfiguration} of a framework $\GP$ is a continuous curve in its configuration space: a function $P(t)=(p_1(t), \ldots, p_n(t))$ defined over some interval of time, satisfying $P(0) = P$ and
\begin{equation}
\eqLab{length-fixed}
 \|p_i(t)-p_j(t)\| = \ell_{ij}, \forall ij\in E
\end{equation}
We avoid trivial flexes by assuming that a certain edge
is fixed: $p_1(t)=p_1(0)$ and $p_2(t)=p_2(0)$, for all $t$. 

\subsubsection*{Infinitesimal rigidity} Classical results in real algebraic geometry imply that any flexible framework $\GP$ admits a differentiable flex $P(t)$. Taking the derivative with respect to the time parameter $t$ in the equations \eqRef{length-fixed}, and denoting by $v_i :=\dot p_i(0)$, we obtain:
\begin{equation}
  \eqLab{rigid}
\langle p_i-p_j, v_i-v_j \rangle = 0, \forall ij\in E
\end{equation}

An \emph{infinitesimal motion} or \emph{infinitesimal flex} of a framework $\GP$
is a family of \emph{velocity vectors} $v=(v_1, \ldots,
v_n)$, $v_i\in \reals^2$ which satisfy the equations \eqRef{rigid}; we
say that $v$ preserve the lengths of all edges $ij\in E$,
infinitesimally. As with continuous flexes, we pin down an edge (by
setting $v_1 = v_2=0$) to exclude the 3-dimensional space of \emph{trivial infinitesimal
  motions}: infinitesimal translations and rotations.
 
A framework $\GP$ is \emph{infinitesimally rigid} if it has no
(non-trivial) infinitesimal motion, and \emph{infinitesimally flexible} otherwise. 

It is known that infinitesimally rigid frameworks
are rigid~\cite{asimow:roth:rigidityOfGraphs:1978}. Intuitively, if they were not, the existence of a 
continuous motion would imply the existence of a differentiable, and hence of an infinitesimal motion. The converse is not true, as illustrated by the second embedding in Fig.~\figRef{prism}.
The framework is rigid, but the interior triangle can be ``infinitesimally rotated''
with respect to the exterior one.

\subsubsection*{The rigidity matrix} Infinitesimal rigidity can be expressed in matrix form. The $2n\times m$ coefficient
matrix $M$ of the system of equations \eqRef{rigid} is
called the \emph{rigidity matrix} associated to the framework
$\GP$. The set of infinitesimal motions is the kernel % $ker(M)$
of $M$, and thus a linear subspace in $(\reals^2)^n \simeq
\reals^{2n}$. This set always contains the three-dimensional linear subspace of
\emph{trivial motions} (infinitesimal rotations and translations). We may reformulate
infinitesimally rigidity as being equivalent to
the kernel of its rigidity matrix having
dimension exactly three. 

We say that $\GP$ has $d$ \emph{\textup(internal\textup) infinitesimal degrees of
freedom} if the space of infinitesimal motions has dimension $d+3$.

\subsubsection*{Generic rigidity}
Different frameworks realizing the same abstract graph 
can have different rigidity properties, as illustrated in
Fig.~\figRef{prism}. A framework $\GP$ is called \emph{generic} if its
rigidity matrix has the maximum rank among all possible
embeddings~$P'$. Observe that the non-generic embeddings form an
algebraic subset (some minors of the rigidity matrix become zero), so generic embeddings form an open dense subset in the set of all embeddings.

In particular, an infinitesimally rigid framework whose rigidity matrix
has  the maximum possible rank $2n-3$ is generic.
An abstract graph $G$ is called \emph{generically rigid}, or shortly \emph{rigid}, if there
exists a set of points on which the framework $\GP$ is
infinitesimally rigid. Equivalently, such graphs are infinitesimally
rigid in \emph{almost all} 
 realizations. 
An abstract graph is \emph{minimally rigid}
if it is rigid but the removal of \emph{any} edge invalidates this
property.

The fundamental theorem of rigidity theory in the plane can
now be stated:

\begin{theorem} [Laman \cite{laman:Rigidity:1970}]
	\label{thm:laman}
A graph is minimally rigid iff it is a Laman graph.
\end{theorem}

\subsection{Pointed Pseudo-Triangu\-lation Graphs are Minimally Rigid}
\secLab{generically-rigid}

By combining Theorem \ref{thm:laman} with Corollary~\corRef{Laman} in Section \secRef{pointedPT}, 
we obtain the first fundamental rigidity property of pseudo-triangu\-lations:

\begin{cor}
\corLab{ppts-genericallyrigid}
The underlying graph of a pointed pseudo-triangulation of a point set is 
minimally rigid.
\end{cor}

Next, we proceed to show that not only is the underlying graph of a pointed pseudo-triangulation special, but so is the embedding.

\subsection{Pseudo-Triangulations are Infinitesimally Rigid}
\label{pt-infRigid}
The following classical result relates infinitesimal rigidity of a
framework to the equilibrium stresses introduced in
Section~\secRef{self-stresses}.

\begin{theorem}
\thmLab{infinitesimal} Let $\GP$ be a framework in the
plane with $\nrVert$ vertices, $\nrEdges$ edges, $d$
{infinitesimal degrees of freedom} and
a space of self-stresses of dimension $s$. Then\textup:
\[
\nrEdges = 2\nrVert -3 + (s - d)
\]
\end{theorem}

\begin{proof}
Compare the system~\eqRef{rigid} with the system~\eqRef{self-stress}
of Section~\secRef{self-stresses}: the equilibrium stresses form the kernel 
of the transpose of the $2n\times m$ rigidity matrix~$M$. Let the rank of $M$ 
be $r$. Since $d+3$ and $s$ are the kernel dimensions of $M$ and 
its transpose, we have
$$
r= \nrEdges - s =2 \nrVert - (d+3)
\eqno\qedhere
$$
\end{proof}

The following consequence of this theorem was first proved by Streinu~\cite{streinu:pseudoTriangRigidityMotionPlanning-confAndJour:2005}
for the pointed case and by Orden and Santos~\cite{orden:santos:polytope:2005} in general.
\begin{theorem}
\thmLab{pt-rigid} % \corLab{pt-rigid}
 Every pseudo-triangulation of a
point set is infinitesimally rigid, hence rigid. Pointed
pseudo-triangulations are minimally infinitesimally rigid.
\end{theorem}

\begin{proof}
Our proof is based on the
 Surface Theorem (our Theorem~\thmRef{liftingbasis})
of Aichholzer et al.~\cite{aichholzer:aurenhammer:brass:krasser:pseudoTriangulationsNovelFlip:2003}.
By Corollary~\corRef{lifts-dimension}, together with Maxwell's
Theorem~\thmRef{MC}, the dimension $s$ of the space of equilibrium
stresses of a pseudo-triangulation equals its number $\nrNonPointed$
of non-pointed vertices. Since a pseudo-triangulation on $\nrVert$
vertices has $2\nrVert - 3 + \nrNonPointed$ edges, the previous
theorem implies that it has no non-trivial infinitesimal motions.
\end{proof}

In particular, the pointed case of Theorem~\thmRef{pt-rigid} implies Corollary~\corRef{ppts-genericallyrigid} (and thus provides a different proof of it).
However, for non-pointed pseudo-triangulations we do not have a direct, combinatorial proof of their generic rigidity (that is, of the fact that they contain a spanning Laman subgraph), other than via their infinitesimal rigidity: not every pseudo-triangulation contains a pointed pseudo-triangulation.

\subsection{Pointed Pseudo-Triangu\-lation Mechanisms}
\secLab{pointed-mech}

If a framework is infinitesimally minimally rigid, the removal of
any edge creates a flexible object with one infinitesimal degree of
freedom. In this section we outline an even stronger property
\cite{streinu:pseudoTriangRigidityMotionPlanning-confAndJour:2005},
for the case of a pointed pseudo-triangu\-lation with a convex hull
edge removed: the resulting framework is a mechanism (with one
degree of freedom) which \emph{expands} all distances between its
vertices when the endpoints of the removed convex hull edge are
moved away from each other.
We call these frameworks \emph{\ppsT\ mechanisms}.

\subsubsection*{Infinitesimal expansion and contraction}
For a point set $\pts$ and infinitesimal velocities $v$, we define
the (infinitesimal) \emph{expansion} $\varepsilon_{ij}$ as:
\begin{equation}
  \eqLab{expansion-def}
\varepsilon_{ij} %(p,v)
:=
\langle p_i-p_j, v_i-v_j \rangle
\end{equation}

We say that (infinitesimally) the pair $ij$ of points \emph{expands}
if $\varepsilon_{ij}\geq0$ and \emph{contracts} if
$\varepsilon_{ij}\leq0$. An infinitesimal motion is \emph{expansive}
if all pairs of vertices expand.

\subsubsection*{Expansive mechanisms} 
An (infinitesimal) mechanism is a framework
$\GP$ with a one-dimen\-sion\-al space of infinitesimal motions.
Intuitively, up to reversal, there is only one direction in which  the mechanism can move.
The mechanism is (infinitesimally) \emph{expansive} if,
in this motion (or its reverse), all pairs of vertices expand.

The following lemma results from the fact that the rigidity matrix
is the transpose of the matrix of the system \eqRef{self-stress}
defining self-stresses. The space of expansions $\eps_{ij}$ is the
image of this matrix, and is thus orthogonal to the kernel of the
transpose (cf.~the proof of Theorem~\thmRef{infinitesimal}).

\begin{lemma}
\lemLab{farkasExpansive}
Let $H$ be a Laman framework with an extra edge,
and let $ij$ and $kl$ be two edges of~$H$.
Let $H'=H-\{ij,kl\}$.
 Let $v$ be an
infinitesimal motion of $H'$, and let $\omega$ be a self-stress of
$H$. Then,
$$%\begin{equation}
  %\eqLab{ooo}
\omega_{ij}\eps_{ij} + \omega_{kl}\eps_{kl} = 0
\eqno{\qed}
$$%\end{equation}
\end{lemma}

We can now prove
the main result of this section:
\begin{theorem} %\ \ {\bf (Infinitesimal expansiveness of pseudo-triangu\-lation mechanisms,
[Streinu
\cite{streinu:pseudoTriangRigidityMotionPlanning-confAndJour:2005}]
\thmLab{infMechanism}
A  bar-and-joint framework whose underlying graph is obtained by
removing a convex hull edge from a \ppsT\ is an
expansive mechanism.
\end{theorem}

\begin{proof} Let $G$ be the underlying graph of the pointed
  pseudo-triangu\-lation, $ij$ be the removed convex hull edge (so that $G
  \setminus \{ij\}$ is a \ppsT\ mechanism).
 Let $v$ be an infinitesimal motion
  preserving the edge lengths of $G\setminus\{ij\}$ and increasing the length
  of the edge $ij$, $\eps_{ij}>0$. We want to prove that it also increases (infinitesimally) the
  distance of every other pair of vertices $kl$, i.\,e.,
  $\eps_{kl}\ge0$.

Consider the graph $H=G\cup\{kl\}$. It has $2n-2$ edges, and
therefore, by a dimension argument, it supports a non-trivial
self-stress $\omega$. Since $G$ itself supports no self-stress, we
must have $\omega_{kl}\ne 0$. Without loss of generality, we assume
$\omega_{kl}> 0$. Now, if we have
 $\eps_{kl}<0$, then % \eqRef{ooo}
Lemma~\lemRef{farkasExpansive}
 would imply
$\omega_{ij}> 0$. Hence, it is sufficient to show that $\omega_{ij}>
0$ and $\omega_{kl}> 0$ leads to a contradiction. We will interpret
the self-stress in terms of the
  induced Maxwell lifting to derive this contradiction.

Consider the framework $G\cup \{kl\}$ obtained from $G$ by adding the
extra edge $e=kl$. Since it is no longer a pointed
pseudo-triangu\-lation, either one endpoint or both endpoints of this new edge $kl$ are
non-pointed, or else the edge $kl$ crosses some other edges
of the pseudo-triangu\-lation. If new crossings have been
introduced, we will apply Bow's construction \cite{maxwell:Bow:1876}
to eliminate them: simply replace each crossing by a new
vertex, and split the two crossing edges in two. We obtain a new
planar framework $G'$ on $n'=n+1$ vertices and
$2n'-2$ edges. In either case, the new framework is non-pointed only
at one or both of the endpoints of $e$ or at the crossings of $e$
with other edges (or both). Denote by $P_\nnonpointed$ the non-empty set of non-pointed
vertices: the endpoints of $e$, if non-pointed, and the crossings of
$e$ with other edges, if there are such crossings.

The signs of the self-stresses are preserved by Bow's construction.
Since we assumed a strictly positive self-stress on $ij$ and $e=kl$,
this means that both $ij$ and $e$ are valley edges in a Maxwell
lifting of $G'$ (and this extends to the split edges in case
 Bow's construction has been applied).
 The only
edges which could be mountain edges are the edges of $G$.

Since the convex hull edge $ij$ is a valley edge, the Maxwell lifting contains
points above the outer face $z=0$, so the maximum height be attained
on some set $M$ consisting of vertices, edges, and bounded faces of $G'$.

Let us look at a convex hull vertex $p$ of $M$.  It follows from
Lemma~\lemRef{lifting-basic} that $p$ must be a non-pointed vertex
from $P_\nnonpointed$. Since all vertices of $P_\nnonpointed$ lie on
$e$, we know that $M$ is a union of vertices and edges that lie
above~$e$. The edges lying above $e$ are either $e$ or splittings of
$e$ and have negative stress, therefore lift to valleys. Valleys
cannot be maxima, obviously. Thus the maximum height is attained at
a set $M$ of isolated interior vertices.

To complete the proof, let us focus on one such vertex $p_i$ of maximum
height. %
By
Lemma~\lemRef{mountain-edge},
% every $180^\circ$ region around $M$
%contains at least
the mountain edges, which must be edges of $G$,
surround $p_i$ in a non-pointed manner, contradicting the
assumption that $G$ is pointed.
\end{proof}

Classical considerations from differential equations are further used
in \cite{streinu:pseudoTriangRigidityMotionPlanning-confAndJour:2005}
to show that this infinitesimal motion leads to an actual motion: a
one-dimensional trajectory in configuration space, for a portion of
which all vertices move away from each other. Moreover, the expansion
portion of the trajectory is free of singularities, implying that, in
contrast to what may occur for other types of mechanisms, there are no
branching points along the way: the mechanism is guaranteed to move in
a well-defined direction until it ceases to be expansive. In
Section~\secRef{carpenter}, this property is used for an algorithmic
solution to the Carpenter's Rule Problem.

\begin{theorem} [Motion of pointed pseudo-triangu\-lation mechanisms
\cite{streinu:pseudoTriangRigidityMotionPlanning-confAndJour:2005}]
\thmLab{mainPropPPTMech}
A pointed pseudo-triangu\-lation mechanism moves expansively on a
unique
trajectory, from the moment when a corner angle is zero to
the moment when two extreme edges of a vertex \textup(or, as a special
case, when one of these is the missing convex hull edge\textup) align.
%One of these two aligned edges may be the
%missing convex hull edge.
\qed
\end{theorem}

\subsection{Parallel Redrawings of Pointed Pseudo-Triangu\-lation Mechanisms}
\secLab{parallelPPT}

Two-dimensional rigidity has an alternative, equivalent
model, called the theory of \emph{parallel redrawings}, in which
edges of frameworks are required to maintain their slopes instead of
their lengths
\cite{whiteley:surveyHandbook:2004,whiteley:Matroids:1996}. It turns
out that pointed pseudo-triangulation mechanisms are also special
from the point of view of the theory of parallel redrawings:
all of their parallel redrawings are non-crossing.

A \emph{parallel redrawing} of a geometric graph is a redrawing 
of the underlying abstract graph where
corresponding edges have the same slopes. Any translation or rescaling
is a parallel redrawing, but this is not interesting: we call it
\emph{trivial}. Parallel redrawings are closely related to the
infinitesimal motions of bar-and-joint frameworks:
a family of velocity vectors $w=(w_1, \ldots,
w_n)$, $w_i\in \reals^2$ preserves the direction of an edge $ij\in E$
if $w_i$ and $w_j$ have the same projection on the direction
\emph{perpendicular} to the direction $p_i-p_j$:
\begin{equation}
  \eqLab{parallel-redraw}
\langle w_i-w_j, (p_i-p_j)^\perp \rangle = 0,\ \forall ij\in E
\end{equation}
where $p^\perp$ denotes the counterclockwise rotation of the vector
$p$ by $90^\circ$.
Except for this rotation, the conditions
\eqRef{parallel-redraw} are the same as the equations
\eqRef{rigid} for infinitesimal motions.
(In contrast to fixed-edge rigidity, the conditions
\eqRef{parallel-redraw} preserve directions exactly,
not just infinitesimally.)
Thus, the spaces of parallel redrawings and infinitesimal motions have
the same dimensions, and
$v$ is an infinitesimal motion if and only if
$w = v^\perp$ is a parallel redrawing, where
$v^\perp$ denotes the family of vectors $(v_1^\perp,\ldots,
v_n^\perp)$.
We eliminate trivial parallel redrawings by pinning down a vertex and
rescaling; this turns the space of parallel redrawings of any generic
one-degree-of-freedom mechanism into a 1-dimensional projective space.
Its double-covering gives a circular configuration space of parallel
redrawings. With an appropriate parametrization, this can be
visualized as a continuous circular parallel redrawing motion,
illustrated in Figure~\figRef{snapshots-new} for a pointed
pseudo-triangular mechanism.
%

%%%%%%
\SingleFig{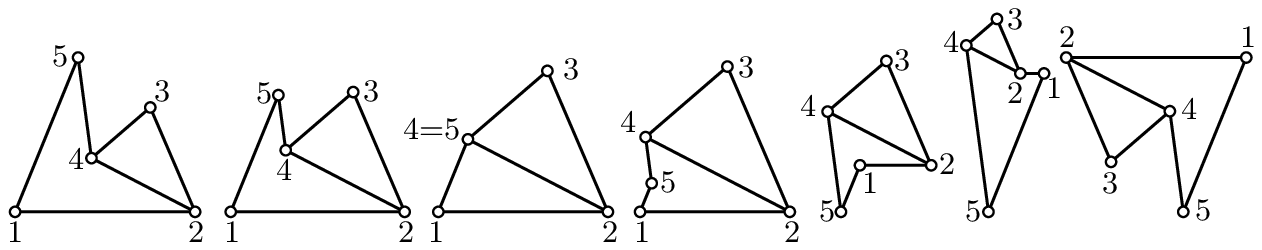}{A few
snapshots of a cyclic trajectory running through all
parallel redrawings of a
pointed pseudo-triangulation mechanism. The first two
snapshots have the same combinatorial structure.
At certain discrete times, like in the third picture,
some edges are reduced to zero length, and
some angle is about to change from
convex to reflex and vice versa.
Each remaining snapshot shows the new combinatorial
structure after a rigid component (possibly just a single edge)
has shrunk to a point.
 The last snapshot is a rotated copy of
the first one, completing half of the cyclic sequence of all parallel redrawings.}

%In particular, we can conclude that pointed pseudo-triangulation
%mechanisms are \emph{flexible} from the point of view of parallel
%redrawings, i.e. they have non-trivial parallel redrawings.  Figure
%\figRef{snapshots-new} illustrates a sequence of pseudo-triangulation
%parallel redrawings.
%
The following result summarizes the properties of pointed
pseudo-triangulation mechanisms, as parallel-redrawing mechanisms.

\begin{theorem} [Streinu
\cite{streinu:parallelPseudoTriangulations:2005}]
\thmLab{pptParallel}
A pointed pseudo-triangulation with a convex hull edge removed has a
one-dimensional projective space of non-trivial parallel redrawings.
All of them are non-crossing,
and, except when some edges are reduced to zero length,
they are
 pointed pseudo-triangulation mechanisms with the
same plane graph facial structure.
%\qed
\end{theorem}
\begin{proof}[Proof (Sketch)]
We only prove here that the redrawings
 are pointed pseudo-triangulation mechanisms, and in particular, that
 they are non-crossing.
The proof nicely connects expansive motions and parallel redrawings.
%We only prove here non-crossing for the parallel redrawings in an appropriately chosen affine part of the projective configuration space. 

Let us pin down a convenient edge, say between vertices $1$ and $2$, by setting
$v_1=v_2=0$ for infinitesimal motions and
$w_1=w_2=0$ for parallel redrawings,
 to exclude trivial parallel redrawings.
%
%; this excludes not only
%all the trivial parallel redrawings, but also the special realization where the pinned edge shrinks to zero-length. This excluded realization gives the velocities by which the vertices move on linear trajectories in the affine part of the projective configuration space, in a process called a {\em parallel redrawing sweep}. 
%
% With these constraints, there is only one non-trivial solution
% $v=(v_1,\ldots,v_n)$ of the system \eqRef{rigid} (infinitesimal
% motion), up to a scalar multiple, by Theorem~\thmRef{infinitesimal},
% and moreover, by Theorem~\thmRef{infMechanism}, this motion is
% expansive:
With these constraints, the space of solutions of the system
\eqRef{rigid} (infinitesimal motions) is one-dimensional by
Theorem~\thmRef{infinitesimal}. Interpreted as parallel redrawings, this gives the (one dimensional) space of realizations with the same edge slopes. Moreover, by
Theorem~\thmRef{infMechanism}, we have a non-trivial solution
$v=(v_1,\ldots,v_n)$ which is an expansive infinitesimal motion:
\begin{equation*}
%  \eqLab{exp0}
\langle v_i-v_j, p_i-p_j \rangle \ge 0,
\end{equation*}
for all $1\le i,j\le n$.
The one-dimensional affine space of parallel redrawings is %can be shown to be
given by $(p_i+t\cdot v_i^\perp)$. As the parameter $t\in \reals$
varies, the vertices move on linear trajectories, in a process called
the \emph{parallel redrawing sweep}.
A straightforward calculation shows that the same velocity vector
$(v_i)$ is an expansive infinitesimal motion for \emph{all} these
parallel redrawings:
\begin{equation*}
\begin{split} 
\langle v_i-v_j, p_i+t\cdot v_i^\perp
-(p_j+t\cdot v_j^\perp) \rangle
=
\langle v_i-v_j, p_i-p_j \rangle
+
t\cdot \langle v_i-v_j,  v_i^\perp
- v_j^\perp \rangle
\\
=
\langle v_i-v_j, p_i-p_j \rangle
+
t\cdot0
=
\langle v_i-v_j, p_i-p_j \rangle,
\end{split}
\end{equation*}
which is 0 for all edges of the graph and non-negative for all pairs
$i,j$,
by assumption.

Thus we know that every parallel redrawing has a non-trivial expansive
motion.  From this we need to conclude that it is
non-crossing.
Graphs whose space of infinitesimal motions is one-dimensional and
contains an expansive motion are also studied below in
Section~\secRef{expansion-cone}:
Theorem~\thmRef{expansion-cone} provides a converse of
Theorem~\thmRef{infMechanism}, and it
implies that such a graph
is a pointed \psT\ mechanism, except that
some \emph{rigid components}
 may have been replaced by other minimally rigid
frameworks on the same point set
(see Fig.~\figRef{collapse} for an example).
To show that these components are in fact non-crossing and pointed,
one uses that fact that rigid components are preserved in parallel
redrawings: they can only be scaled as a whole.
\end{proof}

%!TEX root = article.tex

%-----------------------------------------------------------------------------------------------     Drawing PPT     ------

%\section{Generically Rigid Graphs are Pseudo-Triangu\-lation Graphs}
\section{Planar Rigid Graphs are Pseudo-Triangu\-lations}
\secLab{planarLaman}
%

%\complaint{I am thinking in changing the title of this section to
%  ``Pseudo-triangulability of rigid graphs''. The idea is to
%  ``reverse'' the title of Section 6. What fo you think? Paco, October
%  6, 2007.
%GR says:
%I would not go for such comlicated words like PSTRBLITY.
%If you want uniformity, then rather invert the title of Section 6:
%"PSTRI are rigid."
%Or we could shorten 7 to
%"Rigid Graphs are Ps-triangulations." (a bit imprecise,
%but a section title is not a statement of a thm. }

%We have seen in Theorem~\thmRef{pt-rigid} that pseudo-triangulations of a point set are infinitesimally rigid. In particular, their underlying (abstract) graphs are generically rigid in the plane.
%\complaint{This sentence may need to change, to adapt to the final
%version of Section 6. Paco, October 27, 2006}
In this section we will mainly discuss the following result:

\begin{theorem}
   \thmLab{planarrigid} Every planar rigid \textup(abstract\textup) graph
  can be embedded as a
  pseudo-triangulation.
\end{theorem}

We prove in fact an even stronger statement: every topological plane embedding of $G$ can be
``stretched'' to become a pseudo-triangu\-lation of a point set. A
\emph{topological plane graph} retains only the facial structure of 
a planar graph embedding.
\complaint{Ileana: I rephrased.}

This is a converse to Corollary~\corRef{ppts-genericallyrigid}. It is also
reminiscent of Tutte's theorem that every simple planar 3-connected
graph can be drawn in the plane with convex faces \cite{t-crg-60,t-hdg-63}.
Instead of convex faces, here we use pseudo-triangles, which are {\em as
non-convex as possible}.  In fact, one of the steps in the proof uses
a variation of Tutte's barycentric method for convex drawings~\cite{t-hdg-63}.
 
%reminiscent of Steinitz's Theorem that every simple planar and 3-connected
%abstract graph is the graph of a 3-polytope, hence it can be drawn in the plane
%with convex faces and convex outer face. Instead of convex faces, here we use
%pseudo-triangles, that are ``as non-convex as possible''.
%In fact, one of the steps in the proof uses a variation of Tutte's barycentric method, which, in its original form, has also been
%used for proving Steinitz's Theorem.

% In fact, one of the steps in the proof uses Tutte's barycentric method, originally
% devised for a new proof of  Steinitz's Theorem.
  
%\complaint{changed this paragraph. Paco d'apres Ileana, dec 20, 2006}
%
%With not much extra work we get the following statement, taken from~\cite{orden:santos:servatiusB:servatiusH:combinatorialPT:2006}:

%\begin{cor}
%%
%\corLab{planarrigid}
%%
%For a planar graph $G$ the following conditions are equivalent:
%\begin{enumerate}
%\item $G$ is generically rigid.
%\item $G$ can be embedded as a
%pseudo-triangu\-lation of a point set.
%\item Every topological
%embedding of $G$ can be stretched to become a pseudo-triangu\-lation
%of a point set.
%\end{enumerate}
%Moreover, the number of non-pointed vertices in the embeddings of parts (b) and (c)
%equals the generic dimension of the space of equilibrium stresses of the graph.
%\end{cor}

%\begin{proof}
%The implication from (3) to (2) in is obvious and the implication (2)
%to (1) was proved in Corollary~\corRef{pt-rigid}. The ``moreover'' holds since it is true
%for every pseudo-triangulation. The implication from (1) to (3) 
%is Theorem~\thmRef{planarrigid}.
%\end{proof}

In Sections~\secRef{CPT}--\secRef{cpt-gen-laman} we sketch the proof
of Theorem~\thmRef{planarrigid}. Most of the details of the first two of these three sections are
contained
in~\cite{streinu:haas:etAl:planarMinRigidPseudoTriang-confAndJour:2005},
although the theorem is only proved in full generality
in~\cite{orden:santos:servatiusB:servatiusH:combinatorialPT:2006}.
% Actually,~\cite{streinu:haas:etAl:planarMinRigidPseudoTriang-confAndJour:2005}
In~\cite{streinu:haas:etAl:planarMinRigidPseudoTriang-confAndJour:2005}, 
two different methods are presented, but only for the case of Laman graphs.
We briefly
discuss them in Section~\secRef{cpt-matchings}.

Section~\secRef{graph-drawing} shows the following general result relating the ``level of rigidity'' of a planar (abstract) graph and the ``level of pointedness'' of its embeddings.

\begin{theorem}
  \thmLab{dof-corners} Let $G$ be a connected planar graph with $d$
  generic degrees of freedom (the space of
  infinitesimal motions of any generic embedding of $G$
  in the plane is $d$-dimensional). Then:
\begin{enumerate}
\item In any non-crossing straight-line embedding of $G$, the ``excess
  of corners'' $\bar{\nrCorners}$ \textup(as defined in
  Section~\secRef{pointedPT}\textup) is greater or equal to $d$.
\item $G$ has non-crossing straight-line embeddings in which $\bar{\nrCorners}=d$.
\end{enumerate}
\end{theorem}

This follows easily from Theorems~\thmRef{pt-rigid} and~\thmRef{planarrigid},
which are the ``base cases'' $\bar{\nrCorners}=0$ and $d=0$, respectively.
But this generalization has not been previously published.

\subsection{Combinatorial Pseudo-Triangulations}
\secLab{CPT}

%The proof of Theorem~\thmRef{planarrigid} relies on the concept of
%combinatorial pseudo-triangu\-lation.

\begin{definition}
\defLab{CPT}
% \footnote{The definition is valid in a more
%general setting than what we use in this paper, and works even
%with multiple edges, vertices of degree one and loops. In such a
%case, a vertex of degree one is  incident to a unique  angle,
%labelled {\em big}.}.
A {\it combinatorial pseudo-triangu\-lation \textup(CPT\textup)} on  a plane graph $G$ is an
assignment of labels {\it big} (or {\it reflex}) and {\it small} (or
{\it convex}) to the angles (vertex-face incidences) of $G$
such that:
\begin{itemize}
\item[(i)] Every face except the outer face gets exactly {\em
three}
vertices marked {\it small}. These will be called the {\em corners} of the face.%
\item[(ii)] The outer face gets only {\it big} labels (it has no corners).%
\item[(iii)] Each vertex is incident to at most one angle labeled
{\it big}. The vertices incident to big angles are called {\em pointed}. %
%\item[(iv)] Every vertex of degree smaller than three is pointed.
%
%\complaint{I think we decided to drop condition (iv), and have it be
%part of ``generalized Laman''}
\end{itemize}
\end{definition}

By analogy with pseudo-triangu\-lations, we also define
\emph{non-pointed vertices}, 
\emph{extreme edges} of a pointed vertex, corners and \emph{pseudo-edges} of a pseudo-triangle, etc.
CPT's behave very much like true pseudo-triangulations.
For example:

\begin{lemma}
\lemLab{CPTcount}%
Every combinatorial pseudo-triangu\-lation on $\nrVert$ vertices 
has $2\nrVert -3 +\nrNonPointed$
edges, where $\nrNonPointed$ is the number of non-pointed vertices.
\end{lemma}

\begin{proof}
Use the  counts of 
Theorem~\thmRef{count}, which are purely combinatorial.
%Alternatively, this statement is a special case of
%Lemma~\lemRef{cornercount} below (take $\nrCorners_1=b$).
\end{proof}

%\complaint{Moved definition of pointed CPT to Section 7.4. It makes more sense there (and is never used before). Paco, October 6, 2007}

%\medskip
 We say that a CPT \emph{can be stretched} if there is a
straight-line embedding of the graph $G$ topologically equivalent to
the given one and in which all angles are convex or reflex as
indicated by their labels. Not all combinatorial
pseudo-triangu\-lations can be stretched, as seen in
Figure~\figRef{cptn-2}.
% In both, the angles marked with a dot are the
% combinatorial corners of each bounded face.

\begin{itemize}
\item The one on the left of cannot be stretched since
  its graph is not generically rigid. (It has $2\nrVert - 3$ edges but it is
  not a Laman graph.)

\item The graph on the right is drawn as a true pseudo-triangulation, hence it is generically rigid.
But its assignment of big and small angles is not a stretchable CPT, since
the four boundary vertices and their two
neighbors form a set of $n'=6$ six pointed vertices whose induced subgraph has 
more than $3n'-3=9$ edges. This is in disagreement with what
Corollary~\corRef{generalized-Laman} predicts for a true pseudo-triangulation.
\end{itemize}
%
%\complaint{I redid this figure. Paco.}
\SingleFig{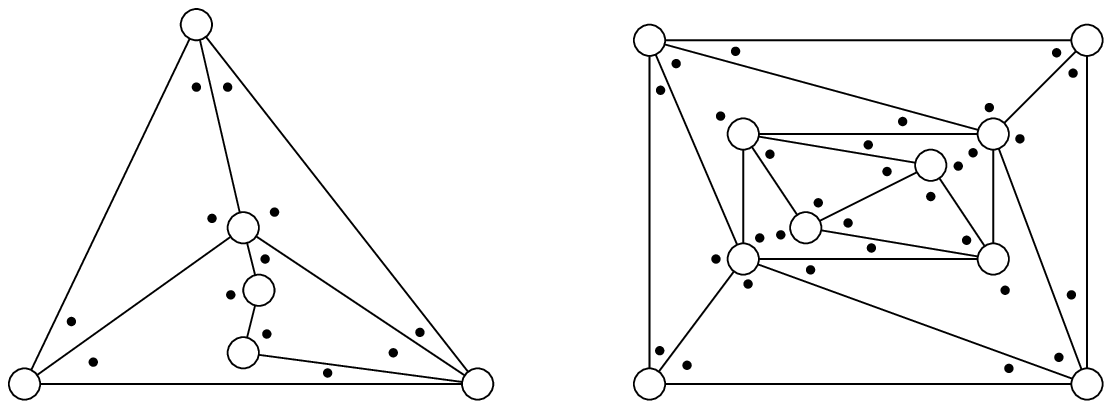}{Two combinatorial pseudo-triangulations that are not stretchable. Dots represent angles labeled ``small''.}

%\begin{figure}[ht]
%\vspace{-0.3in}
%\begin{center}
%\ \psfig{figure=Images/CPTN.eps,width=3.4in}%,width=3.0in}
%\end{center}
%\vspace{-0.3in} \caption{ (a) A plane graph with $2n-3$ edges and no
%CPT assignment possible.
%(b) A plane  graph with a non-stretchable CPT assignment.}%
%\figLab{CPTN}
%\end{figure}

The last example motivates the following definition:
\begin{definition}
\defLab{generalized-Laman}
A combinatorial pseudo-triangulation on a set $V$ of $\nrVert$ vertices is called
a \emph{generalized Laman CPT} if 
for every subset of $k\ge 2$ vertices in $G$, $l$ of them non-pointed, the subgraph induced by these vertices has at most $2k-3 + l$ edges.
\end{definition}

%Sometimes it will be more convenient to reformulate this property as follows:

%\begin{lemma}
%A combinatorial pseudo-triangulation on a set $V$ of $\nrVert$ vertices is
%{generalized Laman} if  and only if
%for every subset of $k\le 2$ vertices in $G$, $l$ of them non-pointed, the 
%number of edges incident to these vertices is at least $2k + l$.
%\end{lemma}

%This reformulation makes it clear that 
%the generalized Laman property in particular implies that every vertex has degree at least two
%and every non-pointed vertex has degree at least three, an obvious necessary condition for stretchability.

Corollary~\corRef{generalized-Laman} implies that being a generalized
Laman CPT is a necessary condition for a CPT to be stretchable. It is also sufficient, which breaks the proof of Theorem~\thmRef{planarrigid} into two parts:

\begin{theorem}
\thmLab{stretchable}%
Every generalized Laman CPT is stretchable to a true pseudo-triangulation.
\end{theorem}

\begin{theorem}
  \thmLab{cpt-gen-laman} For every planar topological embedding of a rigid
  graph, its angles can be labeled as a generalized Laman CPT.
\end{theorem}

We sketch the proofs of Theorems~\thmRef{stretchable} and~\thmRef{cpt-gen-laman} in Sections~\secRef{stretchable} and~\secRef{cpt-gen-laman} respectively. A different (and simpler) proof of~\thmRef{cpt-gen-laman} is given in Section~\secRef{cpt-matchings}  for the case of Laman graphs.

\subsection{Stretching Combinatorial Pseudo-Triangulations}
\secLab{stretchable}

To prove Theorem~\thmRef{stretchable}, we use a variation of Tutte's
\emph{barycentric method} to embed 3-connected planar graphs with
convex faces~\cite{t-hdg-63}.  The first step is to construct, from a
generalized Laman CPT, a certain auxiliary directed graph that is
``3-connected to the boundary'':

\subsubsection*{The 3-connected partially directed graph  of a CPT} 

A \emph{partially directed graph} $D=(V,E,\vec{E})$ is a graph $(V,E)$ together with
an assignment of directions to \emph{some} of its edges. Edges are
allowed to get two directions, one direction only, or remain
undirected. Formally, $\vec{E}$ is a subset of
the set which contains two opposite directed arcs for each edge of~$E$.
In what follows we say that $q$ is an out-neighbor of $p$ if there is an arc directed from $p$ to $q$.
%\complaint{Added last sentence, for referee 2. Paco, September 30, 2007}
% $E\cup (-E)$ where
%$E$ is the set of edges of $G$ with arbitrarily assigned directions.

\begin{lemma}
\lemLab{Gconstruction} For every combinatorial
pseudo-triangu\-lation $G$, we can construct a partially directed graph $D$
satisfying the following conditions:
\begin{enumerate}
\item $D$ is planar.
%\complaint{GR removed
%``and, as a graph, it contains the underlying graph of $G$.''
%E.g. An edge shared by two triangles will not appear at all.}
\item Every interior non-pointed vertex has as out-neighbors all its neighbors in~$G$.
\item For an interior pointed vertex $p$  of $G$, 
 let $\Delta$ be the pseudo-triangle of $G$ containing the big angle
at $p$. The out-neighbors of $p$ are the two neighbors of $p$ in the boundary of $\Delta$
together with one vertex of $\Delta$ not lying in the same pseudo-edge as $p$.
%\item For every pointed vertex $v$ of $\GP$ its out-neighbors in
%$D$ are exactly its neighbors in $\GP$.
\end{enumerate}
\end{lemma}

\begin{proof}
For each bounded face $F$ of $G$, add edges that triangulate $F$ with
ears only at (perhaps not all) its corners.
This ensures that each non-corner vertex
is incident to some added interior edge, which can be
oriented to satisfy condition~(3).
% Then orient the edges of this new graph as stated.
 See Figure \figRef{tutte-new} for an illustration.
%  First we extend $G$ % to be a (topological) triangulation by
%                      % triangulating
%  by adding edges into every pseudo-triangle $\Delta$ that is not a
%  triangle.  For every interior big angle we add a directed edge
%  through $\Delta$ towards some vertex not lying on the same
%  pseudo-edge. This can be done in many different ways; one only needs
%  to avoid crossings.
%%One way of doing this is recursively  dissecting each
%% face with an edge joining a non-corner to the
%% opposite corner (for the purposes of this recursion, the non-corner is then considered a corner of the
%% two faces obtained).
%  The edges on the boundary of $\Delta$ are oriented as required by
%  the statement.  See Figure \figRef{tutte} for an illustration of the
%  edges inserted and orientations given in a face~$\Delta$.  Observe
%  that a boundary edge incident to a corners of $\Delta$ may get a
%  second orientation if the corner is non-pointed or if the corner is
%  pointed and the edge is extremal.
\end{proof}
%\complaint{Made this proof more concise. Paco, October 6, 2007.
%changed afterwards GR}

%%%%%%
\SingleFig{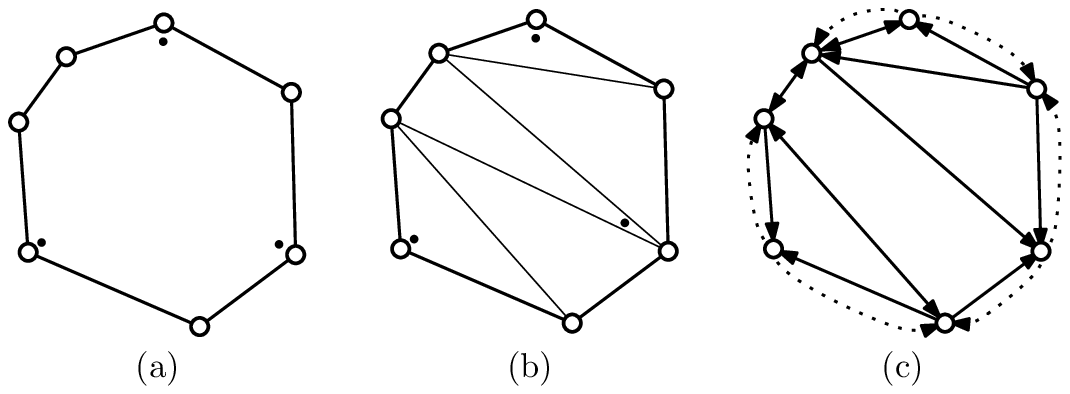}{(a) A CPT face, with 
  black dots indicating its three corners.  (b) % A compatible triangulation of the face. (c)
  A triangulation of it with ears only at corners. (c)
  One possible way for oriending the arcs of the auxiliary directed
  graph~$D$.  The dotted arrows indicate possible additional orientations
  for the boundary edges, depending on adjacent faces.}
  %\complaint{Changed caption, following Bereg's comment. Paco,
  %  September 30 2007. Rechanged + figure. Oct 10 GR}
%%%%%%

%\begin{figure}[ht]
%\vspace{-0.3in}
%\begin{center}
%\ \psfig{figure=Images/tutte.eps,height=1.6in}%,width=3.0in}
%\end{center}
%\vspace{-0.2in}
%\caption{Left: A combinatorial pseudo-triangular face, with a
%small black dot indicating a small angle (big angles are not
%marked). Middle: a compatible triangulation of the face. Right:
%the edges of the auxiliary directed graph.}
%\label{figure:Gconstruction}
%\end{figure}

%\medskip
We say that 
a plane embedding of a partially directed graph $(V,E,\vec{E})$ is
\emph{$3$-connected to the boundary} if from every interior vertex
there are at least three vertex-disjoint directed paths in
$\vec{E}$ ending in three different vertices of the unbounded face.
The following statement follows from
Theorem~6 in~\cite{orden:santos:servatiusB:servatiusH:combinatorialPT:2006} and
(part of) Theorem 7 in~\cite{streinu:haas:etAl:planarMinRigidPseudoTriang-confAndJour:2005}.

\begin{theorem}
%[\cite{streinu:haas:etAl:planarMinRigidPseudoTriang-confAndJour:2005}]
\thmLab{3connected}
If a CPT has the generalized Laman property, then
the auxiliary graph constructed in Lemma~\lemRef{Gconstruction} is $3$-connected to the boundary.
\end{theorem}

\subsubsection*{A directed version of Tutte's equilibrium method}

%The main result of this section can now be stated.

%\begin{theorem}
%\thmLab{stretchable-old}%
%For a combinatorial pseudo-triangu\-lation $G$ with non-degenerate
%\complaint{don't we assume gen. position?\\
%Yes, but this statement is purely combinatorial. We need to state a condition on the abstract graph. Paco.}
%(that is, simple polygonal) faces, the following are equivalent:
%\begin{enumerate}
%\item $G$ can be stretched into a compatible pseudo-triangu\-lation.%
%\item Every subgraph of $G$ with at least three vertices has at
%least three corners. %
%\item Every partially directed graph satisfying the requirements
%of Lemma \lemRef{Gconstruction} is $3$-connected to the boundary.
%\end{enumerate}
%\end{theorem}

%

%The implication from part 1 to part 2 is trivial. If  $G$ is
%embedded as a pseudo-triangu\-lation, there is no loss of generality
%in assuming that the embedding is in general position, so that every
%subgraph with at least three vertices has at least three convex hull
%vertices. And all convex hull vertices of a subgraph of $G$ will be
%corners according to our definition.
%We omit the proof from part 2 of 3, 
%(see the details
%in~\cite{streinu:haas:etAl:planarMinRigidPseudoTriang-confAndJour:2005})
%\complaint{proof of 2$\to$3 omitted, following complaint: ``This is the heart of the argument. Since it is omitted,
%why don't we omit the whole proof?''. Paco, October 2006}
%and now prepare to prove the implication from 3 to 1.

%

%\subsubsection*{\bf Tutte's equilibrium condition.} 

To finish the proof of Theorem~\thmRef{stretchable}
 we use a directed version of Tutte's Theorem on barycentric
embeddings of graphs. An embedding $D(\pts)$ of a partially directed
graph $D=(V,E,\vec{E})$ on a set of points $\pts=\{p_1, \ldots, p_n\}$,
together with an assignment $\omega \colon \vec{E}\rightarrow \reals$ of weights to
the directed edges is said to be in {\em equilibrium} at a vertex
$i\in V$ if
\begin{equation}
  \eqLab{tutte-equilibrium}
\sum_{j : (i,j)\in\vec{E}} \omega_{ij} (p_i-p_j)=0.
\end{equation}

The following is Theorem~8 in~\cite{streinu:haas:etAl:planarMinRigidPseudoTriang-confAndJour:2005}.
Its proof is not very different from the proof of Tutte's Theorem given in
\cite[Theorem~12.2.2, pp.~123--132]{richter-gebert:realizationSpacesPolytopes:1996}.
See the details
in~\cite{vpv-tbmai-03}
or~\cite{streinu:haas:etAl:planarMinRigidPseudoTriang-confAndJour:2005}.

\begin{theorem} [Directed Tutte Embedding]
\thmLab{directedTutte} Let $D=(\{1,\ldots,n\},E,\vec{E})$ be a
partially directed plane graph, $3$-connected to the boundary and
whose boundary cycle has no repeated vertices. Assume $(k+1,\ldots,n)$
is the ordered sequence of vertices in this boundary cycle and let
$p_{k+1},\ldots,p_n$ be the ordered vertices of a convex
$(n-k)$-gon.

Let $ \omega\colon \vec{E'}\rightarrow \reals$ be an assignment of
positive weights to the internal directed edges. Then:
\begin{itemize}
\item[(i)] There are unique positions $p_1,\ldots,p_k\in \reals^2$ for
the interior vertices such that all of them are in equilibrium.
% in the embedding $D(\pts), \ \pts=\{p_1,\ldots,p_n\}$.
\item[(ii)]
These positions yield a straight-line plane embedding of $D$.
All faces of $D$ are strictly convex polygons.
\qed
\end{itemize}
\end{theorem}

Observe that the system \eqRef{tutte-equilibrium} of equations that
define equilibrium is closely related to the system \eqRef{z} that
define the heights $z_i$ in a lifting of a \psT\
(Section~\secRef{pst-liftings}).  In fact, \eqRef{tutte-equilibrium}
decomposes into two independent systems, one for the $x$-coordinates
of the points $p_i$, and one for the $y$-coordinates. Uniqueness of
the heights in the Surface Theorem Surface Theorem
(Theorem~\thmRef{liftingbasis}) can be derived from uniqueness
of the solution of \eqRef{tutte-equilibrium}.

%\begin{proof}
%
%First, in the definition of {\em good
%representation} (Definition~12.2.6, page~126), each point $p_i$ is
%required to be in the relative interior of its {\em
%out-neighbors}, since this is what the directed equilibrium
%condition gives. Second, Claim~1 on page~126 proves that in a good
%representation it is not possible for a vertex $\pts$ that $\pts$ and
%all its neighbors lie in a certain line $\ell$, using
%3-connectedness. The proof can be adapted to use {\em
%$3$-connectedness to the boundary} as follows: consider three
%vertex disjoint paths from $\pts$ to the boundary. Call $q$ any of
%the three end-points, assumed not to lie in the line $\ell$.
%Complete the other two paths to end at $q$ using boundary edges in
%opposite directions. This produces three vertex-disjoint paths
%from $\pts$ to a vertex $q$ not lying on $\ell$. The rest needs no
%change.
%\end{proof}
%\bigskip

\begin{proof} [Proof of Theorem~\thmRef{stretchable}.]
Consider a partially directed graph $D$ constructed from our
CPT in the conditions of Lemma~\lemRef{Gconstruction}.
Choose arbitrary positive weights for the edges of $D$
and embed it in equilibrium, which can be done by
Theorems%s~\thmRef{3corners}, and 
~\thmRef{3connected} and Theorem~\thmRef{directedTutte}.
Taking into account that the equilibrium places
each interior vertex in the relative interior of the convex
hull of its out-neighbors, 
the conditions on $D$ stated in Lemma~\lemRef{Gconstruction}
imply that the
straight-line embedding of $G$ so obtained has convex and reflex
angles consistent with the original CPT labeling.
%: indeed, for a non-pointed vertex, all its
%out-neighbors are taken among the 
%edges of $G$, hence the equilibrium implies that all incident angles are small.
%For a pointed vertex $p_i$, its three out-neighbors are its extreme neighbors in the CPT
%plus a vertex of the same pseudo-triangle that is not a neighbor of $p_i$ in $G$, which implies that 
%the angle of $p_i$ in that pseudo-triangle is bigger than 180 degrees.
 \end{proof}

% \medskip
% \noindent {\bf Time Analysis.} Suppose that we are given a CPT that
% can be stretched. Tutte's theorem actually gives an algorithm to
% find a stretching: construct the auxiliary graph $D$ of Lemma
% \lemRef{Gconstruction}, choose coordinates for the boundary cycle
% in convex position and arbitrary positive weights for the directed
% edges, and then compute the equilibrium positions.

% Everything  can be done in linear time, except for the computation
% of the equilibrium position for the interior vertices. In this
% computation one writes a linear equation for each interior vertex,
% which says that the position of the vertex is the average of its
% (out-)neighbors. The position of the boundary vertices is fixed. It
% has been observed~\cite[Section
% 3.4]{chrobak:goodrich:tamassia:convexGraphDrawing:1996} that the
% planar {\it structure} of this system of equations allows it to be
% solved in $O(n^{3/2})$ time, using the $\sqrt n$-separator theorem
% for planar graphs in connection with the method of Generalized
% Nested Dissection
% (see~\cite{tarjan:lipton:rose:nestedDissection:1979,tarjan:lipton:planarSeparator:1980}
% or \cite[Section 2.1.3.4]{rote:thesis:1988}), or even in time
% $O(M(\sqrt n))$, where $M(n)=O(n^{2.375})$ is the time to multiply
% two $n\times n$ matrices.

\subsection{Generalized Laman CPT Labelings of Rigid Graphs}
\secLab{cpt-gen-laman} 

Here we sketch how to construct  a generalized Laman CPT labeling of the angles of any topologically embedded generically rigid graph $G$. 
In the next section we show an alternative method that works if $G$ is
minimally generically rigid, i.\,e., a Laman graph.
The general idea in the proof is to mimic combinatorially the Henneberg incremental construction of  pseudo-triangulations described in Theorem~\thmRef{henneberg}, using the fact that 
every generically rigid graph contains a Laman spanning subgraph. More precisely:

\subsubsection*{``Henneberg 1'' case:}
If $G$ has some vertex $i_0$ such that $G\setminus i_0$ is still
generically rigid, let $F$ be the face of the embedding of $G\setminus
i_0$ that contains $i_0$. By inductive hypothesis, assume that
$G\setminus i_0$ has been given a generalized Laman CPT labeling. With
a simple case study depending solely on the position of the edges
incident to $i_0$ with respect to the corners of $F$ (if $F$ has
corners; one of the cases is when $F$ is the unbounded face),
\cite[Lemma~5]{orden:santos:servatiusB:servatiusH:combinatorialPT:2006}
shows how to extend the labeling of $G\setminus i_0$ to a generalized
Laman CPT labeling of $G$.
By ``extend'' we mean that angles of $G\setminus i_0$ that are not bisected by the insertion of $i_0$ keep their labels.

\subsubsection*{``Henneberg 2'' case:}
Assume now that $G$ is a generically rigid plane graph on $\nrVert$ vertices, but that the removal of any vertex breaks rigidity. Let $L$ be a  Laman subgraph of $G$, with $2\nrVert -3$
edges. $L$ cannot have vertices of degree one, because they prevent generic rigidity, and it cannot have vertices of degree two, because their removal would leave $L$, hence $G$, generically rigid. Since the average degree is less than four, there is a vertex $i_0$ of degree three.

It can be proved that in these conditions $L\setminus i_0$ has ``one (generic) degree of freedom''; that is, that the linear space of infinitesimal motions of any sufficiently generic embedding of it has dimension one. In particular, for any two vertices $j$ and $j'$ whose distance (infinitesimally) changes in this unique infinitesimal motion, the graph $L\setminus i_0$ together with the edge $jj'$ is infinitesimally, hence generically, rigid. Thus, $G\setminus i_0$ with the same edge added is generically rigid, too. We say that such an edge $jj'$ \emph{restores rigidity}.
%
%Lemma 6 and Corollary 3 of~\cite{orden:santos:servatiusB:servatiusH:combinatorialPT:2006} show that:
The proof of Theorem~\thmRef{cpt-gen-laman} is finished with the
following lemma.
\begin{lemma}
  \lemLab{cpt-gen-laman} An edge $e$ that restores rigidity can always
  be found joining two neighbors \textup(in $G$\textup) of $i_0$ and with the
  property that every generalized Laman CPT labeling of $G\setminus
  i_0 \cup e$ can be extended to one of $G$.
\end{lemma}

The proof of this statement is much more involved than the Henneberg 1 case, and occupies
five pages in~\cite{orden:santos:servatiusB:servatiusH:combinatorialPT:2006} (Lemma 6 and Corollary~3).
The first step is to show that for \emph{every} choice of $e$ that restores rigidity, a CPT labeling of $G\setminus  i_0 \cup e$ can be found. But the difficult part
is to show that for \emph{some} choice of $e$ the generalized Laman property can be kept in the process. This involves planarity, rigidity and pseudo-triangulation arguments.

%\complaint{Section 7.3 could well finish here}
%, and I have mixed feelings about what follows. I want to stress that the proof of Lemma~\lemRef{cpt-gen-laman} is way more complicated than, for example, the one in the Henneberg 1 case. But, as an author of~\cite{orden:santos:servatiusB:servatiusH:combinatorialPT:2006} I don't want to be too self-laudatory. If you feel like cutting or rephrasing this self-laudatio, please do. Paco, October 6, 2007}
%We want to point out, however, that the proof of Lemma~\lemRef{cpt-gen-laman} is not as direct as it might seem. It occupies five pages in~\cite{orden:santos:servatiusB:servatiusH:combinatorialPT:2006} and has the following steps:

%\begin{itemize}
%\item There is at least one edge $e$ that restores rigidity and of the form $e=i_j i_{j+1}$ or $e=i_j i_{j+2}$, where $i_1,\dots,i_k$ denote the neighbors of $i_0$ in $G$ numbered cyclically with respect to the plane embedding of $G$. 
%
%For any such edge, every CPT labeling of $G\setminus i_0 \cup e$ can be extended to a CPT labeling of 
%$G$.

%\item If $e=i_j i_{j+1}$, then the extension preserves the
%generalized Laman property.  

%\item If $e=i_j i_{j+2}$, then either the same is true or $G\setminus i_0$ 
%contains a generically rigid subgraph with certain properties.

%\item If no $e=i_j i_{j+1}$ restores rigidity then the subgraph of the previous sentence cannot 
%exist for every $e=i_j i_{j+2}$ that restores rigidity, by
%arguments that involve planarity, rigidity and pseudo-triangulation properties.
%\end{itemize}

\subsection{Laman Graphs. CPT Labelings via Matchings}
%\complaint{Changed title.}
\secLab{cpt-matchings} 

Lemma~\lemRef{CPTcount} implies that a CPT with exactly $2\nrVert-3$
edges has all its vertices pointed. We call it a {\it pointed
  combinatorial pseudo-triangu\-lation} (or pointed CPT).  But, for a
pointed CPT, the generalized Laman property
(Def.~\defRef{generalized-Laman}) is simply the usual Laman property
of the underlying graph, which comes for free with the fact that
the graph is minimally generically rigid.
Thus, to prove Theorem~\thmRef{cpt-gen-laman} for a minimally
generically rigid graph one need not worry about the generalized Laman
property.

In particular, the Henneberg-like proof of
Theorem~\thmRef{cpt-gen-laman} sketched in the previous section can
 be greatly simplified for the Laman case. Even more so, it can be carried out at the
geometric level, each insertion step producing a straight-line pointed
pseudo-triangulation embedding of the graph in question by adding the
new point into the existing embedding. This is detailed
in~\cite[Section~3]
{streinu:haas:etAl:planarMinRigidPseudoTriang-confAndJour:2005} and
produces a direct proof of Theorem~\thmRef{planarrigid} for Laman
graphs.
%  \complaint{Rewrote Section 7.4 up to this point and
%  restructured it a bit starting at this point. Paco, October 6, 2007}
Bereg~\cite{bereg:certifying-rigid-graphs:2005} has shown that the
sequence of Henneberg steps that build up the graph, which parallels
the Henneberg ordering of Theorem~\thmRef{henneberg} at the abstract
graph level, can be found in $O(n^2)$ time, and hence the embedding
can be constructed in $O(n^2)$ time. 
% \complaint{Changed this to
%  include Bereg's reference. Paco, Sept 29, 2007. made more specific,
%  GR Oct 10}

A second proof of
Theorem~\thmRef{cpt-gen-laman}~\cite{streinu:haas:etAl:planarMinRigidPseudoTriang-confAndJour:2005}
is based on relating CPT labelings with matchings in a certain
bipartite graph constructed from our CPT. We include this proof here
for its simplicity.
% and because it has the advantage of being fully algorithmic.

%In what follows, we use the following convenient rephrasing of the Laman property:

%\begin{lemma}
%\lemLab{rephraseLaman} Let $G$ be a graph with $n$ vertices and
%$2\nrVert-3$ edges. $G$ is Laman if and only if every subset of $l\le \nrVert-2$
%vertices is incident to at least $2l$ edges.
%\end{lemma}

%\begin{proof}
%For the ``if'' part, 
%let $V_1$ be a subset of the vertices of $G$. We need to prove the Laman condition, that is, that
%if $V_1$ has $k\ge 2$  elements then the subgraph induced by it has at most $2k-3$ edges.

%For this, let $V_2$ be the complementary set of vertices, so that $|V_2|=\nrVert - k \le \nrVert -2$ and we can apply to it the
%hypothesis in the statement: $V_2$ is incident to at least $2(\nrVert - k)$ edges. The edges in the graph induced by $V_1$ are precisely those not incident to $V_2$ and, since the total number of edges is $2\nrVert - 3$, this number is at most $2k-3$, as desired.

%For the ``only if'' part we use the same arguments in reverse.
%\end{proof}

Let $G$ be a plane connected graph and consider the following bipartite graph $H$: one part is the set $V$ of vertices of $G$,
and the other part has $d-3$ nodes for each bounded face of degree $d$, and as many nodes for the outer face as its degree. The edges join the node of a vertex $i\in V$ to all nodes corresponding to faces incident to $v$ in the embedding of $G$.  The two parts of $H$ have equal sizes 
%if and only if
since the graph has $2\nrVert - 3$ edges (a necessary condition for the
existence of a pointed CPT labeling, by Lemma~\lemRef{CPTcount}).

\begin{lemma}
\lemLab{match}%
Pointed CPT labelings of a plane connected graph $G$ are in bijection with perfect matchings of $H$.
\end{lemma}

\begin{proof}
The edges in $H$ correspond to the assignment of reflex angles in $G$.
%In a pointed CPT labeling of a plane graph $G$ each vertex gets a reflex angle and each face, of degree say $d$, gets $d-3$ reflex angles, except for the outer face that gets $d$ of them.
\end{proof}

%We
%define a bipartite graph $H$ as follows: one part is $V$, 
%the set of vertices of $G$
%and has $\nrVert$ elements. The set $W$ corresponds to the faces $F$ of
%$G$ taken with certain multiplicities. For an interior face $f\in
%F'$ of degree (number of edges on the face) $d_f$, we will put
%$d_f-3$ vertices in $W$. For the outer face $f_o$ we will put
%$h=d_{f_o}$ nodes in $W$. The total number of elements in $W$ is
%thus $\sum_{f\in F'}(d_f-3) + h = \sum_{f\in F} d_f - 3 |F'|=2 |E| -
%3(\nrVert-2) = 2(2 \nrVert-3) - 3(\nrVert-2) = \nrVert $.

\SingleFig{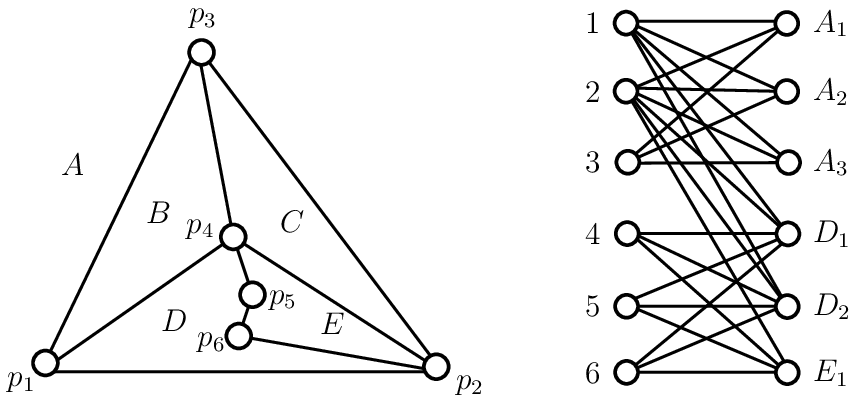}{
The $6$ vertices of the plane graph $G$ on the left, and its $5$ faces of degrees $3$ (outer
face $A$), $3$ (faces $B$ and $C$), $4$ (face $E$) and $5$ (face
$D$) lead to the bipartite graph $H$ on the right, with 
bipartition sets $V=\{1,2,3,4,5,6\}$ and $W=\{A_1,
A_2, A_3, D_1, D_2, E_1\}$.
 }

%\begin{figure}[ht]
%\vspace{-0.3in}
%\begin{center}
%\ \psfig{figure=Images/match.eps,width=3.2in}%,width=3.0in}
%\end{center}
%\vspace{-0.3in}
%\caption{A plane graph with $2n-3$ edges and its associated
%bipartite graph $H$.}%
%\figLab{match}
%\end{figure}
%
%A vertex $v\in V$ is connected in $H$ to the vertices in $W$
%corresponding to the interior faces $f$ of degree larger than $3$ to
%which it belongs in $G$, and to the vertices corresponding to the
%outer face (if it belongs to it). Hence, if $v$ belongs to the
%faces $f_1, f_2, \ldots $, and these faces have multiplicities $d_1,
%d_2, \ldots$ in $W$, then $v$ is connected to $d_1$ copies of the
%vertex for $f_1$, $d_2$ copies for $f_2$, etc.

We illustrate this result in  Figure~\figRef{match}. 
The horizontal edges in the bipartite graph $H$ form a perfect matching,
 and this induces a pointed CPT labeling of $G$. Observe, however, that $G$ is
 not a Laman graph; hence, $G$ cannot be stretched to a pointed pseudo-triangulation.

%Let $G=(V,E,F)$ be a plane graph with vertices $V$, edges $E$ and
%faces $F$. Assume $|V|=\nrVert$ and $|E|=2\nrVert-3$. Euler's relation implies
%that $|F|=\nrVert-1$. Denote by $F'$ the set of interior faces and by
%$f_o$ the outer face (with $h$ vertices), $F=F'\cup\{f_o\}$. 

%The connections (edges) in the bipartite graph $H$ represent
%potential assignments of {\it big} angles, where an {\it angle} is
%viewed as a pair {\it (vertex, face) }. Since each vertex must
%receive a big angle, we want a perfect matching. Since each interior
%face receives all but three big angles, and the outer face receives
%all big angles, the choice of multiplicities reflects just that.
%These considerations lead to the following Lemma.

%plane graphs satisfying the conditions of the previous
%Lemma may or may not have combinatorial pseudo-triangu\-lation
%assignments. (See again Figure \figRef{cptn}). But for Laman graphs,
%we are guaranteed a solution:

\begin{theorem}
[\cite{streinu:haas:etAl:planarMinRigidPseudoTriang-confAndJour:2005}]
\thmLab{match}%
If $G$ is a Laman graph, then $H$ has a perfect matching. Hence $G$
has a pointed CPT labeling. This labeling is automatically a
generalized Laman CPT labeling.
\end{theorem}

%Let us mention that the proof of this result
%in~\cite[Theorem 4]{streinu:haas:etAl:planarMinRigidPseudoTriang-confAndJour:2005}
%has a typographic error:
%the numbers of edges and vertices are swapped in formula
%(1).\complaint{I'd rather omit this. GR.}

\begin{proof}[Proof (Sketch)]
%\complaint{
%I did not trust/understand the previous one. 
%Incidentally, the one in the paper is wrong.
%In equation (1) $m_B$ and $n_B$ are the other way round...
%Paco, October 21.}
Let $W\subset V$ be a subset of $|W|=k$ vertices.
Let $F_W$ be the set of faces incident to the
vertices in $W$, and $R_W$ the union of those (closed) faces.
Let $D=\sum_{f\in F_W}d_f$. 
Hall's condition
for the existence of a perfect matching
 amounts to showing that $k \leq D- 3 |F_W|$.
(Actually, if one of the faces in $F_W$ is the unbounded one, 
Hall's condition would be  $k \leq D- 3 |F_W| + 3$, but 
the stronger inequality holds  anyway, assuming that $F_W$ does not
contain all faces.)
%The given vertices $W$ lie in the interior of $R_W$.
%
%%%%%%
%\SingleFig{faces}{A region $R_W$ which has two face-connected
%  components, shaded differently.  The one on the left has two holes,
%  the white triangle and the white non-convex pentagon.}
%%%%%%

If $R_W$ is not face-connected, we can prove the inequalities $k \leq
D- 3 |F_W|$ for each face-connected component separately and add
them up.  
If $R_W$ is face-connected, the inequality follows from the following three relations:
Euler's formula for the graph of all vertices and edges of $G$ in $R_W$,
Euler's formula to the graph of edges and vertices in the boundary of $R_W$, 
and the Laman property of $G$ applied to the first of these two subgraphs of $G$.
\end{proof}

%OLD PROOF BELOW: (Paco, October 22, 2006).

%\begin{proof}
%We will check Hall's condition to guarantee the existence of a
%perfect matching. Let $A\subset V$ be a subset of $|A|=a$ vertices.
%Let $F_A$ of size $|F_A|=f_a$ be the set of faces incident to the
%vertices in $A$, and let $D=\sum_{f\in F_A}d_f$. We need to show
%that $a\leq D-3 f_A$.

%It suffices to carry out the analysis on different face-connected
%components of $F_A$ separately. See Figure \figRef{faces}. One
%face-connected component is a polygon with (say) $b$ boundary edges,
%$b' \leq b$ boundary vertices, $h\geq 0$ holes and $e_i$ interior
%edges. We have $D=2 e_i +b$. By Euler's relation $(a+b') + (f_A + h
%+ 1) = (b + e_i) + 2$. Hence $f_A = e_i - a + 1 + b - b' - h$, where
%$\Delta := b - b' - h \geq 0$. Laman's condition implies that $b +
%e_i \leq 2(a+b') - 3$, hence $e_i \leq 2a + 2b'-b-3$. Now to show
%$a+3f_A \leq D$ we need $a + 3(e_i - a + 1 + b - b' -h) \leq 2 e_i +
%b$, i.w. $e_i \leq 2 a - 3 + 3 b' - 2b + 3h$. Since we know $e_i
%\leq 2a + 2b'-b-3$, it remains to show $2b' - b \leq 3b' - 2b + 3h$,
%i.e. $b \leq b'+3h$, which is obviously true.
%\end{proof}

%\medskip
%\noindent {\bf Time Analysis.} ????

\subsection{Rigidity versus Pointedness. A General Theorem}
\secLab{graph-drawing}

%We finish this chapter deriving from Theorem~\thmRef{planarrigid} a
%quite general statement about the possible number of pointed vertices
%in straight-line embeddings of a planar graph, in terms of the generic rigidity properties of it.

Let $G$ be a connected (abstract) graph with $\nrVert$ vertices and $\nrEdges$ edges.
The following generic rigidity parameters are associated to $G$:
\begin{itemize}
%\item Its \emph{generic rigidity rank}, $r$. By this we mean the rank of the rigidity map, in any generic 
%(not necessarily non-crossing) embedding of $G$. 

\item The number $d$ of generic \emph{degrees of freedom}, i.\,e., the
  dimension of the space of infinitesimal motions in any generic
  straight-line embedding of $G$.

\item The generic dimension $s$ of its space of \emph{self-stresses}.

\end{itemize}

Theorem~\thmRef{infinitesimal} relates the two as:
\begin{equation}
\eqLab{rigidity}
%r+s = 
\nrEdges = 2\nrVert -3 + s - d.
\end{equation}

Now let $G(\pts)$ be a straight-line embedding of $G$ on a point set $\pts$ in general position.
By Theorem~\thmRef{generalized-Laman} we have the equation
\begin{equation}
\eqLab{pointedness}
\nrEdges = 2\nrVert -3 + \nrNonPointed - \bar{\nrCorners},
\end{equation}
relating the following two pointedness-related parameters of $G(P)$:
\begin{itemize} 
\item The number $\nrNonPointed$ of non-pointed vertices.
\item The \emph{excess of corners} $\bar{\nrCorners}$: The
%\complaint{replace $\bar{\nrCorners}$ by $\bar\nrCorners$? }
total number of convex angles minus three times the number of bounded regions. 
\end{itemize}

{}From~\eqRef{rigidity}
 and~\eqRef{pointedness} we get the following equality between the
 generic rigidity parameters of $G$ and pointedness  parameters of $G(P)$:
\begin{equation}
\eqLab{pointed-vd-rigid}
 \nrNonPointed - \bar{\nrCorners}  = s - d.
\end{equation}
With this notation, Theorem~\thmRef{planarrigid} can be rephrased as: ``if $d=0$, then $G$ has embeddings with 
$\bar{\nrCorners} =0$". This generalizes to \emph{any} planar graph, as follows:

\begin{theorem}
\thmLab{graph-drawing}
\begin{enumerate}
 \item In every non-crossing straight-line embedding $G(\pts)$ of a planar graph $G$,  one has
 \[
 \quad
\nrNonPointed\ge s
\quad
\text{and hence}
\quad
\bar{\nrCorners}\ge d.
\]
 \item
Every planar graph $G$ has embeddings with: 
 \[
 \quad
\nrNonPointed= s
\quad
\text{and hence}
\quad
\bar{\nrCorners}= d.
\]
\end{enumerate}
\end{theorem}

\begin{proof}
%The equality in part 1 comes directly from~\eqRef{rigidity}
% and~\eqRef{pointedness}. 
The inequality $\nrNonPointed\ge s$ follows from the fact that every plane graph can be completed to a pseudo-triangulation with the same set of non-pointed vertices (Theorem~\thmRef{maximal}).
In the completion process, $s$ can only increase, and in the final pseudo-triangulation, it equals 
$\nrNonPointed$ by \eqRef{pointed-vd-rigid}, since in a pseudo-triangulation $d=\bar{\nrCorners}=0$. 
This proves part 1.

For part 2, we first prove by induction on $d$
 that edges can be
added to $G$ to make it generically rigid but keeping its planarity and its
self-stress dimension~$s$.
 (That is, we want every additional edge to decrease $d$ by one, instead of increasing~$s$.)
 If $d=0$ there is nothing to prove.
 If $d>0$, let $G(\pts')$ be any generic plane embedding. Since $G$ is not
 generically rigid, there is an infinitesimal flex in $G(\pts')$. Some face of the
 embedding must be deformed by this flex, or otherwise the flex would be
 trivial.  Thus, there is an edge between two vertices of this face that
 removes one degree of freedom from the space of motions. This edge may
 produce crossings in the embedding $G(\pts')$, but since it joins two
 vertices of a face, the new graph is still planar. 
 
Now, let $G_0$ be the resulting rigid planar graph with the same $s$ as $G$. By Theorem~\thmRef{planarrigid},
$G_0$ can be embedded as a pseudo-triangulation $G_0(\pts)$,
in which $\bar{\nrCorners}=d=0$ and hence $\nrNonPointed=s$.
If we remove the additional edges from this embedding to obtain an
embedding $G(\pts)$ of $G$, the edge removals can certainly not
increase $\nrNonPointed$.
But since $s$ remains constant, any decrease of $\nrNonPointed$ 
would violate part~(i). Hence  $\nrNonPointed=s$ holds for $G(P)$ as well.
\end{proof}

%The embeddings that attain the bounds of part 2 (which are, in general, not unique) can be called
 %\emph{maximally pointed} drawings of $G$.

%!TEX root = article.tex

%-----------------------------------------------------------------------------------------------     Polytope     ------

\section{Polytopes of Pseudo-Triangulations}
\secLab{polytope}
%
% \begin{theorem}
% \label{thm:PPT-polytope} \label{thm:PT-polytope} There is a polytope
% whose vertices are all pseudo-triangulations (pointed or not) of
% a given point set $\pts$ and whose faces correspond to the sets of
% pseudo-triangulations of the several ``pointgons" (including
% non-connected ones) ``contained in'' $\pts$.
% \end{theorem}

%Replace the start of the main file by:
%
%\documentclass[11pt]{article}
%\usepackage{latexsym}
%\usepackage{amsmath}
%\usepackage{amsfonts}
%\usepackage{amsthm}
%\usepackage{graphicx}

%%%%%%%%%%%%%%%%%%%%%%
\renewcommand{\reals}{\mathbb{R}}

\newcommand\expansionprime{\thetag{$\ref{eq:expansion}'$}}
%\newcommand\expansionprime{\eqRef{exprime}} % pdftex does not treat
                                % this correctly; (maybe I have an old
                                % version.)
                                % it has to be done by hand:

% To use with package "hyperref" (did not work in arXiV:
%\newcommand\expansionprime{\thetag{\hyperlink{exprime}{$\ref*{eq:expansion}'$}}}

In this section, we describe several high-dimensional polytopes and
polyhedra, whose skeletons represent the graphs of certain classes
of \psT s. We assume familiarity with basic notions of polyhedral
theory, such as vertices, faces, or extreme rays.
We refer to~\cite{grunbaum:polytopes,ziegler:lectures-on-polytopes:1995} 
for basic concepts in polytope theory.
%\complaint{Added references~\cite{grunbaum:polytopes,ziegler:lectures-on-polytopes:1995} for referee 2. Paco, September 30, 2007.}

 We start with the polyhedral cone of expansive motions of a point
set (the \emph{expansion cone}), as defined in the last section. Its
extreme rays correspond, in a certain way, to \ppsT s
(Theorem~\thmRef{expansion-cone}). There are variations of the
expansion cone in which the (pointed) \psT s are represented more
directly: the \emph{\ppsT\ polyhedron} and \emph{polytope}
(Theorem~\thmRef{expansion-polytope}), and the \emph{\psT\ polytope}
(Theorem~\thmRef{pt-polytope}). Finally, we will mention a polytope
corresponding to the pseudo-triangulations of a pointgon that lift to
locally convex surfaces (Section~\secRef{liftings}), the
\emph{regular \psT\ polytope} (Theorem~\thmRef{regular-polytope}).

\subsection{The Expansion Cone}
\secLab{expansion-cone}

In Section~\secRef{rigidity}, we considered mechanisms which had an
expansive infinitesimal motion, i.\,e., a motion in which all
pairwise distances are nondecreasing while certain other distances
are held fixed. Abstracting from this mechanism, we may study the
space of \emph{all} expansive infinitesimal motions
$v=(v_1,\ldots,v_\nrVert)$ for a point set $P=\{p_1,\ldots,p_\nrVert\}$. These
motions form a polyhedral cone in $(\reals^2)^\nrVert = \reals^{2\nrVert}$,
given by the $\binom \nrVert 2$ homogeneous linear inequalities in the $\nrVert $
vector variables $v_i\in\reals^2$:
\begin{equation}
  \eqLab{expansion}
  \langle p_i-p_j, v_i-v_j \rangle \ge 0,\ \forall i,j
\end{equation}
The rigid motions (translations and rotations) of the set $P$ as a
whole form a three-dimensional subspace of \emph{trivial motions}
 for which all inequalities \eqRef{expansion} are fulfilled as
equations.
%lineality space
To get rid of these trivial motions one can arbitrarily pin $p_1$ by
fixing $v_1$ at 0 (this eliminates the translations) and by
restricting the motion of $p_2$ to the line $p_1p_2$ (which then
eliminates the rotations). Thus we add the following normalizing
equations:
\begin{equation}
  \eqLab{pin}
  v_1= 0,\ \langle v_2, w\rangle = 0,
\end{equation}
where $w$ is a vector perpendicular to $p_2-p_1$. This results in a
($2n-3$)-dimensional polyhedral cone, which we call the
\emph{expansion cone} $\bar X_0=\bar X_0(P)$. The extreme rays of this cone
turn out to be exactly the expansive motions defined by the
pseudo-triangulation mechanisms of Theorem~\thmRef{infMechanism}.
\begin{theorem}
  \thmLab{expansion-cone}
For a point set $P$ in general position,
the {expansion cone} $\bar X_0$ is a pointed polyhedral cone.
Each extreme ray of $\bar X_0$ consists of the expansive
infinitesimal motions of a mechanism that is obtained by removing an
arbitrary convex hull edge from an arbitrary \ppsT\ of $P$, and each
mechanism of this type defines an extreme ray of $\bar X_0$.
\end{theorem}
This theorem will be proved as a consequence of
Theorem~\thmRef{expansion-polytope} below.
Theorem~\thmRef{infMechanism} about \ppsT\ mechanisms is an easy consequence of this theorem,
and thus we have another, very indirect, proof for
Theorem~\thmRef{infMechanism}, via polytopes.
%and hence for the existence of a

\begin{figure}[ht]
\begin{center}
\includegraphics{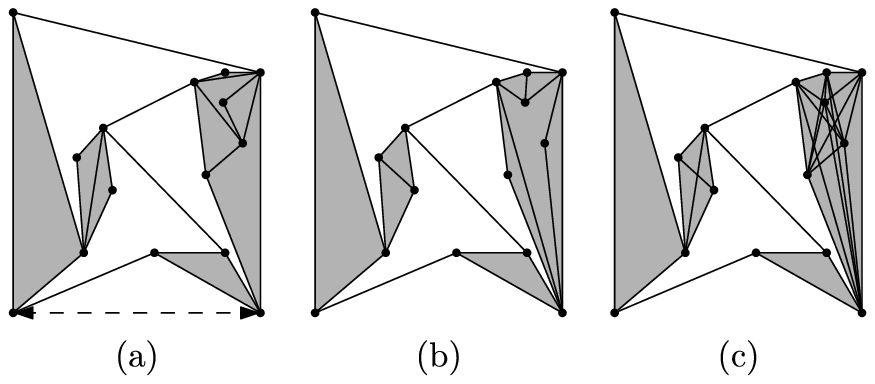}
\end{center}
\caption{ (a) A \ppsT\ mechanism representing an extreme ray $v$ of
the expansion cone. The rigid subcomponents are drawn shaded. (b)
Another \ppsT\ mechanism representing the same extreme ray. (c) The
set $E(v)$ of edges whose length is unchanged would be a canonical
representation of the ray~$v$.} \figLab{collapse}
\end{figure}

The correspondence between the extreme rays and the \ppsT s of $P$
is not one-to-one: consider two \psT s from which the same convex
hull edge has been removed; if they have the same rigid components,
they have the same expansive motions and thus they  define the same
extreme ray, see Figure~\figRef{collapse} for an example.
(See also Figure~\figRef{unfold} in Section~\secRef{carpenter} below.)
The rigid components in a \ppsT\ mechanism $T$ are formed by the
maximal convex regions enclosed by convex cycles in $T$, and they can
be identified in linear time:

\begin{theorem}[\cite{ss-crcpt-05}]
The rigid components of a
\psT\ from which a convex hull edge is removed are
% precisely
the parts of the graph that are enclosed by
the maximal convex polygons in the graph.
\end{theorem}

\subsection{The \PpsT\ Polyhedron and Polytope}
\secLab{ppt-polytope} One can obtain a polyhedron whose vertices are
in one-to-one correspondence with the \ppsT s of $P$ by modifying
the constraints \eqRef{expansion} as follows
\begin{equation*}
% \eqLab{exprime} % pdftex does not treat
                  % this correctly; (complains about duplicate pdf
                  % target labels.) Has to be done by hand:
  \langle p_i-p_j, v_i-v_j \rangle \ge f_{ij},\ \forall i,j,
 \tag  {$\ref{eq:expansion}'$}
 % \leqno
    %\smash{\raise\HyperRaiseLinkLength\hbox{\hypertarget{exprime}{}}}
    %\thetag{$\ref*{eq:expansion}'$}
\end{equation*}
where the quantities $f_{ij}$ are given by the following squared
$2\times 2$ determinants:
\[
%f_{ij} := \left| p_i\ p_j \right|^2
f_{ij} := \left | \begin{matrix}  p_i & p_j \\ \end{matrix} \right|^2.
\]

\begin{remark}
\ 
%\complaint{Added this remark, for referee 2. Paco, September 30, 2007.}
%
This choice of $f_{ij}$'s is not the only one that produces the polyhedron we want. 
Theorem~3.7 in~\cite{rote:santos:streinu:polytopePseudoTriangulations:2003} states the necessary and sufficient conditions that these quantities need to satisfy. In fact, the choice we adopt here is the case 
$a=b=0$ of the following more general choice, valid for any $a,b\in \reals^2$:
\[
f_{ij} = 
\left | \begin{matrix}  1 & 1 & 1 \\ a & p_i & p_j \\ \end{matrix} \right| \cdot
\left | \begin{matrix}  1 & 1 & 1 \\ b & p_i & p_j \\ \end{matrix} \right|.
\]
\end{remark}

The constraints \expansionprime\  and \eqRef{pin} define a
polyhedron $\bar X_f = \bar X_f(P)$, 
%\complaint{Shouldn't we omit the (P)
%everywhere? Say that $P$ is fixed throughout. Done GR. (not EVERYwhere}
 the \emph{\ppsT\ polyhedron}.
Moreover, we obtain a bounded polytope
 by setting some of the
equalities \expansionprime\ to equations
\begin{equation}\eqLab{exp-equal}
   \langle p_i-p_j, v_i-v_j \rangle = f_{ij}
\text{,\quad for all convex hull edges $ij$}
\end{equation}
The resulting polytope, defined by
 \expansionprime, \eqRef{pin} and \eqRef{exp-equal}, is the
 \emph{\ppsT\ polytope}
$X_f = X_f(P)$.

Note that, in the sequel, the term \emph{edges} will occur with two
meanings: edges of a \emph{polytope}, and edges $ij$ in a geometric
\emph{graph} on
the point set~$P$. The meaning
will always be clear from the context. When we speak about
\emph{vertices} in this section, we will always refer to polytope
vertices.

For each point $v\in \bar X_f$ we may define
% a corresponding graph by
the index set of tight inequalities:
\begin{displaymath}
  E(v) := \{\, ij \mid \text{\expansionprime\
holds as an equation for $v$}\,\}
\end{displaymath}
This set $E(v)$ is taken as the edge set of a geometric graph
on~$P$, the {support graph} of~$v$. By this correspondence, we get
precisely the \ppsT s:
 \begin{theorem}
   \thmLab{expansion-polytope} For a set $P$ of $n$ points with $\nrCH$
   points on the convex hull, $ X_f$ is a simple polytope of
   dimension $2n-3-\nrCH$, and $\bar X_f$ is a simple polyhedron of
   dimension $2n-3$.  $ X_f$ and $\bar X_f$ have the same set of
   vertices, and they are in one-to-one correspondence with the \ppsT
   s of $P$.  Two vertices of $ X_f$ or $\bar X_f$ are adjacent
   \textup(on the polyhedron\textup), if the corresponding \psT s are
   related by a diagonal flip.

The extreme rays of
 $\bar X_f$ are in one-to-one
 correspondence with the \ppsT s of $P$ with one convex-hull edge removed.
 \end{theorem}

In particular, the skeleton of $X_f$ is % nothing but
the
graph of \ppsT s defined in
Section~\secRef{graph-of}. %\complaint{added GR Dec 13}

A consequence of this theorem is that
\ppsT s are infinitesimally rigid (Theorem~\thmRef{pt-rigid}, for
we have given a different proof via the Maxwell-Cremona lifting
in Section~\secRef{rigidity}):
consider an arbitrary \ppsT\ $T$, and the corresponding vertex
$v$ with $E(v)=T$.
Since the polytope is simple,
the tight inequalities \expansionprime\ at $v$ must be
linearly independent. This means that the corresponding homogeneous
system of $2n-3$ equations
\begin{displaymath}
  % \eqLab{expansion}
  \langle p_i-p_j, v_i-v_j \rangle = 0,\ \forall ij\in E(v)
\end{displaymath}
together with \eqRef{pin} has only trivial solutions. In other words, the
support graph $E(v)$ is infinitesimally rigid.

The key statement of the proof is the following property of the set
of tight edges:
\begin{lemma}
\lemLab{Ev}
For a point $v\in\bar X_f$, % the support graph
  $E(v)$ cannot contain two crossing edges.
  $E(v)$ cannot contain three edges incident to
  a common point which make this point non-pointed.

In particular, $|E(v)|\le 2n-3$.
\end{lemma}
\SingleFig{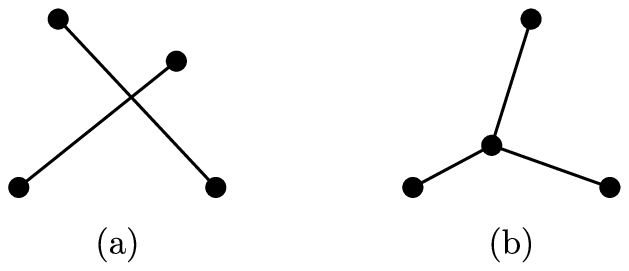}{(a) Four points in convex position (b)
in non-convex position.
The shown edges cannot be simultaneously tight.}
\begin{proof}[Proof (Sketch)]
 The statement of the lemma involves only four points: in
  the first case, they are four points in convex position, and in the
  second case; a triangle with a fourth point in the middle,
see Figure~\figRef{fourpoints}.
To prove the lemma, one has to show that, for
  all sets of four points, the inequalities for $ij\in E(v)$ cannot
  hold as equations while fulfilling the remaining inequalities
  \expansionprime.
This is done by taking an appropriate linear combination of these
constraints and deriving a contradiction, which
 boils down to showing that a certain linear combination of
the quantities $f_{ij}$ is positive. It turns out that this linear
equation is identically equal to~1.

The bound  $|E(v)|\le 2n-3$ follows from
Corollary~\corRef{max-pointed-noncrossing}.
% Theorem~\thmRef{charactPPT}, Conditions (4) and~(5).
%\complaint{Or better yet:
%Insert a corollary directly after the proof of Thm 2.6. Done GR}
\end{proof}

Remarkably, both statements of the lemma reduce to the same identity
involving the bounds $f_{ij}$. Rote et
al.~\cite{rote:santos:streinu:polytopePseudoTriangulations:2003}
give a whole family of alternative expressions for $f_{ij}$ that
satisfy the same identity and that can be used in \expansionprime.
This has only the effect of translating the  polyhedra
 $ X_f$ and
 $\bar X_f$ in $\reals^{2n}$, but it does not change their
combinatorial properties.

\begin{proof}[Proof of Theorem~\thmRef{expansion-polytope}]
The proof proceeds now in a somewhat indirect way.  First we look at
the polyhedron $\bar X_f$. It is easy to see that it has non-empty
interior in the $(2n-3)$-dimensional subspace defined by
 \eqRef{pin}, and it can be shown that it contains no line.
 Thus, it has dimension $2n-3$, and contains at least one
vertex $v_0$.  For every vertex $v$ in a $(2n-3)$-dimensional
polytope, $E(v)$ must contain at least $2n-3$ tight edges, but
Lemma~\lemRef{Ev} implies that $|E(v)|\le 2n-3$.  It follows that
there are exactly $2n-3$ tight inequalities, and $E(v)$ is a \ppsT;
hence $\bar X_f$ is a simple polyhedron.

The proof that \emph{all} \ppsT s appear as vertices of $\bar X_f$
is now somewhat indirect.
 Every vertex $v$ of the  polyhedron
 is incident to exactly
 $2n-3$ polyhedral edges, which lead to adjacent vertices or are
 infinite extreme rays. Each of these edges is characterized
by removing one element from $E(v)$. If this removed element
 $ij$ is
a boundary edge of the convex hull of~$P$,
 there is no polyhedron vertex $v'$
for which $E(v')$ contains $E(v)-\{ij\}$, and therefore the edge
leaving $v$ must be an extreme ray.
On the other hand, if the removed element $ij$ is an interior edge
of~$P$, there is only one other possible polyhedron vertex $v'$ for
which $E(v')$ contains $E(v)-\{ij\}$, namely the \psT\ $E(v')$
obtained by flipping $ij$.  (One can argue that this polyhedron edge
must be a bounded edge, i.\,e., it is not an extreme ray.)

We have thus proved that for every \ppsT\ $E(v)$ represented by a
vertex $v$ on the polyhedron, all its neighbors that are obtained by
flipping an edge are also represented on $\bar X_f$. Since we know
that  $\bar X_f$ has at least one vertex $v_0$,
 it follows
that \emph{all} \ppsT s are represented on $\bar X_f$. Thus, the
theorem is proved as far as $\bar X_f$ is concerned.
%
%\complaint {adapt to notation/terminology for the graph of \psT s!
%Sec.3). I think this is done. GR}

By intersecting $\bar X_f$ with the hyperplanes \eqRef{exp-equal},
one obtains a face $X_f$ of  $\bar X_f$ which contains all vertices
 of  $\bar X_f$ but none of its extreme rays.
Thus, $X_f$ is a bounded polytope that contains the same vertices
and (bounded) edges as $\bar X_f$.
\end{proof}

\begin{proof}[Proof of Theorem~\thmRef{expansion-cone}]
 $\bar X_0$ is obtained from $\bar X_f$ by replacing all right-hand
 sides $f_{ij}$ by 0, thus
 $\bar X_0$ is the recession cone of~$\bar X_f$.
It has a single vertex at the origin, and every extreme ray of
 $\bar X_0$ comes from one or several extreme rays of~$\bar X_f$.
It follows from the definition that every extreme ray of
 $\bar X_0$ is the expansive motion of a \psT\ mechanism.
\end{proof}

\subsection{The \PsT\ Polytope}
\secLab{pt-polytope} 

One can extend the \ppsT\ polytope to a
polytope representing \emph{all} \psT s, pointed or not, by
introducing a variable $t_i\ge 0$ for each point $p_i\in P$, and
modifying equations \expansionprime\ to become
%
% \complaint{There is one important feature of this polytope that is not mentioned: every non-crossing
% graph on $P$ appears as a face in the polytope, since they all can be extended to pseudo-triangulations by
% Theorem~\thmRef{maximal}. This is why the orden-Santos paper calls this ``the polytope of non-crossing graphs on a point set''. Similarly, the ppt-polytope has as faces all the pointed and non-crossing graphs on $P$.
% Paco, Nov 27, 2006. mentioned it. GR dec 11}
%
\begin{equation}
\eqLab{pt-polytope}
  \langle p_i-p_j, v_i-v_j \rangle
 + \|p_i-p_j\|\cdot (t_i + t_j)
 \ge f_{ij},\ \forall i,j,
\end{equation}
with the same values of $f_{ij}$,
and adding the equations
\begin{equation}
\eqLab{t}
t_i\ge 0, \forall i
\end{equation}
with equality for boundary vertices.
% on the boundary, we require equality
% \begin{equation}\eqLab{pt-polytope}
%    \langle p_i-p_j, v_i-v_j \rangle = f_{ij}
% \text{,\quad for all convex hull edges $ij$}
% \end{equation}
The polytope defined by
  \eqRef{t},   \eqRef{pt-polytope}, \eqRef{pin},  and \eqRef{exp-equal},
%and $t_i\ge 0$
 is the
 \emph{\psT\ polytope} $Y_f$.
 By definition, it contains 
the \ppsT\ polytope
 $ X_f$ as the face obtained by setting
 all the extra variables $t_i=0$.
 
%From  \eqRef{pt-polytope} and \eqRef{exp-equal} if follows that
%$t_i$ must be zero for all boundary vertices~$i$.

 \begin{theorem}[\cite{orden:santos:polytope:2005}]
  \thmLab{pt-polytope}
For any set $P$ of $\nrVert$ points with $\nrCH$ of them on the convex hull,
 $ Y_f(\pts)$ is a simple polytope of dimension $3\nrVert-3-2\nrCH$.
Its vertices
 are in one-to-one
 correspondence with the \psT s of $P$.

%In a vertex of $ Y_f$, $t_i=0$ holds if and only if the point $i$
%is pointed and equality in \eqRef{pt-polytope} if and only if $p_ip_j$
%is used.

Two vertices of $Y_f$
 are adjacent \textup(on the polytope\textup), if the corresponding
\psT s are related by a
\textup(diagonal, insertion, or deletion\textup) flip.

Moreover, the faces of the polytope are are in one-to-one correspondence with
the non-crossing graphs on~$P$.
\qed
\end{theorem}

Figure~\figRef{graphsq} in Section~\secRef{setAll}
(p.~\pageref{fig:graphsq}) shows the
4-dimensional polytope of all \psT s of a certain five-point set, in the
form of a Schlegel diagram. More precisely, the solid lines in the
figure form the polytope of pointed \psT s,
of dimension three, which is a wedge of two pentagons.
This is a facet of the 4-polytope of all pseudo-triangulations. The other facets
appear as a polyhedral subdivision of it into: two tetrahedra, two triangular prisms and
two more wedges of two pentagons. These six new facets correspond each to pseudo-triangulations
that use one of the six possible interior edges.
%complaint{Paco, is this geometrically accurate? is it even the
%  Schlegel
%diagram of SOME polytope? Or is it a schematic drawing?
%\\
%Made description more precise. And yes, it is a true Schlegel diagram. Paco.
% The \ppsT s are a face of this polytope; they form a 3-dimensional polytope,
% shown with solid edges.}

Each inequality  \eqRef{pt-polytope} or \eqRef{t}
defines a facet of $Y_f$. Setting an
inequality \eqRef{pt-polytope} to an equation corresponds to
insisting that an edge $ij$ is part of the \psT. Setting a variable
$t_i=0$ means that the corresponding vertex has to remain pointed.
Thus, each face of $Y_f$ corresponds to a set of
\emph{constrained} \psT s where certain edges are required
to belong to the \psT, and certain vertices are required to be pointed. 
This observation has
Theorem~\thmRef{constrained} as a corollary.

The \psT s of a pointgon $(R,P)$ can also be obtained, indirectly, as a face of $Y_f(P)$: for this, arbitrarily triangulate the
exterior of $R$, and consider the pseudo-triangulations of $\pts$
constrained to using the boundary of $R$ and this chosen
triangulation of the exterior. The vertices of the resulting face are in
one-to-one correspondence with the \psT s of $(R,P)$:
%\complaint{adapt notation!}
 
 \begin{theorem}
  \thmLab{pointgon-polytope}
For every pointgon $(R,P)$
% of $n$ points with $\nrCH$ points on the convex hull,
% $ Y_f$ is
there is
 a simple polytope
%of dimension $3n-3-2\nrCH$.
% Its
whose vertices
 are in one-to-one
 correspondence with the \psT s of $(R,P)$.
%Two vertices
% are adjacent (on the polytope), if the corresponding
%\psT s are related by a flip
%(a diagonal flip, an insertion flip, or a deletion flip).
\qed
\end{theorem}

\subsection{The Polytope of Regular \PsT s of a Pointgon}
\secLab{regular-polytope}

For a pointgon $(R,P)$ one can define another polytope whose
vertices represent certain \psT s of $(R,P)$, which are related to
 the locally convex liftings studied in
 Section~\secRef{liftings}.

A \emph{regular} \psT\ $T$ of a pointgon $(R,P)$ is a \psT\ of a pointgon
$(R,P')$ with $P'\subseteq P$ that can be
lifted to a locally convex function on $R$,
in such a way that every interior edge of $T$ is lifted to a strictly
convex edge (no two adjacent faces of $T$ are lifted coplanar).
%\complaint{adapt terminology!}

\begin{theorem}
[Aichholzer et al.~\cite{aichholzer:aurenhammer:brass:krasser:pseudoTriangulationsNovelFlip:2003}]
 \thmLab{regular-polytope}
For a pointgon $(R,P)$, there is a polytope whose vertices are in
one-to-one correspondence with the regular \psT s of $(R,P)$.

Edges on the polytope represent
 diagonal flips,  insertion flips, deletion flips,
 or vertex removal flips and their inverse.
\qed
\end{theorem}
The \emph{vertex removal} flips in this statement consist of the deletion of a vertex of degree two
and its incident edges. The need for this type of flip
comes from the remark we made after Theorem~\thmRef{convexfunction}.
%\complaint{GR I chose vertex REMOVAL in accordance with AAKB et al..,
% also in Section 4. Dec 13}

Note that the class of regular \psT s is quite different from the set
of all \psT s of a pointgon, which are represented as a face of
$Y_f(P)$, as discussed
in Theorem~\thmRef{pointgon-polytope} at the end of the previous subsection.
Firstly, we don't insist that all vertices of $P$ are used
in a regular \psT.  Secondly, a regular \psT\ will have no pointed
interior vertices, by Lemma~\lemRef{lifting-basic} in
Section~\secRef{liftings-subsec}.

When $R$ is convex, the regular \psT s coincide with the regular
triangulations of~$P$. In fact, the proof of the theorem closely
follows
the construction of the \emph{secondary polytope}, a polytope whose
vertices represent the regular triangulations of a point set~$P$
\cite{billera:filliman:sturmfels:constructionsSecondaryPolytopes:1990,
gelfand:kapranov:zelevinsky:discriminantsResultants:1994}.
%\complaint{Paco, is this the correct reference?\\
%Answer: yes, but I added a second one}
We sketch how this polytope is constructed.

For a given \psT\ $T$ of $(R,P')$, specifying a height
$h_i$ for every vertex $i\in P$ leads to a unique lifted surface,
by the Surface Theorem (Theorem~\thmRef{liftingbasis}).
This surface depends only on the heights of the complete vertices,
and it does so in a linear way.
The volume $V$ under the surface (or in other words, the integral of the
function whose graph is the surface, over the region~$R$) is therefore
a linear function of the heights $h_i$:
\begin{displaymath} % \nonumber
  V =  V_T(h_1,\ldots,h_n) = \sum _{i\in P} c_i \cdot h_i
\end{displaymath}
The \rp\ vertices and the vertices that are not used at all in $T$
have coefficients $c_i=0$ in this expression.
We use the vector $(c_1,\ldots,c_n)\in \reals^n$ to represent $T$.
The convex hull of these vectors forms the polytope of
Theorem~\thmRef{regular-polytope}.
Not all vectors will lie on the convex hull of the polytope.
It turns out that the vertices of the polytope
correspond to the regular \psT s.

\subsection{Delaunay \PsT s}
\secLab{delaunay}
The Delaunay triangulation of a point set is a very special
sample within the family of triangulations,
standing out with many remarkable properties.
%\complaint{do we need a reference?}
One would wish to have a similar object in the realm of \ppsT s.

Rote and Schulz~\cite{rs-pdpts-05} proposed a definition of a pointed
``Delaunay'' \psT\ of a \emph{polygon}~$R$.  The definition is based on locally
convex liftings of Section~\secRef{liftings}, but the argument why this is a
``reasonable'' definition uses the \ppsT\ polytope.
% and the analogy with the
% constrained Delaunay \emph{triangulation} of a polygon~\cite{c-cdt-89}.
The same concept has been considered, in a more 
general context, in~\cite{aichholzer:aurenhammer:hackl:preTriangulations:2006}.

For each corner $p_i=(x_i,y_i)$ of $R$, define $z_i^0 := x_i^2 + y_i^2$, and
for each reflex vertex, set $z_i^0$ to a value larger than all values at the
corners.  Then the highest locally convex function $f$ below these values will
induce a \psT~$T$ of $R$, by Theorem~\thmRef{envelope}. The choice of the
values $z_i^0$ for the reflex vertices ensures that $f$ does not achieve these
values, and hence $T$ will be pointed.
We call $T$ the pointed Delaunay \psT\ of~$R$.

To see why this definition might have some justification, let us first
consider a convex polygon $R$.  In this case, triangulations and \psT s
coincide, and one would certainly wish the ``Delaunay \psT'' to coincide with
the Delaunay triangulation.  This is indeed the case.  However, this
coincidence extends to some ``neighborhood'' of convex polygons: to see this,
we have to look at polytopes.  For a point set in convex position, the
vertices of the secondary polytope
\cite{gelfand:kapranov:zelevinsky:discriminantsResultants:1994,
  billera:filliman:sturmfels:constructionsSecondaryPolytopes:1990} represent
all triangulations, and a certain canonical objective function $c$ will select
the vertex corresponding to the Delaunay triangulation of~$R$.  The \ppsT\
polytope $X_f$ is an affine image of the secondary polytope, and hence there
is an analogous objective function $c'$, defined from the geometric parameters
of $R$, which selects the Delaunay triangulation on $X_f$.

\complaint{Better state a lemma here. "If $R$ satisfies ... then ....". Paco, October 14, 2007}
Now, we can simply use the same definition of $c'$ for the case when $R$ is
not convex.  (The secondary polytope and $X_f$ are no longer affinely
equivalent; they even have different combinatorial structures.)  This objective
function picks a vertex of $X_f$, which represents some \ppsT\ $T'$ of the
vertex set of~$R$.  It can be shown that for
% the class of \emph{neighborly visible}
a certain class of
polygons $R$, $T'$ contains the boundary of $R$, and in the interior of
$R$ it coincides with the pointed Delaunay \psT~$T$ defined above.
%  A polygon
% is neighborly visible if any two neighbors of a corner can see each other.
The precise condition for $R$ is as follows:
All corners of $R$ must lie on the convex hull, and any two corners
that are separated by just one corner along the boundary can see each other.
This class contains non-convex polygons, for which $T$ is a ``proper \psT''
and not just a triangulation, but it does not contain very ``convolved''
polygons.

Further properties of these Delaunay \psT s have not been
investigated so far. Also, there is no satisfactory definition of a ``pointed
Delaunay \psT\ of a \emph{point set}''.

%%%%%%%%%%%%% END OF SECTION 8 %%%%%%%%%%%%%%%

%!TEX root = article.tex

%-----------------------------------------------------------------------------------------------    Applications     ------

\section{Applications of Pseudo-Triangu\-lations}
\secLab{applications}

\subsection{Balanced Geodesic Triangulations for Ray Shooting}
\secLab{geodesic}
Balanced geodesic triangulations were introduced in%
%Goodrich and
%\complaint{earlier reference: Hersh et Snoe et ...!}
%Tamassia  have introduced 
~\cite{chazelle:etAl:rayShootingGeodesicTriangulations:1994,goodrich:tamassia:dynamicRayShooting:1997}
 as a data structure that performs ray shooting
queries and shortest path queries in a dynamically changing
connected planar embedded straight-line graph. 
%This work was done independently of
%other work on pseudo-triangu\-lations, and thus,
%pseudo-triangu\-lations don't explicitly appear in their paper.

Ray shooting refers to the problem of finding the first boundary
point that is hit by a query ray.  The classical approach to ray
shooting, say, in a simple polygon $P$ with $n$ vertices, uses a
triangulation of $P$. After locating the starting point of the query
ray in one of the triangles, one follows the ray from triangle to
adjacent triangle until the boundary is hit. The running time,
after the initial point location step, is proportional to the number
of triangles that are transversed, which can be at most $O(n)$,
depending on the triangulation. Similarly, a shortest path between
two query points can be found quite easily after identifying the
unique sequence of triangles which connects the two triangles
containing the query points.
\SingleFig{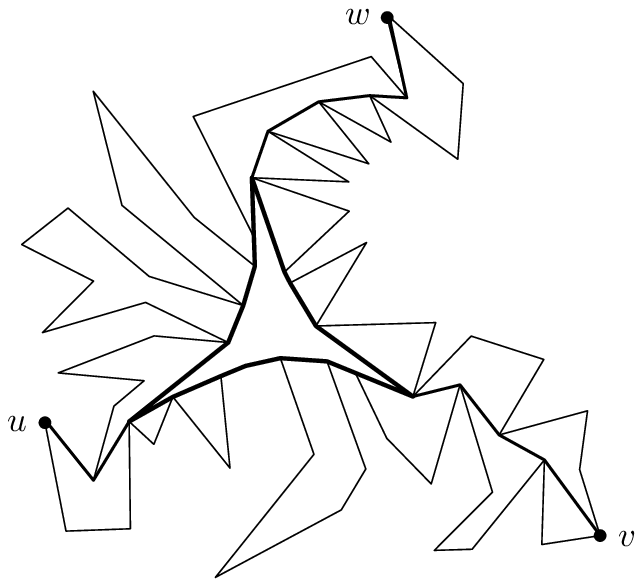}{A geodesic triangle.}

We would like to keep the number of triangles on such a path small and at the same
time, we want to maintain the search structure under changes of the
polygon.
For this purpose, we use geodesic triangulations instead of
triangulations.
%
%A geodesic path $uv$ in a simple $n$-gon $P$ is simply a
%shortest path between two given points $u$ and $v$ within this
% complaint{GR. this should be known by now, see intro}
% polygon.
 Three geodesic paths between $uv$, $vw$, and $uw$ will
form a \emph{geodesic triangle} which has a pseudo-triangle in the
center (which is called a ``deltoid region''
in~\cite{goodrich:tamassia:dynamicRayShooting:1997}) and possibly
some complicated paths at each corner which are shared between two
geodesic paths, see Figure~\figRef{geodesicTriangle}.  Let $\hat P$
be a convex $n$-gon whose vertices correspond to the vertices of $P$
as they appear on the boundary. Consider a triangulation $\hat T$ of
$\hat P$. For every triangle $\hat u\hat v\hat w$ in $\hat T$, we
can consider the corresponding geodesic triangle $uvw$ in $P$. The
set of these geodesic triangles will form a \emph{geodesic
triangulation} of $P$, see Figure~\figRef{geodesicTriangulation}.
\SingleFig{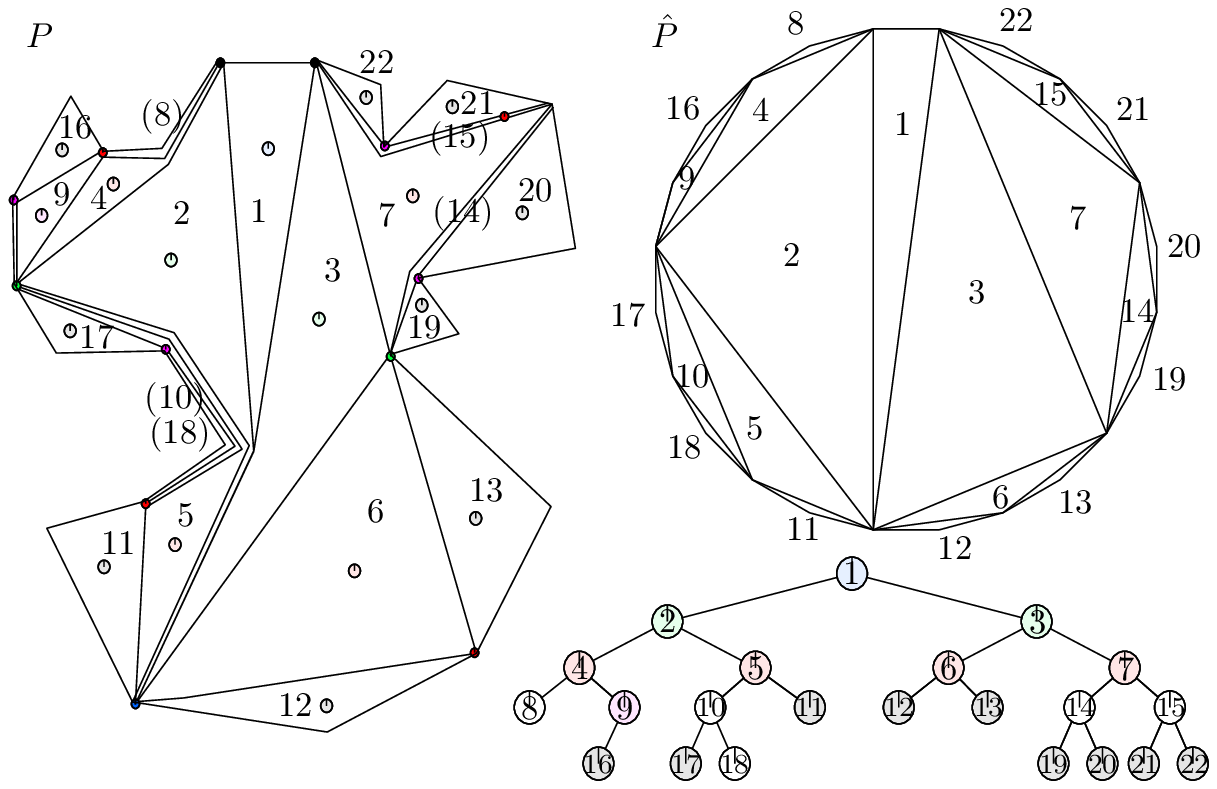}{A geodesic triangulation of a
polygon $P$, the
  corresponding triangulation $\hat T$ of the convex polygon $\hat P$, and the
  binary tree representation. Some geodesic triangles in $P$ have no area;
  their numbers are given in parentheses.}
 It
is a pseudo-triangu\-lation of the interior of $P$, but it stores
additional information about the correspondence between edges and
their original geodesic paths, and about adjacencies in the
triangulation $\hat T$.  The triangulation $\hat T$ of the convex
$n$-gon $\hat P$ can be represented as a binary tree (after
selecting an edge as the ``root'').  To follow the above-mentioned
paradigm for ray shooting queries, one has to walk through a
sequence of geodesic triangles.  Going from one geodesic triangle to
an adjacent geodesic triangle is no longer a constant-time
operation, because geodesic triangles have non-constant size, but it
can be carried out in logarithmic time, with appropriated data
structures. (A step from a geodesic triangle to an adjacent one may
be a zero-length step, in those regions where geodesic triangles
have no thickness.)

The advantage of using geodesic triangles it that is easier to get a
bound in the number of geodesic triangles traversed. If the
triangulation $\hat T$ is \emph{balanced} in the sense that it is
constructed by always splitting the remaining part of the boundary
into roughly equal parts, the number of triangles in any path
between two triangles is $O(\log n)$ (which is clearly best
possible).

% When the underlying polygon changes, it is may be necessary to
% change the structure of the tree by tree rotations in order to
% maintain balance. This results in flips of the underlying
% pseudo-triangu\-lation.  With appropriate data structures, all
% mentioned operations, including modifications of the underlying
% polygon (or more generally, a connected set of edges forming a
% planar subdivision) can be performed in $O(\log^2 n)$ time and
% linear space.

\complaint{The following cmplaint was here, but I guess is obsolete: "start with 9.2 VISIBILITY:
problem = visibility graph". Paco, October 14, 2007.
It results from an old remark that we should start with the *problem
statement* in this section. GR}

\subsection{Pseudo-triangu\-lations of Convex Obstacles}
\secLab{obstacles}

%\complaint{This is an EXTENSION of our stuff (but historically first)}

%In this section, we will briefly review a few problems from the area
%of motion control and visibility which have been solved with the aid
%of pseudo-triangu\-lations. 
In contrast to the previous sections,
which dealt with point sets, here we consider a collection
% $\mathcal O$
of $n$ disjoint convex obstacles in the plane, and define
pseudo-triangulations of them below.
%Fig.~\figRef{pseudoTriangPV}.
%, and we will modify the concept of
%pseudo-triangu\-lation for such a collection.
For simplicity, we assume that the bodies
are smooth. In this setting, a pseudo-triangle is not a polygon. It is
a region  bounded by a Jordan curve
consisting of three smooth inward-concave pieces,
each two meeting at a cusp, where they have a common tangent.
See a pseudo-triangulation of three convex bodies in 
Figure~\figRef{pseudoTriangPV}, where four pseudo-triangles arise. 
\complaint{Modified this first paragraph. Paco, October 14, 2007}

\SingleFig{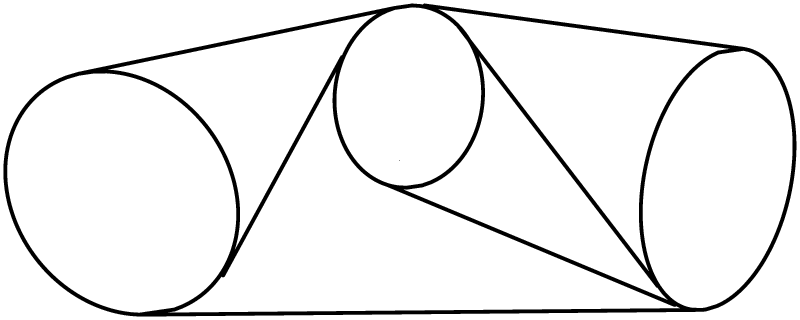}{A pseudo-triangu\-lation for
% of the convex hull of
three smooth convex obstacles.}

The analogue of the general position assumption is  that 
%the obstacles are smooth and strictly convex, and that 
no three obstacles have a common tangent line.
For a line which is tangent to two obstacles, we call the segment
between the two points of tangency a \emph{bitangent}. Each pair of bodies 
defines four bitangents. We call a bitangent \emph{free} it if does not
\complaint{I prefer to call common tangents bitangents even if they meet other bodies, and call
\emph{free} bitangents the ones that do not. Paco, october 14, 2005.}%
intersect any other obstacle. 

A \emph{pseudo-triangu\-lation} for a
set of convex obstacles is a maximal set of non-crossing
free bitangents, see~Figure~\figRef{pseudoTriangPV}. These
pseudo-triangu\-lations
have analogous properties to pointed pseudo-triangu\-lations of
point sets: there is an analogue of Theorem~\ref{thm:charactPPT},
every interior bitangent of a pseudo-triangulation can be flipped, etc.
In particular:

\begin{theorem}[Pocchiola and Vegter~\cite{pocchiola:vegter:visibilityComplex:1996}]
%  \lemLab{pt-objects-properties}
\ \null\   \begin{enumerate}
  \item
Any \psT\ of a set of $n$ convex obstacles
$O_1,O_2,\ldots,O_n$ decomposes the \emph{free space} inside the convex
hull, that is, the region
\[
\operatorname{conv}(O_1,O_2,\ldots,O_n) - (O_1\cup O_2
\cup\cdots\cup O_n),
\]
into
$2n-2$ pseudo-triangles and uses $3n-3$ bitangents.
%
%In this context, a pseudo-triangle is bounded by three smooth convex
%curves which meet at three \emph{cusps} where they have common
%tangents, such that the convex hull of the cusps contains the
%pseudo-triangle.  
\item The boundary of each pseudo-triangle is a sequence of
pieces which strictly alternate between bitangents and pieces of
obstacle boundaries.
%\item
%Every {pseudo-triangu\-lation} of a set of $n$ convex obstacles has
%$2n-2$ pseudo-triangles and $3n-3$ bitangents. 
\hfill\qed
  \end{enumerate}
\end{theorem}

Apart of the visibility applications mentioned in the next section, Pocchiola and Vegter~\cite{pocchiola:vegter:pc:99} later showed that the problem of finding a \emph{polygonal cover} of a collection of disjoint convex bodies is equivalent to finding a pseudo-triangulation of them.
%
%\complaint{Added this sentence for referee 1. Paco, September 30, 2007}

It is easy to extend these concepts to obstacles which
are not smooth, in particular to polygons. Conceptually, one has to
imagine that the polygons are rounded off at the vertices, and the
definitions have to be modified accordingly. (In particular, when 
two bitangents with a common endpoint make that vertex non-pointed, they
are regarded as crossing.)
Alternatively, if $O_1,O_2,\ldots,O_n$ are convex polygons with total vertex set $V$,
all pseudo-triangulations of $O_1,O_2,\ldots,O_n$ can be obtained as pointed pseudo-triangulations
of $V$ constrained to using all the boundary edges of all the $O_i$'s. This point of view allows to translate to triangulations of bodies many of the results about pseudo-triangulations of point sets (but it has to be noted that the work of Pocchiola and Vegter precedes the study of pseudo-triangulations of point sets).
\complaint{Added "Alternatively,...". Paco, October 14, 2007}

Less obvious is the extension to non convex polygonal objects described
by Angelier and Pocchiola in~\cite[p.~103]{angelier:pocchiola:sumofsquares:2003}:
they replace the vertices of the polygons by 
infinitesimal obstacles and the edges of the polygons by suitable bitangents joining them. 
%
%\complaint{Added this sentence for referee 1. Paco, September 30, 2007}

\subsection{The Visibility Complex}
\secLab{visibility-complex} 
The visibility complex 
%is a mathematical
%structure (a polyhedral complex) and also a data structure that
%naturally captures visibility in a set
% $\mathcal O$
%of {convex} obstacles.  It 
was introduced by Pocchiola and
Vegter~\cite{pocchiola:vegter:visibilityComplex:1996}.
We consider the set of directed \emph{visibility segments} or
\emph{maximal
  free} segments, which do not intersect any obstacle in the interior, but
which cannot be extended without cutting into an obstacle, see
Figure~\figRef{visibility}. Such a segment starts and ends at an
obstacle, or it extends to infinity in one or two directions (into
the ``blue sky'', which, for the purposes of this discussion, can be
treated like another obstacle ``at infinity''.)  The segments can be
moved continuously, forming a topological space, the \emph{visibility
complex}. This space is two-dimensional, as a segment can be locally
parameterized by, say, the slope and the signed distance from the
origin.  All segments which can be transformed into each other while
keeping their endpoints on the same two obstacles form a
two-dimensional \emph{face} of the {visibility complex}. A segment
reaches the boundary of a face when it becomes tangent to some
object.  An \emph{edge} of the {visibility complex} is thus formed
by a free segment which is tangent to an object and rotates
around this object while keeping its starting point and its terminal
point on two other objects.  Finally a \emph{vertex} of the
{visibility
  complex} corresponds to a \emph{free bitangent}, 
  in the sense of the previous section.

\SingleFig{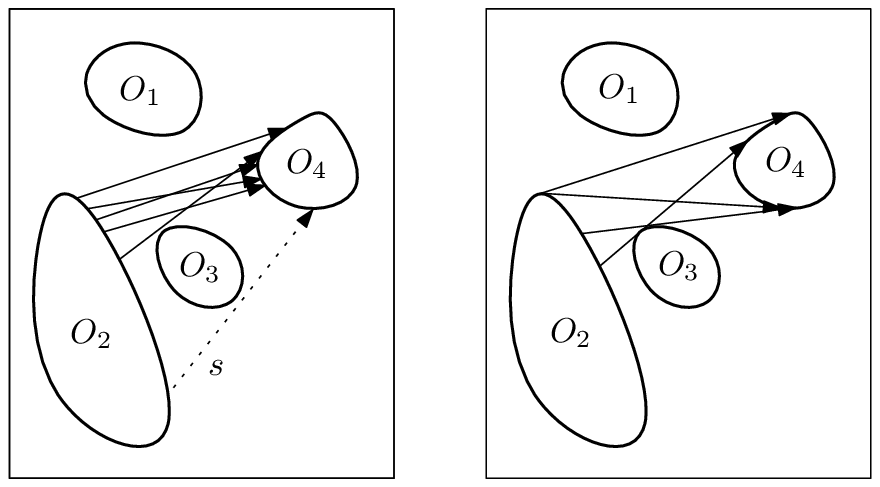} {A collection of some visibility segments
that belong to a common face of the visibility complex. Note that
the dotted segment $s$ does not belong to the same face as the other
segments although it starts and ends at the same two obstacles $O_2$
and $O_4$. The right part of the figure shows the segments corresponding to the
vertices of the face.}

The visibility complex, regarded as a set of vertices, edges, and
faces together with the incidences between them, forms an abstract
polyhedral complex, which can be stored as a data structure. (A
slightly different definition of the visibility complex, which
includes additional three-dimensional faces, was given in a
successor paper by Pocchiola and
Vegter~\cite{pocchiola:vegter:topoSweep:1996}.)

Under the general position assumption that no three obstacles share
a common tangent, the visibility complex has a quite regular
structure: every edge belongs to three faces, and every vertex
belongs to four edges and six faces. The face figure of a vertex
(formed by these edges and faces) has the combinatorial structure of
the graph of a tetrahedron.

The faces of the visibility complex correspond to pseudo-quadrangles
where two opposite sides are special: they are formed by a part of a
single obstacle boundary or by a convex hull edge.

The set of all free bitangents forms the \emph{visibility graph} of the
objects,
% $\mathcal O$,
which is a central concept in the context of visibility and shortest
path problems.  The number $k$ of free bitangents of the visibility graph
can vary in the range between $\Omega(n)$ 
(a tight lower bound of $4n-4$ is proved in~\cite{pocchiola:vegter:mtvg:96})
%
%\complaint{Added parenthetical comment for referee 1. Paco, September 30, 2007}
%
and $O(n^2)$.  The
visibility graph itself is not rich enough to allow the computation
of the set of visible points from a query point (the visibility
region).  This is where the additional structure of the visibility
complex is necessary. The total complexity of
the visibility complex
is $\Theta(k)$ if it has $k$ vertices.

Pocchiola and Vegter~\cite{pocchiola:vegter:visibilityComplex:1996}
have shown that the visibility region of a query point can be
determined from the visibility complex in $O(m\log n)$ time if its
size is $m$. The algorithm simply sweeps a ray of vision around the
query point, and it has to trace a corresponding path through the
visibility complex.

To construct the visibility complex for $n$ obstacles they have 
given~\cite{pocchiola:vegter:topoSweep:1996}  two different algorithms 
that work in $O(n\log n+k)$ time with only $O(n)$
intermediate storage, under the assumption that the common tangents
between two obstacles can be determined in $O(1)$ time. (Previous
algorithms for visibility graphs had achieved the same running time
but needed more than linear storage.) The algorithm of
\cite{pocchiola:vegter:topoSweep:1996}
% uses only simple data structures, but it 
establishes and exploits a partial order
structure on the set of free bitangents with strong properties. Roughly
speaking, two directed bitangents are related in this partial order
when the corresponding free segments can be continuously moved into
each other while always maintaining tangency at some obstacle, and
changing the direction monotonically.  Here the \emph{direction} is
measured as an angle in $R$, and two directed bitangents which are
otherwise the same but whose angle differs by a multiple of $2\pi$
are regarded as \emph{different}.

Pseudo-triangu\-lations arise as maximal antichains in this
order.  The algorithm flips through a sequence of
pseudo-triangu\-lations in a simple ``greedy'' manner and finds all
free bitangents on the way. A constant amortized time method for performing
each  flip was described in~\cite{angelier:pocchiola:sumofsquares:2003}
and has been implemented as part of the CGAL library~\cite{angelier:pocchiola:CGAL:2003}.
%\complaint{Added last sentence for refere 1. Paco, September 30, 2007}

\subsection{Pseudo-Triangu\-lations and Pseudoline Arrangements}
\secLab{pseudo-lines}
\complaint{Rewrote and expanded a bit this section. Paco, october 14, 2007}

The term \emph{pseudo-triangulation} was coined by Pocchiola
and Vegter \cite{pocchiola:vegter:orderType:94} because of an
interesting connection with pseudoline arrangements. 
%We conclude our
%list of applications with their intriguing observation.

\subsubsection*{Pseudoline arrangements} 
A \emph{pseudoline} is a
simple planar curve with end-points at infinity and
which partitions the Euclidean plane into two
parts. A \emph{pseudoline arrangement} is a collection of
pseudolines in which each pair of which has exactly one crossing.
A particular example is a \emph{line arrangement}. The combinatorial
type of a pseudo-line arrangement is determined by its facial
structure,
or equivalently, 
by the relative order in which each pseudo-line is crossed by all the other ones.

\subsubsection*{Tangents to pseudo-triangles. The dual pseudo-line}
For any given direction, a pseudo-triangle has a unique interior tangent
parallel to this direction, in the sense defined in
Section~\secRef{tangent} (before
Lemma~\lemRef{tangents}, p.~\pageref{section:tangent}).
%
%\complaint{should this be mentioned when tangents are defined? GR.
%Ileana: no, that would be 60 pages ago and nobody would remember.
%GR: Therefore added reference to 60 pages ago.}
%
Hence, if we dualize the lines carrying the tangents of a
pseudo-triangle (using the standard concept of point-line duality in
the Euclidean plane), the dual points form an
$x$-monotone curve, which is in particular a pseudo-line.

\subsubsection*{The dual pseudo-line arrangement}
Suppose now that we have a collection of pseudo-triangles with pairwise disjoint interiors
(for example, the ones in a pseudo-triangulation).
Any two of them have exactly one common interior tangent (see
Fig.~\figRef{bitangent}) so the dual pseudo-lines form a pseudo-line arrangement.
%%%%%%
\SingleFig{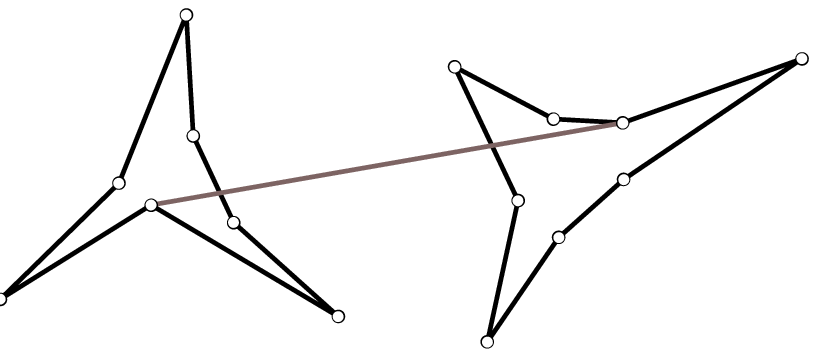}{Two pseudo-triangles and their unique common
tangent.}
%%%%%%

Pocchiola and Vegter~\cite{pocchiola:vegter:topoSweep:1996}
use this pseudo-line arrangement to reinterpret their visibility complex algorithm as a \emph{sweepline} algorithm. But they also
ask exactly what (combinatorial types of) pseudo-line arrangements can be obtained via collections of pseudo-triangles or, more specially, via pseudo-triangulations~\cite{pocchiola:vegter:orderType:94}. 
Any line arrangement can be obtained by considering the points dual to the lines of the arrangement and then changing each point for a tiny pseudo-triangle containing it.
But Pocchiola and Vegter show that also some non-stretchable arrangements can be obtained.

\subsection{Guarding Polygons with \texorpdfstring{$\pi$}{pi}-Guards}
\secLab{guarding}
\emph{Art Galleries} and \emph{Illumination} are a popular category
of geometric problems, where one asks for the number of
\emph{guards}, placed in the interior of a planar region so that
they would entirely cover it, or for light sources that would
illuminate it entirely. Bounds on the necessary number of guards
have been traditionally obtained using decompositions into convex
regions, in particular triangles, which can be covered with exactly
one guard placed at any vertex. Speckmann and T{\'{o}}th
\cite{speckmann:toth:vertexGuardsPseudoTriang:2003} improved the
known bounds for guarding a polygon with \emph{ restricted
visibility} guards by employing pseudo-triangulations. A $\pi$-guard
is a placement of a $\pi$ angle at a vertex of the polygon.

\begin{theorem}[Speckmann and T{\'{o}}th
\cite{speckmann:toth:vertexGuardsPseudoTriang:2003}]
Any simple polygon with $n$ vertices, $k$ of which are convex, can
be monitored with $\lfloor \frac{2n-k}{3} \rfloor$ edge-aligned
$\pi$-guards.
\qed
\end{theorem}

\subsection{Kinetic Data Structures for Collision Detection}
\secLab{KDS}
Pseudo-triangu\-lations were used to maintain a moving set of
objects in such a way that collisions can be detected quickly,
by Basch et al.~\cite{beghz-kcdts-04},
 Agarwal
et~al.~\cite{agarwal:basch:guibas:hershberger:zhang:kineticCollisionDetection:2003},
 and % in a different way
Kirkpatrick, Snoeyink, and Speckmann~%
\cite{speckmann:kirkpatrick:kineticMaintenance:2002,speckmann:kirkpatrick:snoeyink:kineticCollisionDetection:2002,s-kdscd-01}.
Consider a set of convex polygons which move simultaneously, under
external
control or autonomously, % of physical laws,
together with a pseudo-triangu\-lation of the free space between
them, in the sense of Section~\secRef{obstacles}.  As the points
move, the pseudo-triangles will change their shape, but it is easy
to check whether it remains a valid  pseudo-triangu\-lation:
\begin{prop}
  \begin{enumerate}
  \item Consider a pseudo-triangle whose vertices move.  It will be a valid
    pseudo-triangle as long as the following conditions are maintained, see
    Figure~\figRef{kineticPseudotriangle}\textup:
\begin{enumerate}
\item no two adjacent vertices coincide\textup;
\item the three corner angles remain positive\textup;
\item all other angles remain larger than $\pi$.
\end{enumerate}
\item Consider a pseudo-triangu\-lation of a set of convex polygons whose
  vertices move.  It will be a valid pseudo-triangu\-lation as long as
  \begin{enumerate}
  \item all pseudo-triangles remain valid\textup;
  \item all obstacles remain convex polygons\textup;
  \item and all exterior angles at the convex hull vertices remain larger than~$\pi$.
\hfill\qed
  \end{enumerate}
\end{enumerate}
\end{prop}

%%%%%%%%%%%%%
\SingleFig{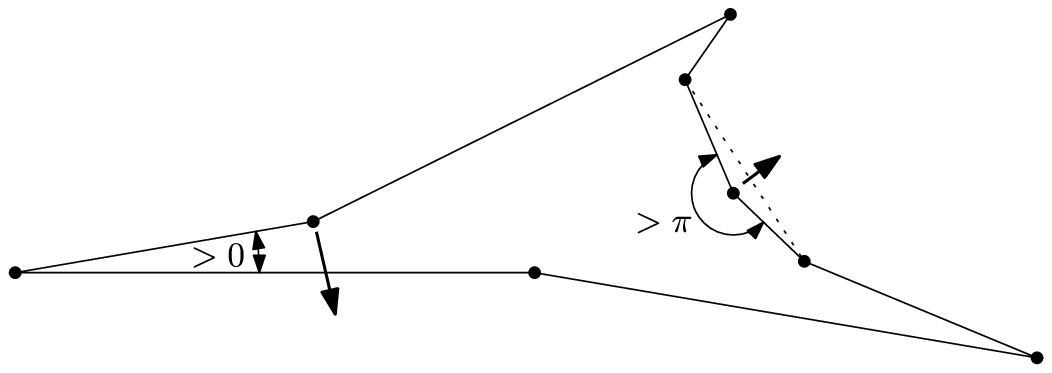}{Possible violations of the
pseudo-triangle condition when vertices move.}
%%%%%%%%%%%%%
Thus, when watching the motion of the obstacles, only the conditions
of the proposition have to be checked. For example, when the obstacles are
rigid or deformable convex polygons, the number of conditions is linear in the
number of obstacles and independent of the total number of vertices
\cite{speckmann:kirkpatrick:kineticMaintenance:2002,s-kdscd-01}.
When a condition becomes violated, a valid
pseudo-triangu\-lation can be restored by a flip.

The approach can be extend to non-convex obstacles, by inserting a
pseudo-triangu\-lation into the pockets of the obstacles
in a balanced way, as in Section~\secRef{geodesic}.
With
special care about how the maintenance and updates % flips
are done, 
this
pseudo-triangu\-lation data structure has good properties in terms
of the framework of kinetic data structures.
In particular, the overhead in running time is sensitive to the
``complexity'' of the scene. Objects whose convex hulls are disjoint can be
handled faster than interlocked pieces.

A special instance of such a setup, where the pseudo-triangu\-lation
is used also to guide the motion, is described
in the next section,
%
%in Section~\ref{sec:carpenter}, \complaint{section reference?}
%
 where it is used to ``unfold'' a polygon.

\subsection{The Carpenter's Rule Problem}
\secLab{carpenter}
The Carpenter's Rule Problem asks whether a simple planar polygonal
linkage can be continuously reconfigured to any other simple planar
configuration with the same edge-lengths, while remaining in the
plane and without creating self-intersections along the way. The
question was answered in the affirmative by Connelly, Demaine and
Rote~\cite{connelly:demaine:rote:carpenterRule:2003}. The
reconfiguring is done by first finding motions that convexify both
configurations with expansive motions (which guarantee non-colliding
motions), then taking one path in reverse.

We sketch now the subsequent algorithm of Streinu
\cite{streinu:pseudoTriangRigidityMotionPlanning-confAndJour:2005},
based on pseudo-triangu\-lation mechanisms. The algebraic details of
the implementation can be found in
\cite{streinu:combinatorialRoadmaps:2003}.

\subsubsection*{Overview of the Convexification  Algorithm} The
convexifying path, seen  as the collection of the $2n$ trajectories
of the $2n$ coordinates $(x_i,y_i), i=1, \cdots, n$ of the vertices
of the polygon, is a finite sequence of algebraic curve segments
(arcs) connecting continuously at their endpoints.

%%%%%%
\SingleFig{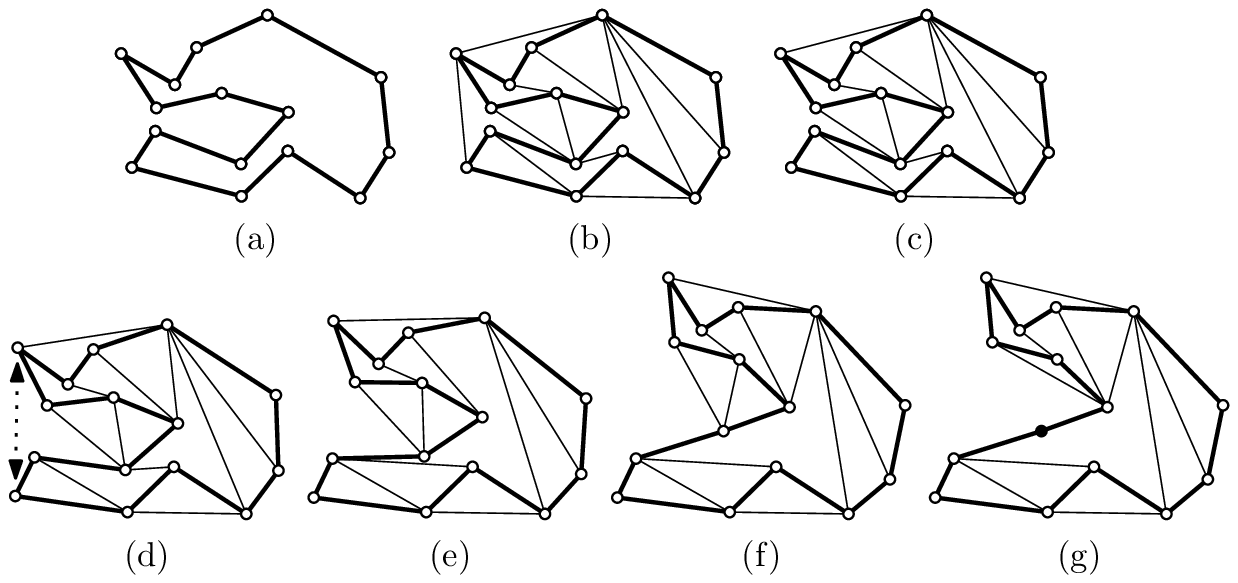}{(a) a simple polygon, (b) one of its pointed
pseudo-triangu\-lations, and (c)--(e) several snapshots in the expansive motion 
of  a pseudo-triangu\-lation mechanism
obtained by removing a convex hull edge. Between (d) and (e),
an alignment  event happens.  (f) Continuing the motion, the next event aligns two
polygon edges. (g) The aligned vertex (black) is frozen, the
pseudo-triangu\-lation is locally restructured and the motion can
continue.}
%%%%%%

Each arc corresponds to the unique free motion of the expansive,
one-degree-of-freedom mechanism induced by a planar pointed
pseudo-triangu\-lation of the given polygon, where a convex hull
edge has been removed and another edge has been pinned down. The
mechanism is constructed by adding $n-4$ bars to the original
polygon in such a way that there are no crossings, each vertex is
incident to an angle larger than $\pi$ and exactly one convex hull
edge is missing. See Fig.~\figRef{unfold}.
% \complaint
% {The different parts of Fig.~\figRef{unfold} are scaled differently.
% This is disturbing since motions are supposed to preserve lengths.
% (Moreover, x and y coordinates seem to be scaled differently,
% this destroys edge lengths as they rotate (but probably not noticeably}
 This
can be done algorithmically in $O(n)$ time. The mechanism is then
set in motion by pinning down one edge and rotating another edge
around one of its joints. The framework now moves {\it expansively},
thus guaranteeing a collision-free trajectory. One step of the
convexification algorithm consists in moving this mechanism until
two incident edges align (see Figure~\figRef{kineticPseudotriangle}).
At this moment it ceases to be a pointed
pseudo-triangu\-lation. We either freeze a joint (if the aligned
edges belong to the polygon) and locally patch a pointed
pseudo-triangu\-lation for a polygon with one less vertex, or
otherwise perform a local flip of the added diagonals. See Fig.
\figRef{unfold}.

\begin{algorithm}
(The Pseudo-Triangu\-lation Road-Map Algorithm)
\begin{enumerate}
\item {\bf Initialization:} Pseudo-triangulate the polygon. Remove
a convex hull edge to obtain a pseudo-triangu\-lation expansive
mechanism.
\item {\bf Repeat until the polygon becomes convex:}
\begin{itemize}
\item {\bf (Next Event)} Move the mechanism until an \emph{alignment event} occurs: two extreme edges at a vertex align.
\item {\bf (Freeze or Flip)} If the aligned edges were polygon
edges, \emph{freeze} them into a single edge by eliminating the
common vertex, and re-compute a compatible pseudo-triangu\-lation
mechanism. If one of the aligned edges is an added edge, drop it and
replace it by the edge extending over the two aligned edges (see
Fig. \figRef{unfold}(a-b)).
\end{itemize}
\end{enumerate}
\end{algorithm}

There are many ways to construct the initial pointed
pseudo-triangu\-lation or to readjust it at an alignment event. For
the sake of the analysis,
\cite{streinu:pseudoTriangRigidityMotionPlanning-confAndJour:2005}
uses a canonical pseudo-triangu\-lation based on shortest-path trees
inside the polygon and its pockets. This helps us maintain a global
integer valued cost function, the total number of bends in the
shortest paths, which is bounded by $O(n^2)$ for $n$ active (not
frozen) vertices. The cost function decreases by at least one at
flip-alignment events and increases at most $n-3$ times, at
freeze events. This analysis bounds the total number of events, and
thus the number of steps induced by simple pseudo-triangu\-lation
mechanism motions, by $O(n^3)$.

%%%%%%
\SingleFig{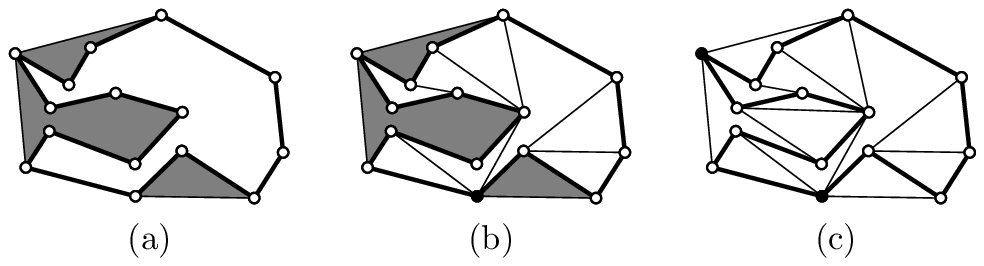}{(a) the pockets of the polygon from
  Fig.~\figRef{unfold}.  (b) a pseudo-triangu\-lation of the interior
  of the polygon formed by a shortest-path tree from the black vertex
  to all corners.  (c) a complete pointed pseudo-triangu\-lation
  obtained by adding shortest-path trees in all pockets.}
%%%%%%

\subsection{Spherical Pseudo-Triangu\-lations and Single-Vertex Origami}
An \emph{origami} is a piece of paper with creases, which is meant to
be folded into a three dimensional shape without bending or
stretching the paper (just folding along the creases). A very
special case, the \emph{single-vertex origami} illustrated in Figure
\figRef{origami}, turns out to be nothing but the Carpenter's Rule
Problem in spherical geometry, see Streinu and Whiteley
\cite{streinu:whiteley:origami:2005}.

%%%%%%
\SingleFig{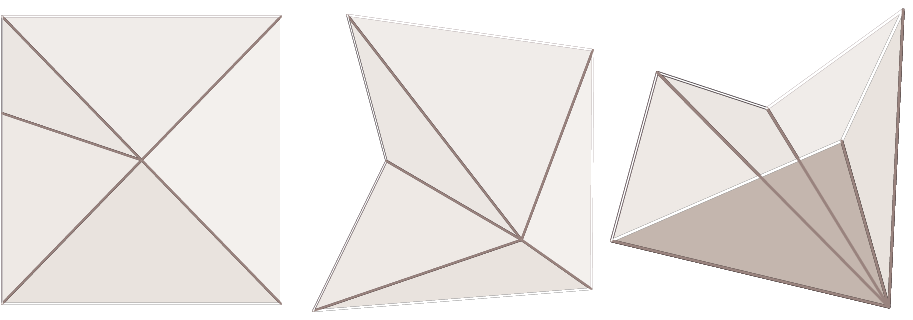}{A single-vertex origami fold: (a) the creased
piece of paper; (b,c) two of its possible folded states.}
%%%%%%

The idea is to associate to every planar framework, via central
projection, a framework on the sphere, with vertices placed as
points on a sphere and with edges along great-circles. This
connection is illustrated in Figure \figRef{sphereAndCone}.
Furthermore, by connecting the vertices with the center of the
sphere, we obtain a series of triangles (which behave in $3$-space
like {rigid panels}) connected along {hinges} into a {conical}
structure, as in Figure~\figRef{sphereAndCone}. This simple sequence
of transformations transforms planar pseudo-triangulations into
spherical or conical structures which inherit all the expansiveness
properties of the planar ones.

%%%%%%
\SingleFig{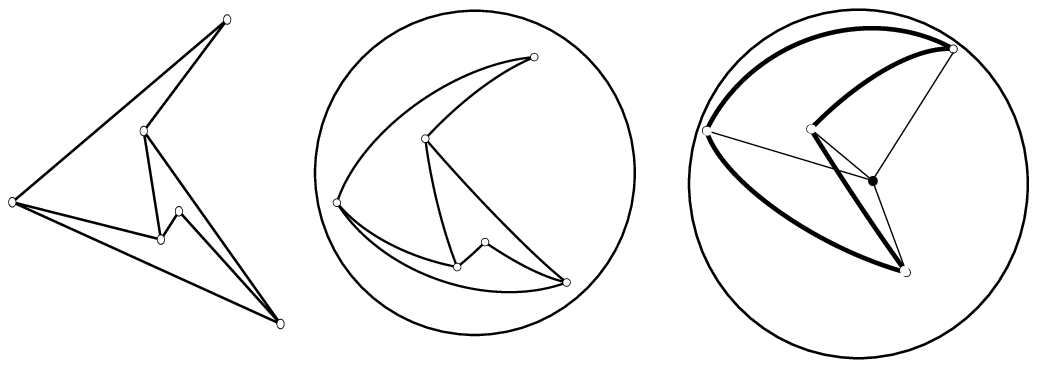}{A planar pseudo-triangulation, its
corresponding spherical version, and a conical panel-and-hinge
structure arising from a (different) pseudo-triangu\-lation.}
%%%%%%

By translating all relevant concepts of the pseudo-triangulation
road-map algorithm for the planar Carpenter's Rule problem, one
obtains:

\begin{theorem} [Streinu and Whiteley
\cite{streinu:whiteley:origami:2005}] Every simple spherical polygon
with perimeter at most $2\pi$ can be convexified in a hemisphere.
Every single vertex origami can be folded from a flat piece of
paper.
\end{theorem}

\subsection{Spherical Pseudo-Triangu\-lations and Convex Geometry}
 In the introduction, we mentioned that a \psT\ in the plane can have at most
one pseudo-$\nrCorners$-gon with $\nrCorners<3$, namely the outer face.  This
property remains true for \psT s on the sphere that are restricted to lie in a
hemisphere, such as the ones that arise in the previous subsection.  However,
if we look at \psT s on the whole sphere, it turns out that one can have an
arbitrary number of pseudo-$2$-gons. % (pseudo-digons).
Therefore, the bound of $2n-3$ on the number of edges of pointed graphs, which
follows from Theorem~\thmRef{generalized-Laman}, does not hold;
pointed graphs with $2n-2$ or more edges exist. As graphs in the
sphere these graphs  support
a self-stress  and have therefore a piecewise linear lifting
(in an appropriate sense).
  Panina~\cite{p-bcah-05} 
has used these liftings to
construct counter-examples to a conjecture of A.~D.~Alexandrov about a
characterization of the sphere among the smooth convex surfaces. The crucial
property, which follows from pointedness, is that the liftings are
\emph{saddle functions}, in a sense analogous to the properties of pointed
vertices $p$ in the liftings of pseudotriangulations with a unique non-pointed
vertex, which were mentioned at the end of Section~\secRef{reciprocal}: at the
lifted point $p'$, there is no supporting plane which intersects the
neighborhood of the lifted surface $f$ only in this point (and leaves the
surface locally on one side of it).
See Panina~\cite{panina:2006} for more information on the connections
between \psT s, saddle function, and so-called hyperbolic virtual polytopes.

\subsection*{Acknowledgments.}
We thank the referees for their extensive comments.

%-----------------------------------------------------------------------------------------------    Bibliography     ------

%\advance\itemsep by 0pt plus 0.5pt
%\bibliographystyle{amsalpha}
\bibliographystyle{abbrv}
%\bibliographystyle{plain}
%\bibliography{article}
\bibliography{pt-survey}

\begin{thebibliography}{10}
\advance\itemsep by 0pt plus 0.5pt

\bibitem{agarwal:basch:guibas:hershberger:zhang:kineticCollisionDetection:2003}
P.~Agarwal, J.~Basch, L.~Guibas, J.~Hershberger, and L.~Zhang.
\newblock Deformable free space tilings for kinetic collision detection.
\newblock {\em International Journal of Robotics Research}, 21:179--197, 2002.

\bibitem{aichholzer:aurenhammer:hackl:preTriangulations:2006}
O.~Aichholzer, F.~Aurenhammer, and T.~Hackl.
\newblock Pre-triangulations and liftable complexes.
\newblock In N.~Amenta and O.~Cheong, editors, {\em Proceedings 22nd Ann. Symp.
  Comput. Geom., Sedona, Arizona, USA, June 5--7, 2006}, pages 282--291. ACM,
  2006.

\bibitem{aahk-tstpt-06}
O.~Aichholzer, F.~Aurenhammer, C.~Huemer, and H.~Krasser.
\newblock Transforming spanning trees and pseudo-triangulations.
\newblock {\em Inf. Process. Lett.}, 97(1):19--22, 2006.

\bibitem{aichholzer:aurenhammer:krasser:adaptingPseudoTriangulations:2003}
O.~Aichholzer, F.~Aurenhammer, and H.~Krasser.
\newblock Adapting (pseudo)-triangulations with a near-linear number of edge
  flips.
\newblock In {\em Proc. 8th International Workshop on Algorithms and Data
  Structures (WADS)}, volume 2748 of {\em Lecture Notes in Computer Science},
  pages 12--24, 2003.

\bibitem{aichholzer:aurenhammer:brass:krasser:pseudoTriangulationsNovelFlip:20%
03}
O.~Aichholzer, F.~Aurenhammer, H.~Krasser, and P.~Brass.
\newblock Pseudo-triangulations from surfaces and a novel type of edge flip.
\newblock {\em SIAM Journal on Computing}, 32:1621--1653, 2003.

\bibitem{aichholzer:aurenhammer:krasser:speckmann:convexityMinimizes:2004}
O.~Aichholzer, F.~Aurenhammer, H.~Krasser, and B.~Speckmann.
\newblock Convexity minimizes pseudo-triangulations.
\newblock {\em Computational Geometry: Theory and Applications}, 28:3--10,
  2004.

\bibitem{aichholzer-etal:number-planegraphs:2007}
O.~Aichholzer, T.~Hackl, C.~Huemer, F.~Hurtado, H.~Krasser, and B.~Vogtenhuber.
\newblock On the number of plane geometric graphs.
\newblock {\em Graphs and Combinatorics}, 23[Suppl.]:67--84, 2007.

\bibitem{aichholzer:hoffmann:speckmann:toth:degreeBoundsConstrainedPseudoTrian%
gulations:2003}
O.~Aichholzer, M.~Hoffmann, B.~Speckmann, and C.~D. T\'oth.
\newblock Degree bounds for constrained pseudo-triangulations.
\newblock In {\em Proc. $15$th Canadian Conference on Computational Geometry
  (CCCG 2003)}, pages 155--158, Halifax, Nova Scotia, Canada, 2003.

\bibitem{aichholzer:hurtado:noy:lowerBoundTriangulations:2004}
O.~Aichholzer, F.~Hurtado, and M.~Noy.
\newblock A lower bound on the number of triangulations of planar point sets.
\newblock {\em Computational Geometry: Theory and Applications}, 29:135--145,
  2004.

\bibitem{aichholzer:orden:santos:speckmann:numberPseudoTriangulationsCertainPo%
intSets:2006+}
O.~Aichholzer, D.~Orden, F.~Santos, and B.~Speckmann.
\newblock On the number of pseudo-triangulations of certain point sets.
\newblock {\em J. Combin. Theory Ser. A}, 2007.
\newblock To appear.

\bibitem{streinu:aichholzer:rote:speckmann:zigZagPath-wads:2003}
O.~Aichholzer, G.~Rote, B.~Speckmann, and I.~Streinu.
\newblock The zig-zag path of a pseudo-triangulation.
\newblock In {\em Proc. 8th International Workshop on Algorithms and Data
  Structures (WADS)}, volume 2748 of {\em Lecture Notes in Computer Science},
  pages 377--388, Ottawa, Canada, 2003. Springer-Verlag.

\bibitem{angelier:pocchiola:sumofsquares:2003}
P.~Angelier and M.~Pocchiola.
\newblock A sum of squares theorem for visibility complexes and applications.
\newblock In B.~Aronov, S.~Basu, J.~Pach, and M.~Sharir, editors, {\em Discrete
  and Computational Geometry -- The Goodman-Pollack Festschrift}, pages
  77--137. Springer-Verlag, Berlin, 2003.

\bibitem{angelier:pocchiola:CGAL:2003}
P.~Angelier and M.~Pocchiola.
\newblock \textsc{Cgal}-based implementation of visibility complexes.
\newblock Technical Report ECG-TR-241207-01, Effective Computational Geometry
  for Curves and Surfaces (ECG), 2003.

\bibitem{asimow:roth:rigidityOfGraphs:1978}
L.~Asimow and B.~Roth.
\newblock The rigidity of graphs.
\newblock {\em Transactions of the American Mathematical Society},
  245:279--289, November 1978.

\bibitem{avis:fukuda:reverse-search:1996}
D.~Avis and K.~Fukuda.
\newblock Reverse search for enumeration.
\newblock {\em Discrete Applied Mathematics}, 65(1-3):21--46, 1996.

\bibitem{beghz-kcdts-04}
J.~Basch, J.~Erickson, L.~Guibas, J.~Hershberger, and L.~Zhang.
\newblock Kinetic collision detection for two simple polygons.
\newblock {\em Computational Geometry: Theory and Applications}, 27:211--235,
  2004.

\bibitem{bereg:transforming-pseudo-triangulations:2004}
S.~Bereg.
\newblock Transforming pseudo-triangulations.
\newblock {\em Information Processing Letters}, 90(3):141--145, 2004.

\bibitem{bereg:certifying-rigid-graphs:2005}
S.~Bereg.
\newblock Certifying and constructing minimally rigid graphs in the plane.
\newblock In {\em Proc. 21st Annu. Sympos. Comput. Geom.}, pages 73--80, 2005.

\bibitem{bereg:enumerating-pseudo-triangulations:2005}
S.~Bereg.
\newblock Enumerating pseudo-triangulations in the plane.
\newblock {\em Comput. Geom. Theory Appl.}, 30(3):207--222, 2005.

\bibitem{billera:filliman:sturmfels:constructionsSecondaryPolytopes:1990}
L.~Billera, P.~Filliman, and B.~Sturmfels.
\newblock Constructions and complexity of secondary polytopes.
\newblock {\em Advances in Mathematics}, 83:155--179, 1990.

\bibitem{bronnimann:kettner:pocchiola:snoeyink:CountingAndEnumeratingPointedPs%
eudoTriang:2006}
H.~Br{\"{o}}nnimann, L.~Kettner, M.~Pocchiola, and J.~Snoeyink.
\newblock Counting and enumerating pointed pseudotriangulations with the greedy
  flip algorithm.
\newblock {\em SIAM Journal on Computing}, 36:721--739, 2006.

\bibitem{chazelle:etAl:rayShootingGeodesicTriangulations:1994}
B.~Chazelle, H.~Edelsbrunner, M.~Grigni, L.~J. Guibas, J.~Hershberger,
  M.~Sharir, and J.~Snoeyink.
\newblock Ray shooting in polygons using geodesic triangulations.
\newblock {\em Algorithmica}, 12:54--68, 1994.
\newblock Preliminary version in Automata, Languages and Programming, Lect.
  Notes. Comput. Sci., Vol.~510, Springer-Verlag (1991), pp.~661--673.

\bibitem{vpv-tbmai-03}
{\'E}.~{Colin de Verdi{\`e}re}, M.~Pocchiola, and G.~Vegter.
\newblock Tutte's barycenter method applied to isotopies.
\newblock {\em Computational Geometry: Theory and Applications}, 26(1):81--97,
  2003.

\bibitem{connelly:demaine:rote:carpenterRule:2003}
R.~Connelly, E.~Demaine, and G.~Rote.
\newblock Straightening polygonal arcs and convexifying polygonal cycles.
\newblock {\em Discrete and Computational Geometry}, 30(2):205--239, August
  2003.

\bibitem{edelsbrunner:shah:incrementalTopologicalFlipping:1996}
H.~Edelsbrunner and N.~R. Shah.
\newblock Incremental topological flipping works for regular triangulations.
\newblock {\em Algorithmica}, 15:223--241, 1996.

\bibitem{gelfand:kapranov:zelevinsky:discriminantsResultants:1994}
I.~M. Gel'fand, M.~M. Kapranov, and A.~V. Zelevinsky.
\newblock {\em Discriminants, Resultants and Multidimensional Determinants}.
\newblock Birkh\"auser, Boston, 1994.

\bibitem{goodrich:tamassia:dynamicRayShooting:1997}
M.~T. Goodrich and R.~Tamassia.
\newblock Dynamic ray shooting and shortest paths in planar subdivisions via
  balanced geodesic triangulations.
\newblock {\em Journal of Algorithms}, 23:51--73, 1997.
\newblock Preliminary version in Proc. 9th Annu. Sympos. Comput. Geom. (1993),
  pp.~318--327.

\bibitem{graver:servatius:servatius:CombinatorialRigidity:1993}
J.~Graver, B.~Servatius, and H.~Servatius.
\newblock {\em Combinatorial Rigidity}.
\newblock Graduate Studies in Mathematics vol. 2. American Mathematical
  Society, 1993.

\bibitem{grunbaum:polytopes}
B.~Gr{\"u}nbaum.
\newblock {\em Convex Polytopes}, volume 221 of {\em Graduate Texts in
  Mathematics}.
\newblock Springer-Verlag, 2003.
\newblock 2nd ed.

\bibitem{streinu:haas:etAl:planarMinRigidPseudoTriang-confAndJour:2005}
R.~Haas, D.~Orden, G.~Rote, F.~Santos, B.~Servatius, H.~Servatius, D.~Souvaine,
  I.~Streinu, and W.~Whiteley.
\newblock Planar minimally rigid graphs and pseudo-triangulations.
\newblock {\em Computational Geometry: Theory and Applications}, 31:31--61,
  2005.

\bibitem{henneberg:graphischeStatik:1911-68}
L.~Henneberg.
\newblock {\em Die graphische {Statik} der starren {Systeme}}.
\newblock Leipzig, 1911.
\newblock Johnson Reprint 1968.

\bibitem{speckmann:kettner:etAl:tightDegreeBounds:2003}
L.~Kettner, D.~Kirkpatrick, A.~Mantler, J.~Snoeyink, B.~Speckmann, and
  F.~Takeuchi.
\newblock Tight degree bounds for pseudo-triangulations of points.
\newblock {\em Computational Geometry: Theory and Applications}, 25:3--12,
  2003.

\bibitem{speckmann:kirkpatrick:snoeyink:kineticCollisionDetection:2002}
D.~Kirkpatrick, J.~Snoeyink, and B.~Speckmann.
\newblock Kinetic collision detection for simple polygons.
\newblock {\em International Journal of Computational Geometry and
  Applications}, 12:3--27, 2002.

\bibitem{speckmann:kirkpatrick:kineticMaintenance:2002}
D.~Kirkpatrick and B.~Speckmann.
\newblock Kinetic maintenance of context-sensitive hierarchical representations
  for disjoint simple polygons.
\newblock In {\em Proc. 18th Ann. Symposium on Computational Geometry (SoCG)},
  pages 179--188, 2002.

\bibitem{laman:Rigidity:1970}
G.~Laman.
\newblock On graphs and rigidity of plane skeletal structures.
\newblock {\em Journal of Engineering Mathematics}, 4:331--340, 1970.

\bibitem{lee:associahedron:1989}
C.~Lee.
\newblock The associahedron and triangulations of the $n$-gon.
\newblock {\em European Journal of Combinatorics}, 10:551--560, 1989.

\bibitem{maxwell:Reciprocal:1864}
J.~C. Maxwell.
\newblock On reciprocal figures and diagrams of forces.
\newblock {\em Philosphical Magazine}, 27:250--261, 1864.

\bibitem{maxwell:Reciprocal:1870}
J.~C. Maxwell.
\newblock On reciprocal figures, frames and diagrams of forces.
\newblock {\em Transactions of the Royal Society Edinburgh}, 26:1--40, 1870.

\bibitem{maxwell:Bow:1876}
J.~C. Maxwell.
\newblock On {B}ow's method of drawing diagrams in graphical statics, with
  illustrations from {P}eaucellier's linkage.
\newblock {\em Cambridge Phil. Soc. Proc.}, 2:407--414, 1876.

\bibitem{mccabe:seidel:newLowerBoundsTriangulations:2002}
P.~McCabe and R.~Seidel.
\newblock New lower bounds for the number of straight-edge triangulations of a
  planar point set.
\newblock In {\em Abstracts of the 20th European Workshop on Computational
  Geometry}, pages 175--176, Seville, Spain, Mar. 2004.

\bibitem{orden:roteEtAl:nonCrossingFrameworks:2004}
D.~Orden, G.~Rote, F.~Santos, B.~Servatius, H.~Servatius, and W.~Whiteley.
\newblock Non-crossing frameworks with non-crossing reciprocals.
\newblock {\em Discrete and Computational Geometry}, 32:567--600, 2004.

\bibitem{orden:santos:polytope:2005}
D.~Orden and F.~Santos.
\newblock The polytope of non-crossing graphs on a planar point set.
\newblock {\em Discrete and Computational Geometry}, 33:275--305, 2005.

\bibitem{orden:santos:servatiusB:servatiusH:combinatorialPT:2006}
D.~Orden, F.~Santos, B.~Servatius, and H.~Servatius.
\newblock Combinatorial pseudo-triangulations.
\newblock {\em Discrete Mathematics}, 307:554--566, 2007.

\bibitem{p-bcah-05}
G.~Panina.
\newblock New counterexamples to {A.~D.~Alexandrov's} hypothesis.
\newblock {\em Adv. Geom.}, 5(2):301--317, 2005.

\bibitem{panina:2006}
G.~Panina.
\newblock Planar pseudo-triangulations, spherical pseudo-tilings and hyperbolic
  virtual polytopes.
\newblock Technical Report math.MG:0607171, arXiv, 2006.

\bibitem{pocchiola:vegter:orderType:94}
M.~Pocchiola and G.~Vegter.
\newblock Order types and visibility types of configurations of disjoint convex
  plane sets.
\newblock Technical Report 94-4, LIENS, January 1994.

\bibitem{pocchiola:vegter:mtvg:96}
M.~Pocchiola and G.~Vegter.
\newblock Minimal tangent visibility graphs.
\newblock {\em Computational Geometry: Theory and Applications}, 6:303--314,
  1996.

\bibitem{pocchiola:vegter:pt:96}
M.~Pocchiola and G.~Vegter.
\newblock Pseudo-triangulations: theory and applications.
\newblock In {\em Proc. 12th Annu. Sympos. Comput. Geom.}, pages 291--300,
  1996.

\bibitem{pocchiola:vegter:topoSweep:1996}
M.~Pocchiola and G.~Vegter.
\newblock Topologically sweeping visibility complexes via
  pseudo-triangulations.
\newblock {\em Discrete \& Computational Geometry}, 16(4):419--453, 1996.

\bibitem{pocchiola:vegter:visibilityComplex:1996}
M.~Pocchiola and G.~Vegter.
\newblock The visibility complex.
\newblock {\em International Journal of Computational Geometry and
  Applications}, 6(3):279--308, 1996.
\newblock Preliminary version in Proc. 9th Annu. Sympos. Comput. Geom. (1993),
  pp.~328--337.

\bibitem{pocchiola:vegter:pc:99}
M.~Pocchiola and G.~Vegter.
\newblock On polygonal covers.
\newblock In B.~Chazelle, J.~Goodman, and R.~Pollack, editors, {\em Advances in
  Discrete and Computational Geometry}, volume 223 of {\em Contemporary
  Mathematics}, pages 257--268. AMS, Providence, 1999.

\bibitem{santos:randall:rote:snoeyink:countingTriangulationsWheels:2001}
D.~Randall, G.~Rote, F.~Santos, and J.~Snoeyink.
\newblock Counting triangulations and pseudo-triangulations of wheels.
\newblock In {\em Proc. $13$th Canadian Conf. Computational Geometry,
  Waterloo}, pages 149--152, 2001.

\bibitem{richter-gebert:realizationSpacesPolytopes:1996}
J.~Richter-Gebert.
\newblock {\em Realization Spaces of Polytopes}.
\newblock Springer-Verlag, 1996.

\bibitem{rote:santos:streinu:polytopePseudoTriangulations:2003}
G.~Rote, F.~Santos, and I.~Streinu.
\newblock Expansive motions and the polytope of pointed pseudo-triangulations.
\newblock In B.~Aronov, S.~Basu, J.~Pach, and M.~Sharir, editors, {\em Discrete
  and Computational Geometry -- The Goodman-Pollack Festschrift}, pages
  699--736. Springer-Verlag, Berlin, 2003.

\bibitem{rs-pdpts-05}
G.~Rote and A.~Schulz.
\newblock A pointed {Delaunay} pseudo-triangulation of a simple polygon.
\newblock In {\em Abstracts of the 21st European Workshop on Computational
  Geometry}, pages 77--80, Eindhoven, Mar. 2005.

\bibitem{rote:wang:wang:xu:constrainedPseudoTriangulations:2003}
G.~Rote, C.~A. Wang, L.~Wang, and Y.~Xu.
\newblock On constrained minimum pseudo-triangulations.
\newblock In T.~Warnow and B.~Zhu, editors, {\em {\em Computing and
  Combinatorics}, Proc. $9$th Ann. Intern. Computing and Combinatorics Conf.
  (COCOON 2003)}, volume 2697 of {\em Lecture Notes in Computer Science}, pages
  445--454, Big Sky, USA, July 2003. Springer-Verlag.

\bibitem{santos-icm06}
F.~Santos.
\newblock Geometric bistellar flips: the setting, the context and a
  construction.
\newblock In {\em International Congress of Mathematicians. Vol. III}, pages
  931--962. Eur. Math. Soc., Z\"urich, 2006.

\bibitem{Sharir+Welzl:RandomTriangulations:2006}
M.~Sharir and E.~Welzl.
\newblock Random triangulations of planar point sets.
\newblock In {\em Proc. 22nd Ann. Symp.~Comput.~Geom.~(SoCG)}, pages 273--281.
  ACM, 2006.

\bibitem{ss-crcpt-05}
J.~Snoeyink and I.~Streinu.
\newblock Computing rigid components of pseudo-triangulation mechanisms in
  linear time.
\newblock In {\em Proc. $17$th Canadian Conference on Computational Geometry
  (CCCG 2005)}, pages 223--226, Windsor, Ontario, Canada, 2005.

\bibitem{s-kdscd-01}
B.~Speckmann.
\newblock {\em Kinetic Data Structures for Collision Detection}.
\newblock PhD thesis, Dept.\ of Computer Science, University of British
  Columbia, 2001.

\bibitem{speckmann:toth:vertexGuardsPseudoTriang:2003}
B.~Speckmann and C.~D. T\'{o}th.
\newblock Allocating vertex $\pi$-guards in simple polygons via
  pseudo-triangulations.
\newblock {\em Discrete and Computational Geometry}, 33:345--364, 2005.

\bibitem{streinu:combinatorialRoadmaps:2003}
I.~Streinu.
\newblock Combinatorial roadmaps in configuration spaces of simple planar
  polygons.
\newblock In S.~Basu and L.~Gonz{\'a}lez-Vega, editors, {\em DIMACS Workshop on
  Algorithmic and Quantitative Aspects of Real Algebraic Geometry in
  Mathematics and Computer Science, March 12-16, 2001, DIMACS Center, Rutgers
  University, Piscataway, NJ, USA}, pages 181--206. American Mathematical
  Society, 2001.

\bibitem{streinu:pseudoTriangRigidityMotionPlanning-confAndJour:2005}
I.~Streinu.
\newblock Pseudo-triangulations, rigidity and motion planning.
\newblock {\em Discrete and Computational Geometry}, 34:587--635, December
  2005.
\newblock Preliminary version in Proc. $41$st Annual IEEE Symp. Found. Comput.
  Sci. (FOCS 2000), pp.~443--453.

\bibitem{streinu:parallelPseudoTriangulations:2005}
I.~Streinu.
\newblock Parallel redrawing mechanisms, pseudo-triangulations and kinetic
  planar graphs.
\newblock In P.~Healy and N.~S. Nikolov, editors, {\em Proc. 13th International
  Symposium on Graph Drawing (GD 2005), Limerick, Ireland, September 12--14,
  2005, Revised Papers}, volume 3843 of {\em Lecture Notes in Computer
  Science}, pages 421--433, 2006.

\bibitem{streinu:whiteley:origami:2005}
I.~Streinu and W.~Whiteley.
\newblock Single-vertex origami and spherical expansive motions.
\newblock In J.~Akiyama and M.~Kano, editors, {\em Proc. Japan Conf. Discrete
  and Computational Geometry (JCDCG 2004)}, volume 3742 of {\em Lecture Notes
  in Computer Science}, pages 161--173, Tokai University, Tokyo, 8-11 October
  2004 2005. Springer-Verlag.

\bibitem{t-crg-60}
W.~T. Tutte.
\newblock Convex representations of graphs.
\newblock {\em Proc. London Math. Soc., III Ser.}, 10:304--320, 1960.

\bibitem{t-hdg-63}
W.~T. Tutte.
\newblock How to draw a graph.
\newblock {\em Proc. London Math. Soc., III Ser.}, 13:743--768, 1963.

\bibitem{whiteley:Matroids:1996}
W.~Whiteley.
\newblock Some matroids from discrete applied geometry.
\newblock In J.~O. J.~Bonin and B.~Servatius, editors, {\em Matroid Theory},
  volume 197 of {\em Contemporary Mathematics}, pages 171--311. American
  Mathematical Society, 1996.

\bibitem{whiteley:surveyHandbook:2004}
W.~Whiteley.
\newblock Rigidity and scene analysis.
\newblock In J.~E. Goodman and J.~O'Rourke, editors, {\em Handbook of Discrete
  and Computational Geometry}, chapter~60, pages 1327--1354. CRC Press, Boca
  Raton New York, second edition, 2004.

\bibitem{ziegler:lectures-on-polytopes:1995}
G.~Ziegler.
\newblock {\em Lectures on Polytopes}, volume 152 of {\em Graduate Texts in
  Mathematics}.
\newblock Springer-Verlag, 1995.

\end{thebibliography}

%\input{biblio}

%\clearpage

%\tableofcontents

\end{document}